\documentclass[12pt]{amsart}
\usepackage{amsmath,amssymb,amscd}
\usepackage[dvips]{graphicx}
\usepackage[mathscr]{eucal}

\theoremstyle{plain}
\newtheorem{thm}{Theorem}[section]
\newtheorem{cor}[thm]{Corollary}
\newtheorem{lem}[thm]{Lemma}
\newtheorem{prop}[thm]{Proposition}

\newtheorem{defn}[thm]{Definition}

\newcommand{\CC}{\operatorname{ch}}
\newcommand{\END}{\operatorname{End}}
\newcommand{\TR}{\operatorname{Tr}}
\newcommand{\TTR}{\operatorname{tr}}
\newcommand{\KER}{\operatorname{Ker}}
\newcommand{\IIm}{\operatorname{Im}}
\newcommand{\CCGS}{\CC ^{GS}}
\newcommand{\AAR}{\operatorname{Ar}}
\newcommand{\HOM}{\operatorname{Hom}}
\newcommand{\IId}{\operatorname{id}}
\newcommand{\SK}{\operatorname{sk}}
\newcommand{\TD}{\operatorname{Td}}

\newcommand{\CUB}{\operatorname{Cub}}
\newcommand{\CONE}{\mathit{Cone}}
\newcommand{\SGN}{\operatorname{sgn}}

\newcommand{\ACYC}{\text{$\varphi $-ac}}

\newcommand{\RE}{\operatorname{Re}}
\newcommand{\DEG}{\operatorname{deg}}

\newcommand{\sbinom}[2]{{\textstyle \binom{#1}{#2}}}
\newcommand{\sfrac}[2]{{\textstyle \frac{#1}{#2}}}

\begin{document}

\title[\null ]{Higher arithmetic $K$-theory}
\author[\null ]{Yuichiro Takeda}

\address{Fuculty of Mathematics, Kyushu University 33, 
Fukuoka, 812-8581, Japan}
\email{yutakeda@comp.metro-u.ac.jp}

\maketitle

\vskip 1pc
\begin{center}
{\sc Introduction}
\end{center}
\vskip 1pc

The aim of this paper is to provide a new definition 
of higher $K$-theory in Arakelov geometry and to show that it enjoys 
the same formal properties as the higher algebraic $K$-theory of schemes. 

Let $X$ be a proper arithmetic variety, that is, let $X$ be a regular scheme 
which is flat, proper and of finite type over $\mathbb Z$, 
the ring of integers.
In the research on arithmetic Chern characters of hermitian vector 
bundles on $X$, Gillet and Soul\'{e} defined the arithmetic $K_0$-group 
$\widehat{K}_0(X)$ of $X$ \cite{gilsoul}.
It can be viewed as an analogue in Arakelov geometry of 
the $K_0$-group of vector bundles on a scheme. 

After the advent of $\widehat{K}_0(X)$, its higher extension was 
discussed in some papers, such as \cite{burwang, deligne, soule}.
In these papers one common thing was suggested that higher arithmetic 
$K$-theory should be obtained as the homotopy group of the homotopy 
fiber of the Beilinson's regulator map. 
That is to say, it was predicted that there would exist a group $KM_n(X)$ for 
each $n\geq 0$ satisfying the long exact sequence
$$
\cdots \to K_{n+1}(X)\overset{\rho }{\to }\underset{p}{\oplus }
H^{2p-n-1}_{\mathscr D}(X, \mathbb R(p))\to 
KM_n(X)\to K_n(X)\to \cdots , 
$$
where $H^n_{\mathscr D}(X, \mathbb R(p))$ is the real 
Deligne cohomology and $\rho $ is the Beilinson's regulator map. 

To get the homotopy fiber, a simplicial description of the regulator map 
is necessary.
And it was given by Burgos and Wang in \cite{burwang}. 
For a compact complex manifold $M$, they defined 
an exact cube of hermitian vector bundles on $M$ and associated with it 
a differential form called a {\it higher Bott-Chern form}. 
This leads us to a homomorphism of complexes 
$$
\CC :\mathbb Z\widehat{S}_*(M)\to \mathscr D^*(M, p)[2p+1]
$$
from the complex $\mathbb Z\widehat{S}_*(M)$ associated 
with the S-construction of the category of hermitian vector bundles 
on $M$ to the complex $\mathscr D^*(M, p)$ computing the real 
Deligne cohomology of $M$, which is defined in \cite{burgos1}. 
It is the main theorem of \cite{burwang} that the following map 
coincides with the higher Chern character map with values 
in Deligne cohomology: 
$$
\rho :K_n(M)\simeq \pi _{n+1}(\widehat{S}(M))
\overset{Hurewicz}{\longrightarrow }H_{n+1}(\widehat{S}(M))
\overset{H(\CC )}{\longrightarrow }H^{2p-n}_{\mathscr D}(M, \mathbb R(p)).
$$
Applying the theory of higher Bott-Chern forms to an arithmetic variety, 
we can obtain a simplicial 
description of the regulator map.

In this paper, we will give another definition of higher arithmetic 
$K$-theory of a proper arithmetic variety by means of higher 
Bott-Chern forms, which differs from the one coming from 
the homotpy fiber of the regulator map. 
One of the remarkable features of our arithmetic $K$-theory is that 
it is given as an extension of the usual $K$-theory by 
the cokernel of the regulator map.

Before explaining our method, let us recall the definition of 
$\widehat{K}_0(X)$.
For a proper arithmetic variety $X$, let $\mathscr A^{p,p}(X)$ be 
the space of real $(p,p)$-forms $\omega $ on $X(\mathbb C)$ such that 
$\overline{F}_{\infty }^*\omega =(-1)^p\omega $ for the complex 
conjugation $F_{\infty }:X(\mathbb C)\to X(\mathbb C)$ and 
let $\widetilde{\mathscr A}(X)=\underset{p}{\oplus }\mathscr A^{p,p}(X)
/\IIm \partial +\IIm \overline{\partial }$. 
Then $\widehat{K}_0(X)$ is defined as a factor group of the free abelian 
group generated by pairs $(\overline{E}, \omega )$ of a hermitian vector 
bundle $\overline{E}$ on $X$ and 
$\omega \in \widetilde{\mathscr A}(X)$.
The relation on pairs is given by each short exact sequence 
$\mathscr E:0\to \overline{E^{\prime }}\to \overline{E}\to 
\overline{E^{\prime \prime }}\to 0$ as follows:
$$
(\overline{E^{\prime }}, \omega ^{\prime })+
(\overline{E^{\prime \prime }}, \omega ^{\prime \prime })=
(\overline{E}, \omega ^{\prime }+\omega ^{\prime \prime }+
\widetilde{\CC }(\mathscr E)),
$$
where $\widetilde{\CC }(\mathscr E)$ is the Bott-Chern secondary 
characteristic class of $\mathscr E$ and $\omega ^{\prime }, 
\omega ^{\prime \prime }\in \widetilde{\mathscr A}(X)$. 

We can interpret the above definition of 
$\widehat{K}_0(X)$ in terms of loops and homotopies on 
$\vert \widehat{S}(X)\vert $, the topological realization of 
the $S$-construction of the category of hermitian vector bundles 
on $X$.
Let us consider a pair $(l, \omega )$, where $l$ is a pointed simplicial 
loop on $\vert \widehat{S}(X)\vert $ and 
$\omega \in \widetilde{\mathscr A}(X)$. 
Two pairs $(l, \omega )$ and $(l^{\prime }, \omega ^{\prime })$ are 
defined to be homotopy equivalent if there is a cellular homotopy 
$H:S^1\times I/\{*\}\times I\to \vert \widehat{S}(X)\vert $ from $l$ to 
$l^{\prime }$ such that the Bott-Chern secondary characteristic class 
$\widetilde{\CC }(H)$ of $H$, which is defined in a natural way, 
is equal to the difference $\omega ^{\prime }-\omega $.
Let $\widehat{\pi }_1(\vert \widehat{S}(X)\vert , \widetilde{\CC })$ 
denote the set of all equivalence classes of such pairs.
Then it carries the structure of an abelian group and 
the map 
$$
\widehat{K}_0(X)\to 
\widehat{\pi }_1(\vert \widehat{S}(X)\vert , \widetilde{\CC })
$$
defined by $(\overline{E}, \omega )\mapsto (l_{\overline{E}}, -\omega )$, 
where $l_{\overline{E}}$ is the simplicial loop on 
$\vert \widehat{S}(X)\vert $ determined by $\overline{E}$, is shown 
to be bijective.

In this paper, we will generalize the above construction of 
$\widehat{\pi }_1$ to higher homotopy groups to define higher arithmetic 
$K$-groups.
Strictly speaking, we will define the $n$-th arithmetic $K$-group 
$\widehat{K}_n(X)$ as the set of all homotopy equivalence classes of 
pairs $(f, \omega )$, where $f:S^{n+1}\to \vert \widehat{S}(X)\vert $ 
is a pointed cellular map and $\omega $ is a real differential form on 
$X(\mathbb C)$ modulo exact forms.
We will employ the theory of higher Bott-Chern forms in defining 
a homotopy equivalence relation on these pairs.

Let us describe the content of the paper in more detail.

In \S 1 we introduce some materials used in the present paper, such as 
$S$-construction, cubes and higher Bott-Chern forms.
Furthermore, by renormalizing the higher Bott-Chern forms, we obtain 
a homomorphism of complexes
$$
\CCGS :\mathbb Z\widehat{S}_*(X)\to \mathscr A_*(X)[1],
$$
where $\mathscr A_*(X)$ is a complex of vector spaces of 
real differential forms on $X(\mathbb C)$ satisfying a certain Hodge 
theoretic condition.
In \S 2 we propose the notion of modified homotopy groups, which 
is a higher generalization of the above $\widehat{\pi }_1$.
In \S 3 we define higher arithmetic $K$-groups $\widehat{K}_*(X)$ 
as the homotopy groups of $\vert \widehat{S}(X)\vert $ modified by 
the homomorphism $\CCGS $.
We show the following exact sequence concerning $\widehat{K}_n(X)$:
$$
K_{n+1}(X)\to \widetilde{\mathscr A}_{n+1}(X)
\to \widehat{K}_n(X)\to K_n(X)\to 0,
$$
where $\widetilde{\mathscr A}_{n+1}(X)=
\mathscr A_{n+1}(X)/\IIm d_{\mathscr A}$. 
When $n=0$, this exact sequence has already been obtained 
in \cite{gilsoul}.
Moreover, we define the Chern form map 
$$
\CCGS _n:\widehat{K}_n(X)\to \mathscr A_n(X) 
$$
in the same way as $\CC :\widehat{K}_0(X)\to \mathscr A(X)$ 
in \cite{gilsoul}.
The group $KM_n(X)$, which is characterized as the long exact 
sequence as mentioned before, is realized as the kernel of $\CCGS _n$.
When we fix an $F_{\infty }$-invariant K\"{a}hler metric $h_X$ 
on $X(\mathbb C)$, we define Arakelov $K$-group $K_n(\overline{X})$ 
of $\overline{X}=(X, h_X)$ as the subgroup of $\widehat{K}_n(X)$ 
consisting of all elements $x\in \widehat{K}_n(X)$ such that 
$\CCGS _n(x)$ is harmonic with respect to $h_X$.
The group $K_n(\overline{X})$ fits into the exact sequence 
$$
K_{n+1}(X)\overset{\rho }{\to }\underset{p}{\oplus }
H^{2p-n-1}_{\mathscr D}(X, \mathbb R(p))\to 
K_n(\overline{X})\to K_n(X)\to 0, 
$$
where $\rho $ is the Beilinson's regulator map.

The next two sections concern product structure on $\widehat{K}_*(X)$.
In \S 4 we prove a product formula for higher Bott-Chern forms. 
It provides an alternative proof of the fact that the regulator map 
respects the products. 
In \S 5, we define a product in higher arithmetic $K$-theory.
We show that this product is graded commutative up to $2$-torsion.
A striking property of the product is the lack of the associativity.
In other words, for $x, y, z\in \widehat{K}_*(X)$, $(x\times y)\times z$ 
is not equal to $x\times (y\times z)$ in general.
We compute this difference explicitly. 
Moreover, we show that there is an associative product in the Arakelov 
$K$-theory.

In \S 6, we define a direct image morphism in higher arithmetic 
$K$-theory.
To do this we employ a higher analytic torsion form of an exact metrized 
cube defined by Roessler \cite{roessler}.
Moreover we establish the projection formula.

\vspace{5mm}

\setcounter{tocdepth}{1}
\tableofcontents

\newpage

\vskip 2pc
\section{Preliminaries}
\vskip 1pc

\subsection{Conventions on complexes} \ \ 
Let us first settle some conventions on complexes.
For a complex $A^*=(A^n, d_A)$ of an abelian category 
$\mathfrak A$ and $n\in \mathbb Z$, the $n$-th translation $A[n]^*$ 
is defined as $A[n]^k=A^{n+k}$ and $d_{A[n]}=(-1)^nd_A$.
For a morphism of complexes $u:A^*\to B^*$, the {\it mapping cone} 
$Cone(u)$ is defined by 
$$
Cone(u)^k=A^{k+1}\oplus B^k
$$
and the differential $d:A^{k+1}\oplus B^k\to A^{k+2}\oplus B^{k+1}$ 
is defined by 
$$
d(a,b)=(-d_A(a), u(a)+d_B(b)). 
$$

A {\it homological complex} is a family $\{A_n\}_{n\in \mathbb Z}$ 
of objects of $\mathfrak A$ with morphisms \linebreak 
$\partial _A:A_n\to A_{n-1}$ such that $\partial _A^2=0$.
For a homological complex $(A_n, \partial _A)$, we can define a complex 
$A^*$ by $A^n=A_{-n}$ and $d_A=\partial _A$.
The $n$-th translation of a homological complex $A_*$ is defined by 
$A[n]_k=A_{k-n}$ and $\partial _{A[n]}=(-1)^n\partial _A$.

\vskip 1pc
\subsection{$S$-construction} \ \ 
In this subsection we recall $S$-construction developed by 
Waldhausen \cite{wal}.
Let $[n]$ be the finite ordered set $\{0, 1, \cdots , n\}$ and 
$\AAR [n]$ the category of arrows of $[n]$.
For a small exact category $\mathfrak A$, let $S_n\mathfrak A$ 
be the set of functors 
$$
E:\AAR [n]\to \mathfrak A, \ i\leq j\mapsto E_{i,j}
$$
satisfying the following conditions: 
\begin{enumerate}
\item  \ $E_{i, i}=0$ for any $0\leq i\leq n$. 
\item  \ For any $i\leq j\leq k$, $E_{i,j}\to E_{i,k}
\to E_{j,k}$ is a short exact sequence of $\mathfrak A$.
\end{enumerate}
For example, $S_0\mathfrak A=\{0\}$, $S_1\mathfrak A$ is the 
set of objects of $\mathfrak A$ and $S_2\mathfrak A$ is the set 
of short exact sequences of $\mathfrak A$.
The functor 
$$
S\mathfrak A:[n]\mapsto S_n\mathfrak A
$$
becomes a simplicial set with the base point $0\in S_0\mathfrak A$. 

The set $S_n\mathfrak A$ can be identified with the set of sequences 
of injections 
$$
E_{0,1}\rightarrowtail \cdots \rightarrowtail E_{0,n}
$$
with quotients $E_{i,j}\simeq E_{0,j}/E_{0,i}$ for each $i<j$.
By using this identification, we can describe the boundary maps and 
the degeneracy maps of $S\mathfrak A$ as follows:
$$
\partial _k(E_{0,1}\rightarrowtail \cdots \rightarrowtail E_{0,n})=
\begin{cases}
E_{1,2}\rightarrowtail \cdots \rightarrowtail E_{1,n}, &k=0, \\
E_{0,1}\rightarrowtail \cdots \rightarrowtail E_{0,k-1}
\rightarrowtail E_{0,k+1}\rightarrowtail \cdots 
\rightarrowtail E_{0,n}, &k\geq 1   \end{cases}
$$
and 
$$
s_k(E_{0,1}\rightarrowtail \cdots \rightarrowtail E_{0,n})=
\begin{cases}
0\rightarrowtail E_{0,1}\rightarrowtail \cdots 
\rightarrowtail E_{0,n}, &k=0,   \\
E_{0,1}\rightarrowtail \cdots \rightarrowtail E_{0,k}
\overset{id}{\rightarrowtail}E_{0,k}\rightarrowtail \cdots 
\rightarrowtail E_{0,n}, &k\geq 1.  \end{cases}
$$

\vskip 1pc
\begin{thm}\cite[\S 1.9]{wal} \ \ 
There is a homotopy equivalence map between the topological 
realizations of simplicial sets
$$
\vert S\mathfrak A\vert \simeq \vert BQ\mathfrak A\vert , 
$$
where $BQ\mathfrak A$ is the classifying space of the Quillen's 
$Q$-construction of $\mathfrak A$ \cite{quillen}.
Therefore for $n\geq 0$, it follows that 
$$ 
\pi _{n+1}(\vert S\mathfrak A\vert , 0)\simeq 
K_n(\mathfrak A).
$$
\end{thm}

\vskip 1pc
\subsection{Exact $n$-cubes} \ \ 
Let us recall the notion of an exact $n$-cube.
For more details, see \cite{burgos2, burwang}.
Let $<-1, 0, 1>$ be the ordered set consisting of three elements 
and $<-1, 0, 1>^n$ its $n$-th power.
For a small exact category $\mathfrak A$, a functor 
$\mathcal F:<-1, 0, 1>^n\to \mathfrak A$ is called an 
$n$-{\it cube} of $\mathfrak A$.
Let $\mathcal F_{\alpha _1, \cdots , \alpha _n}$ denote the image of 
an object $(\alpha _1, \cdots , \alpha _n)$ of $<-1, 0, 1>^n$.
For integers $i$ and $j$ satisfying $1\leq i\leq n$ and $-1\leq j\leq 1$, 
an $(n-1)$-cube $\partial _i^j\mathcal F$ is defined by 
$(\partial _i^j\mathcal F)_{\alpha _1, \cdots , \alpha _{n-1}}=
\mathcal F_{\alpha _1, \cdots , \alpha _{i-1}, j, \alpha _i, 
\cdots , \alpha _{n-1}}$. 
It is called a {\it face} of $\mathcal F$. 
For an object $\alpha $ of $<-1, 0, 1>^{n-1}$ and an integer $i$ 
satisfying $1\leq i\leq n$, a $1$-cube 
$\partial _{i^c}^{\alpha }\mathcal F$ called an {\it edge} of 
$\mathcal F$ is defined by 
$$
\mathcal F_{\alpha _1, \cdots , \alpha _{i-1},-1,\alpha _i, 
\cdots , \alpha _{n-1}}\to \mathcal F_{\alpha _1, \cdots , 
\alpha _{i-1},0, \alpha _i, \cdots , \alpha _{n-1}}\to 
\mathcal F_{\alpha _1, \cdots , 
\alpha _{i-1},1,\alpha _i, \cdots , \alpha _{n-1}}.
$$
An $n$-cube $\mathcal F$ is said to be {\it exact} if all edges of 
$\mathcal F$ are short exact sequences.
Let $C_n\mathfrak A$ denote the set of all exact $n$-cubes of 
$\mathfrak A$.
If $\mathcal F$ is an exact $n$-cube, then any face 
$\partial _i^j\mathcal F$ is also exact.
Hence $\partial _i^j$ induces a map
$$
\partial _i^j:C_n\mathfrak A\to C_{n-1}\mathfrak A.
$$

Let $\mathcal F$ be an exact $n$-cube of $\mathfrak A$.
For an integer $i$ satisfying $1\leq i\leq n+1$, 
let $s_i^1\mathcal F$ be an exact $(n+1)$-cube defined as follows: 
$$
(s_i^1\mathcal F)_{\alpha _1, \cdots , \alpha _{n+1}}=\begin{cases}
0, &\alpha _i=1,  \\
\mathcal F_{\alpha _1, \cdots , \alpha _{i-1}, \alpha _{i+1}, 
\cdots , \alpha _{n+1}}, \ &\alpha _i\not= 1 
\end{cases}
$$
and the morphism $(s_i^1\mathcal F)_{\alpha _1, \cdots , \alpha _{i-1}, -1, 
\alpha _{i+1}, \cdots , \alpha _{n+1}}\to (s_i^1\mathcal F)_{\alpha _1, 
\cdots , \alpha _{i-1}, 0, \alpha _{i+1}, \cdots , \alpha _{n+1}}$ 
is the identity of $\mathcal F_{\alpha _1, \cdots , 
\alpha _{i-1}, \alpha _{i+1}, \cdots , \alpha _{n+1}}$. 
In addition, let $s_i^{-1}\mathcal F$ be an exact $(n+1)$-cube defined 
as follows: 
$$
(s_i^{-1}\mathcal F)_{\alpha _1, \cdots , \alpha _{n+1}}=\begin{cases}
0, &\alpha _i=-1,  \\
\mathcal F_{\alpha _1, \cdots , \alpha _{i-1}, \alpha _{i+1}, 
\cdots , \alpha _{n+1}}, \ &\alpha _i\not= -1 
\end{cases}
$$
and the morphism $(s_i^{-1}\mathcal F)_{\alpha _1, \cdots , \alpha _{i-1}, 0, 
\alpha _{i+1}, \cdots , \alpha _{n+1}}\to (s_i^{-1}\mathcal F)_{\alpha _1, 
\cdots , \alpha _{i-1}, 1, \alpha _{i+1}, \cdots , \alpha _{n+1}}$ 
is the identity of $\mathcal F_{\alpha _1, \cdots , 
\alpha _{i-1}, \alpha _{i+1}, \cdots , \alpha _{n+1}}$. 
Then the map 
$$
s_i^j:C_n\mathfrak A\to C_{n+1}\mathfrak A
$$
is also defined.
An exact cube written as $s_i^j\mathcal F$ is said to be 
{\it degenerate}.

Let $\mathbb ZC_n\mathfrak A$ be the free abelian group generated by 
$C_n\mathfrak A$.
Let $D_n\subset \mathbb ZC_n\mathfrak A$ be the subgroup generated by 
all degenerate exact $n$-cubes. 
Let $\widetilde{\mathbb Z}C_n\mathfrak A=\mathbb ZC_n\mathfrak A/D_n$ 
and 
$$
\partial =\underset{i=1}{\overset{n}{\sum }}
\underset{j=-1}{\overset{1}{\sum }}(-1)^{i+j+1}\partial _i^j: 
\widetilde{\mathbb Z}C_n\mathfrak A\to 
\widetilde{\mathbb Z}C_{n-1}\mathfrak A.
$$
Then $\widetilde{\mathbb Z}C_*\mathfrak A=
(\widetilde{\mathbb Z}C_n\mathfrak A, \partial )$ becomes a 
homological complex. 

We can construct an exact $(n-1)$-cube $\CUB (E)$ associated with 
an element $E\in S_n\mathfrak A$ in an inductive way.
In the case of $n=1$, let us define $\CUB (E)=E$ for 
$E\in S_1\mathfrak A$. 
When an exact $(m-1)$-cube is associated with any element of 
$S_m\mathfrak A$ for $m<n$, an $(n-1)$-cube $\CUB (E)$ associated 
with $E\in S_n\mathfrak A$ is defined as 
\begin{gather*}
\partial _1^{-1}\CUB (E)=s_{n-2}^1\cdots s_1^1(E_{0,1}), \\
\partial _1^0\CUB (E)=\CUB (\partial _1E), \\
\partial _1^1\CUB (E)=\CUB (\partial _0E).
\end{gather*}
Then $\CUB :S_n\mathfrak A\to C_{n-1}\mathfrak A$ 
induces a homomorphism of complexes
$$
\CUB :\mathbb ZS_*\mathfrak A[1]\to \widetilde{\mathbb Z}C_*\mathfrak A.
$$

\vskip 1pc
{\it Remark}: \ 
The differential of $\widetilde{\mathbb Z}C_*\mathfrak A$ 
defined above is the minus of the one defined in \cite{burwang}.
But since the differential of $\mathbb ZS_*\mathfrak A[1]$ is also 
the minus of the one of $\mathbb ZS_*\mathfrak A$, the homomorphism 
$\CUB $ is compatible with the differentials of complexes.

\vskip 1pc
\subsection{Deligne cohomology} \ \ 
Let $M$ be a complex algebraic manifold of dimension $n$, that is, 
let $M$ be the analytic space consisting of all $\mathbb C$-valued points 
of a smooth algebraic variety of dimension $n$ over $\mathbb C$.
Then higher Chern character map from higher algebraic 
$K$-theory to a reasonable cohomology theory 
$$
\CC _i:K_i(M)\to \underset{j}{\oplus }H^{2j-i}(M, \Gamma (j))
$$
has been developed in \cite{gillet}.
In the case of $i=0$, it is given by the Chern characters of 
vector bundles.
But when $i>0$ and $M$ is compact, $\CC _i$ with values in the singular 
cohomology is known to be trivial.

In order to obtain a nontrivial higher Chern character map, we should 
consider Deligne cohomology. 
From now on we assume that $M$ is compact.
Let $\Omega ^*_M$ be the complex of analytic sheaves of holomorphic 
differential forms and let $F^p\Omega ^*_M$ be the de Rham 
filtration of $\Omega ^*_M$, that is, 
$$
F^p\Omega ^*_M=(0\to \Omega ^p_M\to \Omega ^{p+1}_M\to 
\cdots \to \Omega ^n_M).
$$
The {\it real Deligne complex} $\mathbb R(p)_{\mathscr D}$ on $M$ is 
a complex of sheaves defined as follows:
$$
\CONE (\mathbb R(p)\oplus F^p\Omega ^*_M\overset{\varepsilon -\iota }
{\longrightarrow }\Omega ^*_M)[-1],
$$
where $\varepsilon $ and $\iota $ are natural inclusions. 
It is obvious that $\mathbb R(p)_{\mathscr D}$ is quasi-isomorphic 
to the complex 
$$
0\to \mathbb R(p)\to \mathcal O_M\to \Omega ^1_M\to \cdots \to 
\Omega ^{p-1}_M\to 0,
$$
where $\mathbb R(p)$ is of degree zero. 
The {\it real Deligne cohomology} $H_{\mathscr D}^i(M, \mathbb R(p))$ 
is defined as the hypercohomology of $\mathbb R(p)$, that is, 
$$
H^i_{\mathscr D}(M, \mathbb R(p))=\mathbf H^i(M, \mathbb R(p)_{\mathscr D}).
$$

In fact, Deligne cohomology can be defined when $M$ is neither smooth 
nor compact \cite{esnault}.
Moreover, it can be a target of 
a nontrivial higher Chern character map, that is, there is 
a homomorphism 
$$
\CC _i:K_i(M)\to \underset{p}{\oplus }H^{2p-i}_{\mathscr D}
(M, \mathbb R(p)),
$$
which is far from trivial.

Let $\mathscr E^p_{\mathbb R}(M)$ be the space of real smooth 
differential forms of degree $p$ on $M$ and 
$\mathscr E^p(M)=\mathscr E^p_{\mathbb R}(M)
\underset{\mathbb R}{\otimes }\mathbb C$. 
Let $\mathscr E^{p,q}(M)$ be the space of complex differential 
forms of type $(p, q)$ on $M$.
We set 
$$
\mathscr D^n(M, p)=\begin{cases}
\mathscr E^{n-1}_{\mathbb R}(M)(p-1)\cap 
\underset{p^{\prime }<p, q^{\prime }<p}
{\underset{p^{\prime }+q^{\prime }=n-1}{\oplus }}
\mathscr E^{p^{\prime }, q^{\prime }}(M),  &  n<2p,  \\
\mathscr E^{2p}_{\mathbb R}(M)(p)\cap 
\mathscr E^{p, p}(M)\cap \KER d,  &  n=2p, \\
0, & n>2p
\end{cases}
$$
and define a differential 
$d_{\mathscr D}:\mathscr D^n(M, p)\to \mathscr D^{n+1}(M, p)$ by 
$$
d_{\mathscr D}(\omega )=\begin{cases}
-\pi (d\omega ),  &  n<2p-1,  \\
-2\partial \overline{\partial }\omega , & n=2p-1,  \\
0, & n>2p-1,  \end{cases}
$$
where $\pi :\mathscr E^n(M)\to \mathscr D^n(M, p)$ is the canonical 
projection.

\vskip 1pc
\begin{thm}\cite[Thm.2.6]{burgos1} \ \ 
Let $M$ be a compact complex algebraic manifold.
Then \linebreak 
$(\mathscr D^*(M, p), d_{\mathscr D})$ becomes a complex of 
$\mathbb R$-vector spaces and there is a canonical isomorphism 
$$
H^n(\mathscr D^*(M, p), d_{\mathscr D})\simeq 
H^n_{\mathscr D}(M, \mathbb R(p))
$$
if $n\leq 2p$.
\end{thm}

\vskip 1pc
\subsection{Higher Bott-Chern forms} \ \ 
In this subsection we recall the higher Bott-Chern forms 
developed by Burgos and Wang.
For more details, the reader should consult the original 
paper \cite{burwang} or a survey \cite{burgos2}.
A {\it hermitian vector bundle} $\overline{E}=(E, h)$ on a complex 
algebraic manifold $M$ is an algebraic vector bundle $E$ on $M$ 
with a smooth hermitian metric $h$.
On a hermitian vector bundle $\overline{E}$, there is a unique 
connection $\nabla _{\overline{E}}$ that is compatible with 
both the metric and the complex structure.
Let $K_{\overline{E}}$ denote the curvature form of 
$\nabla _{\overline{E}}$.
The Chern form of $\overline{E}$ is defined as 
$$
\CC _0(\overline{E})=\TR (\exp (-K_{\overline{E}}))\in 
\underset{p}{\oplus }\mathscr D^{2p}(M, p).
$$

Hereafter we assume that $M$ is compact. 
An {\it exact metrized $n$-cube} on $M$ is an exact $n$-cube made 
of hermitian vector bundles on $M$.
Let $\mathcal F=\{\overline{E}_{\alpha }\}$ be an exact metrized 
$n$-cube on $M$.
For an $n$-tuple $\alpha =(\alpha _1, \cdots , \alpha _n)$ with 
$-1\leq \alpha _k\leq 1$ and an integer $i$ satisfying $\alpha _i=1$, 
there is a surjection 
$\overline{E}_{\alpha _1, \cdots , \alpha _{i-1}, 0, \alpha _{i+1}, 
\cdots , \alpha _n}\to \overline{E}_{\alpha }$.
$\mathcal F$ is called an {\it emi-$n$-cube} if 
the metric on any $\overline{E}_{\alpha }$ with $\alpha _i=1$ 
coincides with the metric induced from 
$\overline{E}_{\alpha _1, \cdots , \alpha _{i-1}, 0, \alpha _{i+1}, 
\cdots , \alpha _n}$ by the above surjection.

Let $(x: y)$ be the homogeneous coordinate of $\mathbb P^1$.
Let $\overline{\mathcal O(1)}$ be the tautological line bundle on 
$\mathbb P^1$ with the Fubini-Study metric.
For a hermitian vector bundle $\overline{E}$ on $M$, let 
$\overline{E}(1)$ be the hermitian vector bundle on 
$M\times \mathbb P^1$ given by $\pi _1^*\overline{E}\otimes 
\pi _2^*\overline{\mathcal O(1)}$.
Let $\sigma _x, \sigma _y\in H^0(\mathbb P^1, \mathcal O(1))$ 
be global sections of $\mathcal O(1)$ determined by 
$x$ and $y$ respectively.

For an emi-$1$-cube $\mathscr E: \overline{E}_{-1}\to \overline{E}_0\to 
\overline{E}_1$, a map 
$\psi :\overline{E}_{-1}\to \overline{E}_0(1)\oplus \overline{E}_{-1}(1)$ 
on $M\times \mathbb P^1$ is defined by 
$e\mapsto (\iota (e)\otimes \sigma _x, e\otimes \sigma _y)$, where 
$\iota $ is the injection $\overline{E}_{-1}\to \overline{E}_0$.
The {\it $1$-transgression bundle} $\TTR _1(\mathscr E)$ of 
$\mathscr E$ is a hermitian vector bundle on $M\times \mathbb P^1$ 
given by the cokernel of $\psi $ with the induced metric.
It follows from the definition that 
$\TTR _1(\mathscr E)\vert _{M\times \{x=0\}}\simeq \overline{E}_0$ and 
\linebreak 
$\TTR _1(\mathscr E)\vert _{M\times \{y=0\}}\simeq \overline{E}_{-1}
\oplus \overline{E}_1$.
Hence $\TTR _1(\mathscr E)$ can be viewed as a family of hermitian 
vector bundles on $M$ parametrized by $\mathbb P^1$ connecting 
$\overline{E}_0$ with $\overline{E}_{-1}\oplus \overline{E}_1$.

The $n$-transgression bundle of an emi-$n$-cube is defined by 
iterating the above process.
The $1$-transgression $\TTR _1(\mathcal F)$ of an emi-$n$-cube 
$\mathcal F$ is an emi-$(n-1)$-cube on $M\times \mathbb P^1$ defined as 
$$
\TTR _1(\mathcal F)_{\alpha }=\TTR _1(\partial _{n^c}^{\alpha }
(\mathcal F))
$$
for $\alpha \in <-1, 0, 1>^{n-1}$.
The {\it $n$-transgression bundle} of $\mathcal F$ is defined as 
$$
\TTR _n(\mathcal F)=
\overbrace{\TTR _1\TTR _1\cdots \TTR _1}^{\text{$n$ times}}
(\mathcal F),
$$
which is a hermitian vector bundle on $M\times (\mathbb P^1)^n$.

Let $z=x/y$ be the Euclidean coordinate of $\mathbb P^1$ and $z_i$ 
the $i$-th Euclidean coordinate of $(\mathbb P^1)^n$.
For an integer $i$ satisfying $1\leq i\leq n$, a differential form 
with logarithmic poles $S_n^i$ on $(\mathbb P^1)^n$ is defined as 
$$
S_n^i=\underset{\sigma \in \mathfrak S_n}{\sum }(-1)^{\sigma }
\log |z_{\sigma (1)}|^2\frac{dz_{\sigma (2)}}{z_{\sigma (2)}}\wedge 
\cdots \wedge \frac{dz_{\sigma (i)}}{z_{\sigma (i)}}\wedge 
\frac{d\bar{z}_{\sigma (i+1)}}{\bar{z}_{\sigma (i+1)}}\wedge \cdots 
\wedge \frac{d\bar{z}_{\sigma (n)}}{\bar{z}_{\sigma (n)}}.
$$
Let us define the Bott-Chern form of an emi-$n$-cube $\mathcal F$ as 
$$
\CC _n(\mathcal F)=\frac{1}{(2\pi \sqrt{-1})^n}\int_{(\mathbb P^1)^n}
\CC _0(\TTR _n(\mathcal F))\wedge T_n\in \underset{p}{\oplus }
\mathscr D^{2p-n}(M, p),
$$
where
$$
T_n=\frac{(-1)^n}{2n!}\underset{i=1}{\overset{n}{\sum }}(-1)^iS_n^i.
$$

\vskip 1pc
\begin{prop} \ \ 
Let $\mathcal F$ be an emi-$n$-cube on $M$.
Then 
$$
\CC _{n+1}(s^j_i\mathcal F)=0
$$
for $1\leq i\leq n+1$ and $j=\pm 1$.
\end{prop}

{\it Proof}: \ 
Let us first consider the case of $n=0$.
For a hermitian vector bundle $\overline{E}$ on $M$, it follows that 
$s^{-1}_1\overline{E}=\left(0\to \overline{E}\overset{\text{id}}{\to }
\overline{E}\right)$.
Hence the $1$-transgression bundle $\TTR (s^{-1}_1\overline{E})$ 
is isometric to $\overline{E}(1)$.
If $r:\mathbb P^1\to \mathbb P^1$ is an involution given by 
$r^*(z)=z^{-1}$, then $r^*(\overline{E}(1))=\overline{E}(1)$.
Hence we have 
\begin{align*}
\CC _1(s^{-1}_1\overline{E})&=\frac{1}{2\pi \sqrt{-1}}\int_{\mathbb P^1}
\CC _0(\overline{E}(1))\log |z|^2  \\
&=\frac{1}{2\pi \sqrt{-1}}\int_{\mathbb P^1}
r^*\left(\CC _0(\overline{E}(1))\log |z|^2\right)  \\
&=\frac{-1}{2\pi \sqrt{-1}}\int_{\mathbb P^1}
\CC _0(\overline{E}(1))\log |z|^2  \\
&=-\CC _1(s^{-1}_1\overline{E}),
\end{align*}
therefore $\CC _1(s^{-1}_1\overline{E})=0$.
In the same way we can show $\CC _1(s^1_1\overline{E})=0$.

Let us move on to the general case.
For $1\leq i\leq n+1$, let 
$r_i:(\mathbb P^1)^{n+1}\to (\mathbb P^1)^{n+1}$ denote an 
involution given by 
$$
r_i^*z_j=\begin{cases} z_j, \ i\not= j,  \\  z_j^{-1}, \ i=j. \end{cases}
$$
Then we have 
$r_i^*(\TTR _{n+1}(s^j_i\mathcal F))=\TTR _{n+1}(s^j_i\mathcal F)$
and $r_i^*T_{n+1}=-T_{n+1}$. 
Therefore we can show $\CC _{n+1}(s^j_i\mathcal F)=0$ in the same 
way as above.
\qed

\vskip 1pc
Let us extend the definition of the Bott-Chern form to an arbitrary 
exact metrized $n$-cube.
Let $\mathcal F$ be an exact metrized $n$-cube, not necessarily emi.
For an integer $i$ satisfying $1\leq i\leq n$, $\lambda _i^1\mathcal F$ 
is defined as 
$$
(\lambda _i^1\mathcal F)_{\alpha }=
\begin{cases}
\overline{E}_{\alpha }, &\alpha _i=-1 \ \text{or} \ 0,    \\
\overline{E^{\prime }}_{\alpha }, &\alpha _i=1,
\end{cases}
$$
where $\overline{E^{\prime }}_{\alpha }$ is the same vector bundle as 
$\overline{E}_{\alpha }$ with the metric induced from 
$\overline{E}_{\alpha _1, \cdots , \alpha _{i-1},0, \alpha _{i+1}, 
\cdots , \alpha _n}$.
Let $\lambda _i^2\mathcal F$ be an exact metrized $n$-cube defined as 
$$
(\lambda _i^2\mathcal F)_{\alpha }=
\begin{cases}
\overline{E}_{\alpha _1, \cdots , \alpha _{i-1}, 1, \alpha _{i+1}, \cdots , 
\alpha _n},   &\alpha _i=-1,  \\
\overline{E^{\prime }}_{\alpha _1, \cdots , \alpha _{i-1}, 1, 
\alpha _{i+1}, \cdots , \alpha _n}, &\alpha _i=0,  \\
0, &\alpha _i=1
\end{cases}
$$
and the morphism $\overline{E}_{\alpha _1, \cdots , \alpha _{i-1}, 1, 
\alpha _{i+1}, \cdots , \alpha _n}\to \overline{E^{\prime }}_{\alpha _1, 
\cdots , \alpha _{i-1}, 1, \alpha _{i+1}, \cdots , \alpha _n}$ is 
the identity.
We set $\lambda _i\mathcal F=\lambda _i^1\mathcal F\oplus 
\lambda _i^2\mathcal F$.
An emi-$n$-cube $\lambda \mathcal F$ is defined as 
$$
\lambda \mathcal F=\begin{cases}
\lambda _n\lambda _{n-1}\cdots \lambda _1
\mathcal F, &n\geq 1,  \\
\mathcal F, & n=0.
\end{cases}
$$

\vskip 1pc
\begin{defn} \ \ 
The Bott-Chern form of an exact metrized $n$-cube $\mathcal F$ 
is an element of $\underset{p}{\oplus }\mathscr D^{2p-n}(M, p)$ 
defined as 
\begin{align*}
\CC _n(\mathcal F)=\frac{1}{(2\pi \sqrt{-1})^n}\int_{(\mathbb P^1)^n}
\CC _0(\TTR _n(\lambda \mathcal F))\wedge T_n.
\end{align*}
\end{defn}

\vskip 1pc
{\it Remark}: \ 
For an emi-$n$-cube $\mathcal F$, the Bott-Chern form 
$\CC _n(\mathcal F)$ defined in Def.1.4 is the same form 
as the one we have defined before, because there is a decomposition 
$$
\lambda \mathcal F=\mathcal F\oplus \text{(a degenerate emi-$n$-cube)}.
$$

\vskip 1pc
For a compact complex algebraic manifold $M$, let 
$\widehat{\mathscr P}(M)$ be the category of hermitian vector 
bundles on $M$ and $\widehat{S}(M)$ the $S$-construction of 
$\widehat{\mathscr P}(M)$.
We note that $\widehat{S}(M)$ is homotopy equivalent to the 
$S$-construction of vector bundles on $M$ by the map forgetting 
metrics.
In particular, it follows that 
$$
\pi _{n+1}(\vert \widehat{S}(M)\vert )\simeq K_n(M).
$$
Let $\widetilde{\mathbb Z}\widehat{C}_*(M)=
\widetilde{\mathbb Z}C_*\widehat{\mathscr P}(M)$ and 
$\widetilde{\mathbb Z}\widehat{C}^{emi}_*(M)$ the subcomplex of 
$\widetilde{\mathbb Z}\widehat{C}_*(M)$ generated by emi-cubes on $M$.
Then $\mathcal F\mapsto \lambda \mathcal F$ induces a homomorphism 
of complexes 
$$
\lambda :\widetilde{\mathbb Z}\widehat{C}_*(M)\to 
\widetilde{\mathbb Z}\widehat{C}^{emi}_*(M).
$$

\vskip 1pc
\begin{thm}\cite{burwang} \ \ 
If $\mathcal F$ is an exact metrized $n$-cube on $M$, then we have 
$$
d_{\mathscr D}\CC _n(\mathcal F)=\CC _{n-1}(\partial F).
$$
Hence the higher Bott-Chern forms induce a homomorphism of complexes 
$$
\CC :\widetilde{\mathbb Z}\widehat{C}_*(M)\to 
\underset{p}{\oplus }\mathscr D^*(M, p)[2p].
$$
Moreover, the following map 
$$
K_n(M)=\pi _{n+1}(\widehat{S}(M))
\overset{Hurewicz}{\longrightarrow }
H_{n+1}(\widehat{S}(M))\overset{\CUB }{\longrightarrow }
H_n(\widetilde{\mathbb Z}\widehat{C}_*(M))\overset{\CC }{\to }
\underset{p}{\oplus }H_{\mathscr D}^{2p-n}(M, p)
$$
agrees with the higher Chern character with values in the Deligne 
cohomology.
\end{thm}

\vskip 1pc
This is the main theorem of \cite{burwang}.
Here let us prove 
$d_{\mathscr D}\CC _n(\mathcal F)=\CC _{n-1}(\partial \mathcal F)$ 
in a different way from \cite{burwang}.
To do this we introduce another description of the logarithmic forms 
$S_n^i$.

For integers $\alpha _1, \cdots , \alpha _k$ with 
$1\leq \alpha _i\leq n$, a $k$-form with logarithmic poles 
$(\alpha _1, \cdots , \alpha _k)$ on $(\mathbb P^1)^n$ is given as 
$$
(\alpha _1, \cdots , \alpha _k)=d\log |z_{\alpha _1}|^2\wedge 
\cdots \wedge d\log |z_{\alpha _k}|^2. 
$$
Let $(\alpha _1, \cdots , \alpha _k)^{(i,k-i)}$ denote 
the $(i,k-i)$-part of $(\alpha _1, \cdots , \alpha _k)$.
Then we have 
$$
S_n^i=(i-1)!(n-i)!\sum_{\alpha =1}^n(-1)^{\alpha +1}\log |z_{\alpha }|^2
(1, \cdots , \widehat{\alpha }, \cdots , n)^{(i-1,n-i)}.
$$
 
\vskip 1pc
\begin{lem} \ \ 
It follows that 
\begin{gather*}
\partial S_n^i=i!(n-i)!(1 \cdots n)^{(i, n-i)}+(n-i)\sum^n_{\alpha =1}
(-1)^{\alpha }\partial \overline{\partial }\log |z_{\alpha }|^2
S_{n-1, \widehat{\alpha }}^i,   \\
\overline{\partial }S_n^i=(i-1)!(n-i+1)!(1 \cdots n)^{(i-1, n-i+1)}
-(i-1)\sum^n_{\alpha =1}(-1)^{\alpha }\partial \overline{\partial }
\log |z_{\alpha }|^2S_{n-1, \widehat{\alpha }}^{i-1},
\end{gather*}
where $S_{n-1, \widehat{\alpha }}^i$ is the logarithmic form 
on $(\mathbb P^1)^n$ with the same expression as $S_{n-1}^i$ for 
the coordinate 
$(z_1, \cdots , z_{\alpha -1}, z_{\alpha +1}, \cdots , z_n)$.
\end{lem}

\vskip 1pc
We will prove more general identities in Lem.4.3.
Let us return to the proof of the identity 
$d_{\mathscr D}\CC _n(\mathcal F)=\CC _{n-1}(\partial \mathcal F)$.
Since $\lambda :\widetilde{\mathbb Z}\widehat{C}_*(M)\to 
\widetilde{\mathbb Z}\widehat{C}^{emi}_*(M)$ is a homomorphism of 
complexes, we may assume that $\mathcal F$ is emi.
Lem.1.6 implies 
\begin{align*}
\sum_{i=1}^n(-1)^idS_n^i&=\sum_{i=0}^n(-1)^i(\partial S_n^i-
\overline{\partial }S_n^{i+1})  \\
&=n\sum_{i=1}^{n-1}(-1)^i\sum_{\alpha =1}^n(-1)^{\alpha }
\partial \overline{\partial }\log |z_{\alpha }|^2
S_{n-1, \widehat{\alpha }}^i.
\end{align*}
If $\delta _{\{z_{\alpha }=0\}}$ 
(resp.~$\delta _{\{z_{\alpha }=\infty \}}$) is the current on 
$(\mathbb P^1)^n$ given by the integration on the subvariety 
$\{z_{\alpha }=0\}$ (resp.~$\{z_{\alpha }=\infty \}$), 
then $\partial \overline{\partial }\log |z_{\alpha }|^2=
-(2\pi \sqrt{-1})(\delta _{\{z_{\alpha }=0\}}-
\delta _{\{z_{\alpha }=\infty \}})$.
Hence we have 
\begin{align*}
d_{\mathscr D}\CC _n(\mathcal F)
&=\frac{(-1)^{n+1}}{2(2\pi \sqrt{-1})^nn!}\int_{(\mathbb P^1)^n}
\CC _0(\TTR _n\lambda \mathcal F)\sum_{i=1}^n(-1)^idS_n^i  \\
&=\frac{(-1)^{n+1}}{2(2\pi \sqrt{-1})^n(n-1)!}\int_{(\mathbb P^1)^n}
\CC _0(\TTR _n\lambda \mathcal F)\sum_{i=1}^{n-1}(-1)^i
\sum_{\alpha =1}^n(-1)^{\alpha }\partial \overline{\partial }
\log |z_{\alpha }|^2S_{n-1, \widehat{\alpha }}^i \\
&=\frac{(-1)^{n+1}}{2(2\pi \sqrt{-1})^{n-1}(n-1)!}
\int_{(\mathbb P^1)^{n-1}}\CC _0(\TTR _{n-1}(\lambda \partial \mathcal F))
\sum_{i=1}^n(-1)^idS_{n-1}^i  \\
&=\CC _{n-1}(\partial \mathcal F).
\end{align*}

\vskip 1pc
{\it Remark}: \ 
Burgos and Wang defined the Bott-Chern forms of metrized cubes 
on complex algebraic manifolds that are not necessarily compact 
in \cite{burwang}.
They indeed developed the theory of higher Bott-Chern forms in 
a different way from the above and obtained an integral expression of them 
in the compact case.
In this paper, we adopted their integral expression as the definition 
of Bott-Chern forms.
But we should notice that our definition is the minus of their expression.

\vskip 1pc
\subsection{Renormalization} \ \ 
The Bott-Chern secondary classes used in Arakelov geometry 
is not equal to $\CC _1(\mathscr E)$.
This disagreement stems from the definition of Chern forms.
Gillet and Soul\'{e} have defined the Chern form of a hermitian 
vector bundle $\overline{E}$ as 
$$
\CCGS _0(\overline{E})=\TR \left(\exp 
\left(-\frac{K_{\overline{E}}}{2\pi \sqrt{-1}}\right)\right)
\in \underset{p}{\oplus }\mathscr E_{\mathbb R}^{2p}(M)\cap 
\mathscr E^{p,p}(M)\cap \KER d.
$$
Then it follows that $\CCGS _0(\overline{E})^{(p,p)}=
\frac{1}{(2\pi \sqrt{-1})^p}\CC _0(\overline{E})^{(p,p)}$.
Therefore to make it possible to apply the theory of higher Bott-Chern 
forms to Arakelov geometry, we have to redefine them as to be 
real forms.

We set 
$$
\mathscr A^n(M, p)=\begin{cases}
\mathscr E^{n-1}_{\mathbb R}(M)\cap 
\underset{p^{\prime }<p, q^{\prime }<p}
{\underset{p^{\prime }+q^{\prime }=n-1}{\oplus }}
\mathscr E^{p^{\prime }, q^{\prime }}(M),  &  n<2p,  \\
\mathscr E^{2p}_{\mathbb R}(M)\cap 
\mathscr E^{p, p}(M)\cap \KER d,  &  n=2p,  \\
0, & n>2p
\end{cases}
$$
and define a differential 
$d_{\mathscr A}:\mathscr A^n(M, p)\to \mathscr A^{n+1}(M, p)$ as 
$$
d_{\mathscr A}(\omega )=\begin{cases}
-\pi (d\omega ),  &  n<2p-1,  \\
dd^cx=\frac{-1}{2\pi \sqrt{-1}}\partial \overline{\partial}\omega , &n=2p-1,\\
0, & n>2p-1.  \end{cases}
$$
Then $(\mathscr A^*(M, p), d_{\mathscr A})$ becomes a complex. 
The renormalization operator $\Theta _{p,n}:\mathscr D^n(M, p)\to 
\mathscr A^n(M, p)$ is defined as 
$$
\Theta _{p,n}(\omega )=\begin{cases}
\frac{1}{(2\pi \sqrt{-1})^p}\omega , &n=2p,  \\
\frac{2}{(2\pi \sqrt{-1})^{p-1}}\omega , &n<2p.
\end{cases}
$$
Then $\Theta =\Theta _{p,*}:\mathscr D^*(M, p)\to \mathscr A^*(M, p)$ 
becomes an isomorphism of complexes.

The renormalized Bott-Chern form of an exact metrized $n$-cube 
$\mathcal F$ is defined as 
$$
\CCGS _n(\mathcal F)=\Theta (\CC _n(\mathcal F))\in \underset{p}{\oplus }
\mathscr A^{2p-n}(M, p).
$$
Then we have a homomorphism of complexes
$$
\CCGS :\widetilde{\mathbb Z}\widehat{C}_*(M)\to 
\underset{p}{\oplus }\mathscr A^*(M, p)[2p].
$$

\vskip 1pc
\begin{prop} \ \ 
The Bott-Chern secondary characteristic class of $\mathscr E$ defined 
in \cite{gilsoul} is equal to $-\CCGS _1(\mathscr E)$ modulo 
$\IIm d_{\mathscr A}$.
\end{prop}

{\it Proof}: \ 
For an exact metrized $1$-cube $\mathscr E$, we have 
\begin{align*}
\CCGS _1(\mathscr E)^{(p-1,p-1)}&=\frac{2}{(2\pi \sqrt{-1})^{p-1}}
\CC _1(\mathscr E)^{(p-1,p-1)}   \\
&=\frac{2}{(2\pi \sqrt{-1})^p}\int_{\mathbb P^1}
\CC _0(\TTR _1\lambda \mathscr E)^{(p,p)}\times \frac{1}{2}\log |z|^2 \\
&=\int_{\mathbb P^1}\CCGS _0(\TTR _1\lambda \mathscr E)^{(p,p)}\log |z|^2,
\end{align*}
which completes the proof.
\qed

\vskip 1pc
\subsection{The case of arithmetic varieties} \ \ 
Let $X$ be a proper arithmetic variety, that is, we assume that $X$ is 
proper over $\mathbb Z$. 
Let $X(\mathbb C)$ denote the compact complex manifold consisting of 
$\mathbb C$-valued points on $X$ and $F_{\infty }$ the complex 
conjugation on $X(\mathbb C)$. 
The real Deligne cohomology of $X$ is defined as 
$$
H^n_{\mathscr D}(X, \mathbb R(p))=
H^n_{\mathscr D}(X(\mathbb C), \mathbb R(p))^{\overline{F}^*_{\infty }=\IId }.
$$
Hence if we set 
$$
\mathscr D^n(X, p)=
\mathscr D^n(X(\mathbb C), p)^{\overline{F}^*_{\infty }=\IId },
$$
then we have an isomorphism 
$$
H^n(\mathscr D^*(X, p), d_{\mathscr D})\simeq 
H_{\mathscr D}^n(X, \mathbb R(p)).
$$

A {\it hermitian vector bundle} on $X$ is a pair $\overline{E}=(E, h)$ 
of a vector bundle $E$ on $X$ and an $F_{\infty }$-invariant 
hermitian metric $h$ on the holomorphic vector bundle $E(\mathbb C)$ on 
$X(\mathbb C)$.
An {\it exact metrized $n$-cube} on $X$ is an exact $n$-cube made  
of hermitian vector bundles on $X$.
Since the Chern form $\CC _0(\overline{E})$ is in 
$\underset{p}{\oplus }\mathscr D^{2p}(X, p)$, 
the Bott-Chern form of an exact metrized $n$-cube 
is in $\underset{p}{\oplus }\mathscr D^{2p-n}(X, p)$. 

Let $\widehat{\mathscr P}(X)$ be the category of hermitian vector 
bundles on $X$ and $\widehat{S}(X)$ the $S$-construction of 
$\widehat{\mathscr P}(X)$.
Then we have a canonical isomorphism 
$\pi _{n+1}(\vert \widehat{S}(X)\vert )\simeq K_n(X)$.
Let $\widetilde{\mathbb Z}\widehat{C}_*(X)=\widetilde{\mathbb Z}
C_*\widehat{\mathscr P}(X)$. 
Then we can define a homomorphism of complexes 
$$
\CC :\mathbb Z\widehat{S}_*(X)\overset{\CUB }{\longrightarrow }
\widetilde{\mathbb Z}\widehat{C}_*(X)[1]
\overset{\CC [1]}{\longrightarrow }\underset{p}{\oplus }
\mathscr D^*(X, p)[2p+1].
$$
If we set 
\begin{align*}
\mathscr A^n(X, p)&=\Theta _{p,n}(\mathscr D^n(X, p))  \\
&=\begin{cases}
\mathscr A^n(X(\mathbb C), p)^{F_{\infty }^*=(-1)^{p-1}}, &n<2p,  \\
\mathscr A^{2p}(X(\mathbb C), p)^{F_{\infty }^*=(-1)^p}, &n=2p,   \\
0, &n>2p,   \end{cases}
\end{align*}
then we can also obtain a homomorphism of complexes 
$$
\CCGS :\mathbb Z\widehat{S}_*(X)\overset{\CUB }{\longrightarrow }
\widetilde{\mathbb Z}\widehat{C}_*(X)[1]
\overset{\CCGS [1]}{\longrightarrow }\underset{p}{\oplus }
\mathscr A^*(X, p)[2p+1].
$$

\vskip 2pc
\section{Modified homotopy groups}
\vskip 1pc

\subsection{Definition of modified homotopy groups} \ \ 
In this section we develop a general framework used later in this paper.
Let $I$ be the closed interval 
$[0, 1]$ equipped with the usual CW-complex structure.
Throughout this paper we identify the sphere $S^n$ with 
$I^n/\partial I^n$.
Therefore $S^n$ consists of two cells as a CW-complex and any point 
of $S^n$ except the base point is expressed by $n$-tuple of real 
numbers $(t_1, \cdots , t_n)$ with $0<t_i<1$.  

Let $T$ be a pointed CW-complex with the base point $*\in T$.
The $n$-th homotopy group of $T$ is defined as the set of 
homotopy equivalence classes of pointed continuous maps 
$S^n\to T$.
However the cellular approximation theorem implies that continuous 
maps and homotopies can be replaced with cellular ones in 
the definition of homotopy groups.

Let $\SK _n(T)$ be the $n$-th skeleton of $T$ when 
$n\geq 0$ and $\SK _{-1}(T)=\{*\}$.
For $n\geq 0$, we set $C_n(T)=H_n(\SK _n(T), \SK _{n-1}(T))$, the $n$-th 
relative homology group of the pair \linebreak 
$(\SK _n(T), \SK _{n-1}(T))$.
Let $\partial : C_n(T)\to C_{n-1}(T)$ be the connecting 
homomorphism for the triple $(\SK _n(T), \SK _{n-1}(T), \SK _{n-2}(T))$. 
Then $(C_*(T), \partial )$ is a homological complex whose homology 
groups are isomorphic to the reduced homology groups of $T$.

Suppose that we have  a homological complex of abelian groups 
$(W_*, \partial )$ and a homomorphism of complexes 
$\rho :C_*(T)\to W_*$.
Let $\widetilde{W}_n=W_n/\IIm \partial $.
Let us consider a pair $(f, \omega )$, where $f:S^n\to T$ is a pointed 
cellular map and $\omega \in \widetilde{W}_{n+1}$.
A {\it cellular homotopy} from one pair $(f, \omega )$ to another pair 
$(f^{\prime }, \omega ^{\prime })$ 
is a pointed cellular map  $H:S^n\times I/\{*\}\times I\to T$ 
satisfying the following: 
\begin{enumerate}  
\item $H(x, 0)=f(x)$ and $H(x, 1)=f^{\prime }(x)$.
\item Let $[S^n\times I]\in C_{n+1}(S^n\times I)$ denote the fundamental 
chain of $S^n\times I$, where the orientation on $S^n\times I$
is inherited from the canonical orientation of the interval $I$.
Then we have 
$$
\omega ^{\prime }-\omega =(-1)^{n+1}\rho H_*([S^n\times I]).
$$
\end{enumerate}
It is shown that the cellular homotopy gives an equivalence 
relation on the set of pairs. 
Such two pairs are said to be {\it homotopy equivalent}.
We denote by $\widehat{\pi }_n(T, \rho )$ the set of all homotopy 
equivalence classes of pairs.

Let us define a multiplication on the set $\widehat{\pi }_n(T, \rho )$.
Let $T\vee T=\{(x, y)\in T\times T; x=* \ \text{or} \ y=*\}$ for 
a pointed CW-complex $T$.
Then there is a natural map $T\vee T\to T$ by $(x, *)\mapsto x$ 
and $(*, y)\mapsto y$.
A comultiplication map $\mu :S^n\to S^n\vee S^n$ is defined as 
$$
\mu (t_1, \cdots , t_n)=\begin{cases}
\left((t_1, t_2, \cdots , 2t_n), *\right), &0<t_n\leq \frac{1}{2}, \\
\left(*, (t_1, t_2, \cdots , 2t_n-1)\right), &\frac{1}{2}\leq t_n<1  
\end{cases}
$$
and a homotopy inverse map $\nu :S^n\to S^n$ by 
$\nu (t_1, \cdots , t_{n-1}, t_n)=(t_1, \cdots , t_{n-1}, 1-t_n)$.
For two pointed cellular maps $f, g:S^n\to T$, a 
multiplication $f*g$ is defined as 
$$
f*g: S^n\overset{\mu }{\to }S^n\vee S^n\overset{f\vee g}{\longrightarrow }
T\vee T\to T
$$
and an inverse $f^{-1}$ as 
$$
f^{-1}: S^n\overset{\nu }{\to }S^n\overset{f}{\to }T.
$$
For two pairs $(f, \omega )$ and $(g, \tau )$ where 
$f, g:S^n\to T$ are pointed cellular maps and 
$\omega , \tau \in \widetilde{W}_{n+1}$, a multiplication 
$(f, \omega )*(g, \tau )$ is defined as 
$$
(f, \omega )*(g, \tau )=(f*g, \omega +\tau ).
$$
It is easy to show that the multiplication $*$ is compatible with 
the homotopy equivalence relation of pairs.
Hence it gives rise to a multiplication on 
$\widehat{\pi }_n(T, \rho )$, which is also denoted by $*$.

We next verify the associativity of the multiplication $*$.
For three pointed cellular maps $f, g, h:S^n\to T$, a cellular homotopy 
$H_1:S^n\times I/\{*\}\times I\to T$ from $(f*g)*h$ to $f*(g*h)$ is 
given as follows:
$$
H_1(t_1, \cdots , t_{n-1}, t_n, u)=\begin{cases}
  f(t_1, \cdots , t_{n-1}, \frac{4t_n}{u+1}), 
       & 0<t_n\leq \frac{u+1}{4}, \\
  g(t_1, \cdots , t_{n-1}, 4t_n-u-1), 
      & \frac{u+1}{4}\leq t_n\leq \frac{u+2}{4}, \\
  h(t_1, \cdots , t_{n-1}, \frac{4t_n-2-u}{2-u}), 
      & \frac{u+2}{4}\leq t_n<1.
 \end{cases}
$$
Since the image of $H_1$ is contained in $sk_n(T)$, we have 
$(H_1)_*([S^n\times I])=0$ in $C_{n+1}(T)$.
Hence $H_1$ becomes a cellular homotopy from 
$\left((f, \omega )*(g, \tau )\right)*(h, \eta )$ to 
$(f, \omega )*\left((g, \tau )*(h, \eta )\right)$ for any 
$\omega , \tau , \eta \in \widetilde{W}_{n+1}$. 

Finally we show the existence of unit and inverse in 
$\widehat{\pi }_n(t, \rho )$ with respect to the multiplication $*$.
Let $0:S^n\to T$ be the map defined by $0(S^n)=*$.
For a pointed cellular map $f:S^n\to T$, a homotopy $H_2$ from $f*0$ to 
$f$ is given as 
$$
H_2(t_1, \cdots , t_{n-1}, t_n , u)=\begin{cases}
  f(t_1, \cdots , t_{n-1}, \frac{2t_n}{u+1}), 
      & 0<t_n\leq \frac{u+1}{2}, \\
  *, & \frac{u+1}{2}\leq t_n<1.
 \end{cases}
$$
A homotopy $H_3$ from $0*f$ to $f$ can be given similarly. 
Moreover, a homotopy $H_4$ from 
$f*f^{-1}$ to $0$ is given as 
$$
H_4(t_1, \cdots , t_{n-1}, t_n, u)=\begin{cases}
  f(t_1, \cdots , t_{n-1}, \frac{2t_n}{1-u}), 
       & 0<t_n\leq \frac{1-u}{2}, \\
  *,  & \frac{1-u}{2}\leq t_n\leq \frac{1+u}{2}, \\
  f(t_1, \cdots , t_{n-1}, \frac{-2t_n+2}{1-u}), 
      & \frac{u+1}{2}\leq t_n<1,
 \end{cases}
$$
and a homotopy $H_5$ from $f^{-1}*f$ to $0$ can be given similarly.
These homotopies are all cellular and their images are 
contained in $sk _n(T)$.
Hence $(f, \omega )*(0, 0)$ and $(0, 0)*(f, \omega )$ 
are homotopy equivalent to $(f, \omega )$ and 
$(f, \omega )*(f^{-1}, -\omega )$ and $(f^{-1}, -\omega )*(f, \omega )$ 
are homotopy equivalent to $(0, 0)$.

\vskip 1pc
\begin{thm} \ \ 
If $n\geq 1$, $\widehat{\pi }_n(T, \rho )$ is a group by 
the multiplication $*$ and when 
$n\geq 2$, it becomes an abelian group.
\end{thm}

{\it Proof}: \ 
We have already shown that the multiplication 
$$
[(f, \omega )]*[(g, \tau )]=[(f*g, \omega +\tau )]
$$
provides $\widehat{\pi }_n(T, \rho )$ with the structure of a group.
For two pointed cellular maps $f, g:S^n\to T$, $f*g$ and $g*f$ are 
homotopy equivalent if $n\geq 2$ and a homotopy between them is 
given in every textbook of homotopy theory.
Although it is too complicated to write it down, it is easy to see  
that the image of this homotopy is contained in $sk _n(T)$. 
Hence $(f*g, 0)$ is homotopy equivalent to $(g*f, 0)$, 
therefore $\widehat{\pi }_n(T, \rho )$ is an abelian group if $n\geq 2$. 
\qed

\vskip 1pc
\begin{defn} \ \
The group $\widehat{\pi }_n(T, \rho )$ is called the $n$-th homotopy 
group of $T$ modified by the homomorphism $\rho $.
\end{defn}

\vskip 1pc
Let $\zeta :\widehat{\pi }_n(T, \rho )\twoheadrightarrow \pi _n(T)$ 
denote the surjection obtained by forgetting elements of 
$\widetilde{W}_{n+1}$.
Then we have the following:

\vskip 1pc
\begin{thm} \ \ 
There is an exact sequence 
$$
\pi _{n+1}(T)\overset{\rho }{\to }\widetilde{W}_{n+1}
\overset{a}{\to }\widehat{\pi }_n(T, \rho )\overset{\zeta }{\to }
\pi _n(T)\to 0,
$$
where the map $\rho $ is defined as 
$$
\rho :\pi _{n+1}(T)\overset{Hurewicz}{\longrightarrow }
H_{n+1}(T)\overset{H_{n+1}(\rho )}{\longrightarrow }H_{n+1}(W_*)
\subset \widetilde{W}_{n+1}
$$
and the map $a$ as 
$$
a(\omega )=[(0, \omega )]\in \widehat{\pi }_n(T, \rho ).
$$
\end{thm}

{\it Proof}: \ 
The cellular approximation theorem implies $\IIm a=\KER \zeta $.
Hence we have only to prove $\KER a=\IIm \rho $.
For an element $\omega \in \widetilde{W}_{n+1}$, the pair 
$(0, \omega )$ is homotopy equivalent to $(0, 0)$ if and only if 
there is a cellular homotopy $H:S^n\times I/\{*\}\times I\to T$ 
from $0$ to $0$ such that $(-1)^{n+1}\rho H_*([S^n\times I])=\omega $.
Since $H(S^n\times \partial I)=*$, $H$ induces a pointed cellular map 
$H^{\prime }:S^{n+1}\to T$.
Then $\omega $ is equal to the image of 
$(-1)^{n+1}[H^{\prime }]\in \pi _{n+1}(T)$ by $\rho $, therefore 
we have $\KER a\subset \IIm \rho $.
The opposite inclusion $\IIm \rho \subset \KER a$ can be obtained 
by viewing a pointed cellular map 
$S^{n+1}\to T$ as a cellular homotopy between the collapsing map $0$.
\qed

\vskip 1pc
\subsection{A homomorphism from a modified homotopy group} \ \ 
For a pair $(f, \omega )$ as in the last subsection, 
$\rho (f, \omega )\in W_n$ is defined by 
$$
\rho (f, \omega )=\rho f_*([S^n])+\partial \omega .
$$

\vskip 1pc
\begin{prop} \ \ 
The above $\rho (f, \omega )$ gives rise to the homomorphism 
$$
\rho :\widehat{\pi }_n(T, \rho )\to W_n
$$
and $\IIm \rho $ is contained in $\KER (\partial :W_n\to W_{n-1})$.
\end{prop}

{\it Proof}: \ 
Let us show that $\rho (f, \omega )$ is compatible with the homotopy 
equivalence relation. 
If $H:S^n\times I/\{*\}\times I\to T$ is a cellular homotopy 
from $(f, \omega )$ to $(f^{\prime }, \omega ^{\prime })$, then we have 
$$
\partial H_*([S^n\times I])=(-1)^n(f^{\prime }_*([S^n])-f_*([S^n]))
$$
in $C_n(T)$ and $\rho H_*([S^n\times I])=(-1)^{n+1}
(\omega ^{\prime }-\omega )$.
Hence we have 
\begin{align*}
\rho (f, \omega )&=\rho f_*([S^n])+\partial \omega    \\
&=\rho f^{\prime }_*([S^n])+(-1)^{n+1}
\partial \rho H_*([S^n\times I]))+\partial \omega    \\
&=\rho f^{\prime }_*([S^n])+\partial (\omega ^{\prime }-\omega )
+\partial \omega  \\
&=\rho (f^{\prime }, \omega ^{\prime }).
\end{align*}
The inclusion $\IIm \rho \subset \KER (\partial :W_n\to W_{n-1})$ and 
the claim that $\rho $ is a homomorphism of groups are obvious.
\qed

\vskip 1pc
The exact sequence in Thm.2.3 implies the following corollaries:

\vskip 1pc
\begin{cor} \ \ 
There is an exact sequence 
$$
\pi _{n+1}(T)\overset{\rho }{\to }H_{n+1}(W_*, \partial )
\overset{a}{\to }\widehat{\pi }_n(T, \rho )
\overset{\zeta \oplus \rho }{\longrightarrow }\pi _n(T)\oplus 
\KER \partial \overset{cl}{\to }H_n(W_*, \partial )\to 0,
$$
where $\KER \partial =\KER (\partial :W_n\to W_{n-1})$ and 
$cl$ is defined as $cl(x, \omega )=\rho (x)-[\omega ]$.
\end{cor}

\vskip 1pc
\begin{cor} \ \ 
For $n\geq 1$, let 
$$
\widehat{\pi }_n(T, \rho )_0=\KER (\rho :\widehat{\pi }_n(T, \rho )
\to W_n).
$$
Then there is a long exact sequence 
$$
\cdots \overset{\zeta }{\to }\pi _{n+1}(T)\overset{\rho }{\to }
H_{n+1}(W_*, \partial )\overset{a}{\to }\widehat{\pi }_n(T, \rho )_0
\overset{\zeta }{\to }\pi _n(T)\overset{\rho }{\to }\cdots .
$$
\end{cor} 

{\it Proof}: \ 
The exactness at $\widehat{\pi }_n(T, \rho )_0$ and 
$H_{n+1}(W_*, \partial )$ has already been verified in Thm.2.3.
For $[(f, \omega )]\in \widehat{\pi }_n(T, \rho )_0$, we have 
$$
\rho \zeta ([(f, \omega )])=\rho f_*([S^n])=-[\partial \omega ]=0.
$$
Conversely, if a pointed cellular map $f:S^n\to T$ satisfies 
$\rho f_*([S^n])=\partial \omega $ for some $\omega \in W_{n+1}$, 
then $[(f, -\omega )]\in \widehat{\pi }_n(T, \rho )_0$ and 
$\zeta ([(f, -\omega )])=[f]\in \pi _n(T)$.
Hence the exactness at $\pi _n(T)$ follows.
\qed

\vskip 1pc
\subsection{A Functorial property of modified homotopy groups} \ \ 
Let $T$ and $T^{\prime }$ be pointed CW-complexes. 
Let $\rho :C_*(T)\to W_*$ and 
$\rho ^{\prime }:C_*(T^{\prime })\to {W^{\prime }}_*$ be homomorphisms 
of complexes. 
Given a pointed cellular map $\alpha :T\to T^{\prime }$ and 
a homomorphism of complexes $\beta :W_*\to {W^{\prime }}_*$ that 
make the diagram 
$$
\begin{CD}
C_*(T) @>{\alpha _*}>> C_*(T^{\prime })  \\
@VV{\rho }V @VV{\rho ^{\prime }}V  \\
W_* @>{\beta }>> {W^{\prime }}_* 
\end{CD}
$$
commutative, a homomorphism of modified homotopy groups 
$$
(\alpha , \beta )_*:\widehat{\pi }_n(T, \rho )\to 
\widehat{\pi }_n(T^{\prime }, \rho ^{\prime })
$$
is defined by $[(f, \omega )]\mapsto [(\alpha f, \beta (\omega ))]$.

We can however define a homomorphism of modified homotopy groups 
under a weaker assumption than the strict commutativity of the diagram. 
We assume that the above diagram is commutative up to homotopy. 
In other words, we assume the existence of a homomorphism 
$$
\Phi :C_*(T)\to {W^{\prime }}_{*+1}
$$
satisfying $\rho ^{\prime }\alpha _*-\beta \rho =
\partial \Phi +\Phi \partial $. 

\vskip 1pc
\begin{prop} \ \ 
Under the above notations, we can define a homomorphism 
$$
(\alpha , \beta , \Phi )_*:\widehat{\pi }_n(T, \rho )\to 
\widehat{\pi }_n(T^{\prime }, \rho ^{\prime })
$$
by $[(f, \omega )]\mapsto [(\alpha f, \beta (\omega )
-\Phi f_*([S^n]))]$. 
This homomorphism enjoys the following functorial property: 
Let $\alpha :T\to T^{\prime }$ and $\alpha ^{\prime }:T^{\prime }\to 
T^{\prime \prime }$ be pointed cellular maps and let 
$\beta :W_*\to W_*^{\prime }$ and $\beta ^{\prime }:W_*^{\prime }\to 
W_*^{\prime \prime }$ be homomorphisms of complexes.
We assume that the squares
$$
\begin{CD}
C_*(T) @>{\alpha _*}>> C_*(T^{\prime }) @>{\alpha ^{\prime }_*}>> 
C_*(T^{\prime \prime }) \\
@VV{\rho }V @VV{\rho ^{\prime }}V @VV{\rho ^{\prime \prime }}V \\
W_* @>{\beta }>> {W^{\prime }}_* @>{\beta ^{\prime }}>> 
{W^{\prime \prime }}_*
\end{CD}
$$
are commutative up to homotopy and let $\Phi $ and $\Phi ^{\prime }$ 
be homotopies of these squares.
Then we have 
$$
(\alpha ^{\prime }, \beta ^{\prime }, \Phi ^{\prime })_*
(\alpha , \beta , \Phi )_*=(\alpha ^{\prime }\alpha , 
\beta ^{\prime }\beta , \beta ^{\prime }\Phi 
+\Phi ^{\prime }\alpha _*)_*:\widehat{\pi }_n(T, \rho )\to 
\widehat{\pi }_n(T^{\prime \prime }, \rho ^{\prime \prime }).
$$
\end{prop}

{\it Proof}: \ 
Let $f, f^{\prime }:S^n\to T$ be pointed cellular maps and 
$\omega , \omega ^{\prime }\in \widetilde{W}_{n+1}$.
If $H:S^n\times I/\{*\}\times I\to T$ is a cellular homotopy from 
$(f, \omega )$ to $(f^{\prime }, \omega ^{\prime })$, then we have 
\begin{align*}
(-1)^{n+1}\rho ^{\prime }\alpha _*H_*([S^n\times I])&=(-1)^{n+1}
\beta \rho H_*([S^n\times I])+(-1)^{n+1}\partial \Phi H_*([S^n\times I]) \\*
&\hskip 1pc +(-1)^{n+1}\Phi \partial H_*([S^n\times I])  \\
&\equiv (\beta (\omega ^{\prime })-\Phi f^{\prime }_*([S^n]))-
(\beta (\omega )-\Phi f_*([S^n])) 
\end{align*}
modulo $\IIm \partial $. 
This tells that the map $\alpha H:S^n\times I/\{*\}\times I\to 
T^{\prime }$ is a cellular homotopy from 
$(\alpha f, \beta (\omega )-\Phi f_*([S^n]))$ to 
$(\alpha f^{\prime }, \beta (\omega ^{\prime })
-\Phi f^{\prime }_*([S^n]))$.
Hence $(\alpha , \beta , \Phi )_*$ is well-defined.
The last identity can be verified by an easy calculation.
\qed

\vskip 1pc
\begin{prop} \ \ 
Under the above notations, we have a commutative diagram 
$$
\begin{CD}
\widehat{\pi }_n(T, \rho ) @>{\rho }>>  W_n  \\
@VV{(\alpha , \beta , \Phi )_*}V @VV{\beta }V  \\
\widehat{\pi }_n(T^{\prime }, \rho ^{\prime }) @>{\rho ^{\prime }}>> 
W^{\prime }_n.
\end{CD}
$$
\end{prop}

{\it Proof}: \ 
For a pointed cellular map $f:S^n\to T$, we have 
$$
\rho ^{\prime }\alpha _*f_*([S^n])-\beta \rho f_*([S^n])=
\partial \Phi f_*([S^n]).
$$
Hence we have 
\begin{align*}
\rho ^{\prime }(\alpha , \beta , \Phi )_*([(f, \omega )])
&=\rho ^{\prime }([(\alpha f, \beta (\omega )-\Phi f_*([S^n]))])  \\
&=\rho ^{\prime }\alpha _*f_*([S^n])
+\partial (\beta (\omega )-\Phi f_*([S^n])) \\
&=\beta (\rho f_*([S^n])+\partial \omega )  \\
&=\beta \rho ([(f, \omega )]), 
\end{align*}
which completes the proof.
\qed

\vskip 2pc
\section{Definition of arithmetic $K$-groups}
\vskip 1pc

\subsection{Triviality of the Bott-Chern form of a degenerate $n$-cell} \ \ 
Let $\mathfrak A$ be a small exact category and $\mathfrak S_n$ 
the $n$-th symmetric group.
For $\sigma \in \mathfrak S_n$ and an exact $n$-cube $\mathcal F$ 
of $\mathfrak A$, an exact $n$-cube $\sigma \mathcal F$ is defined as
$$
(\sigma \mathcal F)_{\alpha_1, \cdots , \alpha _n}=
\mathcal F_{\alpha _{\sigma (1)}, \alpha _{\sigma (2)}, \cdots , 
\alpha _{\sigma (n)}}.
$$

\vskip 1pc
\begin{lem} \ \ 
For any $E\in S_n\mathfrak A$, we have 
\begin{align*}
\CUB (s_0E)&=s_1^{-1}\CUB (E),   \\
\CUB (s_nE)&=s_n^1\CUB (E)   \\
\intertext{and if $1\leq i\leq n-1$, then we have} 
\CUB (s_iE)&=\tau _i\CUB (s_iE), 
\end{align*}
where $\tau _i\in \mathfrak S_n$ is the transposition of $i$ and $i+1$.
\end{lem}

{\it Proof}: \ 
To prove the lemma, we use the following identities, which is proved in 
Prop.4.5 in \cite{burwang}:
For $E\in S_n\mathfrak A$, it follows that 
\begin{align*}
\partial _i^{-1}\CUB (E)&=s_{n-2}^1\cdots s_i^1\CUB (\partial _{i+1}
\cdots \partial _nE),   \\
\partial _i^0\CUB (E)&=\CUB (\partial _iE),    \\
\partial _i^1\CUB (E)&=s_{i-1}^{-1}\cdots s_1^{-1}\CUB (\partial _0
\cdots \partial _{i-1}E)
\end{align*}
for $1\leq i\leq n-1$.
The first identity of the lemma follows from 
\begin{align*}
\partial ^{-1}_1\CUB (s_0E)&=s^1_{n-1}\cdots s^1_1\CUB 
(\partial _2\cdots \partial _{n+1}s_0E)  \\
&=s^1_{n-1}\cdots s^1_1\CUB (s_0\partial _1\cdots \partial _nE)  \\
&=0,  \\
\partial _1^0\CUB (s_0E)&=\CUB (\partial _1s_0E)=\CUB (E)  \\
\intertext{and} 
\partial _1^1\CUB (s_0E)&=\CUB (\partial _0s_0E)=\CUB (E).
\end{align*}
The second one follows from 
\begin{align*}
\partial _n^{-1}\CUB (s_nE)&=\CUB (\partial _{n+1}s_nE)=\CUB (E),  \\
\partial _n^0\CUB (s_nE)&=\CUB (\partial _ns_nE)=\CUB (E)  \\
\intertext{and} 
\partial ^1_n\CUB (s_nE)&=s^{-1}_{n-1}\cdots s^{-1}_1\CUB 
(\partial _0\cdots \partial _{n-1}s_nE)  \\
&=s^{-1}_{n-1}\cdots s^{-1}_1\CUB (s_0\partial _0\cdots \partial _{n-1}E) \\
&=0.
\end{align*}

We turn to the last one.
If $1\leq i\leq n-1$, then we have 
{\allowdisplaybreaks
\begin{align*}
\partial _i^{-1}\CUB (s_iE)&=s_{n-1}^1\cdots s_i^1\CUB (\partial _{i+1}
\cdots \partial _{n+1}s_iE)   \\
&=s_{n-1}^1\cdots s_i^1\CUB (\partial _{i+1}\cdots \partial _nE)   \\
&=s_{n-1}^1\cdots s_{i+1}^1\CUB (s_i\partial _{i+1}\cdots \partial _nE) \\
&=s_{n-1}^1\cdots s_{i+1}^1\CUB (\partial _{i+2}\cdots \partial _{n+1}
s_iE)   \\
&=\partial _{i+1}^{-1}\CUB (s_iE), \\
\partial _i^0\CUB (s_iE)&=\CUB (E)=\partial _{i+1}^0\CUB (s_iE)  \\
\intertext{and}
\partial _i^1\CUB (s_iE)&=s_{i-1}^{-1}\cdots s_1^{-1}\CUB (\partial _0
\cdots \partial _{i-1}s_iE)   \\
&=s_{i-1}^{-1}\cdots s_1^{-1}\CUB (s_0\partial _0\cdots \partial _{i-1}E) \\
&=s_{i-1}^{-1}\cdots s_1^{-1}s_1^{-1}
\CUB (\partial _0\cdots \partial _{i-1}E)   \\
&=s_i^{-1}\cdots s_1^{-1}\CUB (\partial _0\cdots \partial _{i-1}
\partial _is_iE)   \\
&=\partial _{i+1}^1\CUB (s_iE).
\end{align*}}
The last identity follows from them.
\qed

\vskip 1pc
Let $S_n\subset \mathbb ZC_n\mathfrak A$ be the subgroup generated by 
exact $n$-cubes $\mathcal F$ such that $\tau _i\mathcal F=\mathcal F$ 
for some integer $i$ satisfying $1\leq i\leq n-1$.
We set 
$$
\CUB _n(\mathfrak A)=\mathbb ZC_n\mathfrak A/(D_n+S_n).
$$

\vskip 1pc
\begin{lem} \ \ 
We have $\partial S_n\subset S_{n-1}$. 
Hence $\CUB _*(\mathfrak A)=(\CUB _n(\mathfrak A), \partial )$ 
becomes a homological complex.
\end{lem} 

{\it Proof}: \ 
Let $\mathcal F$ be an exact $n$-cube satisfying 
$\tau _i\mathcal F=\mathcal F$. 
If $k<i$, then $\partial _k^j\mathcal F=\partial _k^j\tau _i\mathcal F=
\tau _{i-1}\partial _k^j\mathcal F$ and if $k>i+1$, then 
$\partial _k^j\mathcal F=\partial _k^j\tau _i\mathcal F
=\tau _i\partial _k^j\mathcal F$.
Furthermore, $\tau _i\mathcal F=\mathcal F$ implies 
$\partial _i^j\mathcal F=\partial _{i+1}^j\mathcal F$.
Hence we have 
$$
\partial \mathcal F=\sum_{k\not= i,i+1}\sum_{j=-1}^1(-1)^{k+j+1}
\partial _k^j\mathcal F\in S_{n-1},
$$
which completes the proof.
\qed

\vskip 1pc
\begin{lem} \ \ 
Let $\mathcal F$ be an exact metrized $n$-cube on a complex 
algebraic manifold $M$.
For any $\sigma \in \mathfrak S_n$, there is a canonical isometry 
of exact metrized $n$-cubes
$$
\sigma (\lambda \mathcal F)\simeq \lambda (\sigma \mathcal F).
$$
\end{lem}

{\it Proof}: \ 
It follows from the definition of $\lambda _i$ that 
$\sigma (\lambda _i\mathcal F)=
\lambda _{\sigma (i)}(\sigma \mathcal F)$. 
Therefore we have 
$$
\sigma (\lambda \mathcal F)=\lambda _{\sigma (n)}\cdots 
\lambda _{\sigma (1)}(\sigma \mathcal F).
$$
Hence it is sufficient to show the existence of a canonical isometry 
$\lambda _i\lambda _j\simeq \lambda _j\lambda _i$.
For simplicity, we prove it only in the case of $n=2$.

For an exact metrized $2$-cube $\mathcal F=\{\overline{E_{i,j}}\}$, 
$\lambda _2\lambda _1\mathcal F$ is given as follows:
$$
\footnotesize{
\begin{CD}
\overline{E_{-1,-1}}\oplus \overline{E_{1,-1}}\oplus \overline{E_{-1,1}}
\oplus \overline{E_{1,1}} @>>> \overline{E_{-1,0}}\oplus \overline{E_{1,0}}
\oplus \overline{E^{\prime }_{-1,1}}\oplus \overline{E^{\prime }_{1,1}}
@>>> \overline{E^{\prime }_{-1,1}}\oplus \overline{E^{\prime }_{1,1}} \\
@VVV @VVV @VVV \\
\overline{E_{0,-1}}\oplus \overline{E^{\prime }_{1,-1}}\oplus 
\overline{E_{0,1}}\oplus \overline{E^{\prime }_{1,1}}@>>> \overline{E_{0,0}}
\oplus \overline{E^{\prime }_{1,0}}\oplus \overline{E^{\prime }_{0,1}}
\oplus \overline{E^{\prime }_{1,1}} @>>> \overline{E^{\prime }_{0,1}}
\oplus \overline{E^{\prime }_{1,1}} \\
@VVV @VVV @VVV \\
\overline{E^{\prime }_{1,-1}}\oplus \overline{E^{\prime }_{1,1}} @>>> 
\overline{E^{\prime }_{1,0}}\oplus \overline{E^{\prime }_{1,1}} @>>> 
\overline{E^{\prime }_{1,1}}, 
\end{CD}
}
$$
where $\overline{E^{\prime }_{i,j}}$ is the same vector bundle as 
$\overline{E_{i,j}}$ equipped with the metric 
induced from $\overline{E_{0,0}}$. 
On the other hand, $\lambda _1\lambda _2\mathcal F$ 
is given as follows:
$$
\footnotesize{
\begin{CD}
\overline{E_{-1,-1}}\oplus \overline{E_{-1,1}}\oplus \overline{E_{1,-1}}
\oplus \overline{E_{1,1}} @>>> \overline{E_{-1,0}}\oplus 
\overline{E^{\prime }_{-1,1}}\oplus \overline{E_{1,0}}\oplus 
\overline{E^{\prime }_{1,1}}@>>> \overline{E^{\prime }_{-1,1}}\oplus 
\overline{E^{\prime }_{1,1}} \\
@VVV @VVV @VVV \\
\overline{E_{0,-1}}\oplus \overline{E_{0,1}}\oplus 
\overline{E^{\prime }_{1,-1}}\oplus \overline{E^{\prime }_{1,1}} @>>> 
\overline{E_{0,0}}\oplus \overline{E^{\prime }_{0,1}}\oplus 
\overline{E^{\prime }_{1,0}}\oplus \overline{E^{\prime }_{1,1}} @>>> 
\overline{E^{\prime }_{0,1}}\oplus \overline{E^{\prime }_{1,1}} \\
@VVV @VVV @VVV \\
\overline{E^{\prime }_{1,-1}}\oplus \overline{E^{\prime }_{1,1}} @>>> 
\overline{E^{\prime }_{1,0}}\oplus \overline{E^{\prime }_{1,1}} @>>> 
\overline{E^{\prime }_{1,1}}.
\end{CD}
}
$$
Hence an isometry $\lambda _2\lambda _1\mathcal F\simeq 
\lambda _1\lambda _2\mathcal F$ is given by appropriate 
permutations of direct \linebreak 
summands.
\qed

\vskip 1pc
\begin{thm} \ \ 
The Bott-Chern form of a degenerate element of $\widehat{S}_n(M)$ 
is zero.
\end{thm}

{\it Proof}: \ 
For an integer $i$ satisfying $1\leq i\leq n-1$, let 
$t_i:(\mathbb P^1)^n\to (\mathbb P^1)^n$ 
denote the involution interchanging the $i$-th and the $(i+1)$-th 
component.
Then by Prop.2.1 in \cite{roessler} and Lem.3.3, there is an isometry 
$$
t_i^*\TTR _n(\lambda \mathcal F)=
\TTR _n(\lambda \tau _i\mathcal F).
$$
Furthermore, it follows from the definition of $T_n$ that 
$t_i^*T_n=-T_n$.
Hence if $\tau _i\mathcal F=\mathcal F$, then we have 
\begin{align*}
\CC _n(\mathcal F)&=\int_{(\mathbb P^1)^n}\CC _0
(\TTR _n(\lambda \mathcal F))\wedge T_n   \\
&=\int_{(\mathbb P^1)^n}t_i^*(\CC _0(\TTR _n(\lambda \mathcal F))
\wedge T_n) \\
&=-\int_{(\mathbb P^1)^n}\CC _0(\TTR _n(\lambda \mathcal F))\wedge T_n \\
&=-\CC _n(\mathcal F), 
\end{align*}
therefore $\CC _n(\mathcal F)=0$. 

By Lem.3.1, the cube $\CUB (E)$ associated with a degenerate element 
$E\in \widehat{S}_n(X)$ is either a degenerate cube or a cube 
satisfying $\tau _i\mathcal F=\mathcal F$ for some $1\leq i\leq n-2$. 
Hence it follows from the above calculation and Prop.1.3 that 
$\CC _{n-1}(E)=0$.
\qed

\vskip 1pc
For a compact complex algebraic manifold $M$, let 
$\widehat{\CUB }_*(M)=\CUB _*(\widehat{\mathscr P}(M))$.
We set 
$$
\mathscr D_n(M)=\underset{p}{\oplus }\mathscr D^{2p-n}(M, p).
$$
Then $(\mathscr D_*(M), d_{\mathscr D})$ becomes a homological complex.
It follows from Thm.3.4 that the higher Bott-Chern forms induce 
a homomorphism 
$$
\CC :\widehat{\CUB }_*(M)\overset{\CC }{\longrightarrow }
\mathscr D_*(M).
$$

\vskip 1pc
\subsection{Definition of higher arithmetic $K$-theory} \ \ 
In this subsection we give the definition of higher arithmetic $K$-theory. 
Let $X$ be a proper arithmetic variety and $\mathscr A_n(X)=
\underset{p}{\oplus }\mathscr A^{2p-n}(X, p)$.
Then $(\mathscr A_*(X), d_{\mathscr A})$ is a homological complex 
and Thm.3.4 yields the homomorphism 
$$
\CCGS :C_*(\vert \widehat{S}(X)\vert )
\overset{\CUB }{\longrightarrow }\widehat{\CUB }_*(X)[1]
\overset{\CCGS }{\longrightarrow }\mathscr A_*(X)[1],
$$
where $\widehat{\CUB }_*(X)=\CUB _*(\widehat{\mathscr P}(X))$.

\vskip 1pc
\begin{defn} \ \ 
The $n$-th arithmetic $K$-group $\widehat{K}_n(X)$ of a proper 
arithmetic variety $X$ is defined as 
the $(n+1)$-th homotopy group of $\vert \widehat{S}(X)\vert $ 
modified by the renormalized higher Bott-Chern forms, that is, 
$$
\widehat{K}_n(X)=\widehat{\pi }_{n+1}(\vert \widehat{S}(X)\vert , \CCGS ).
$$
\end{defn}

\vskip 1pc
By Thm.2.3, we have the following:

\vskip 1pc
\begin{thm} \ \ 
There is an exact sequence
$$
K_{n+1}(X)\overset{\rho }{\to }\widetilde{\mathscr A}_{n+1}(X)\to 
\widehat{K}_n(X)\to K_n(X)\to 0,
$$
where $\widetilde{\mathscr A}_{n+1}(X)=\mathscr A_{n+1}(X)/
\IIm d_{\mathscr A}$ and $\rho $ is the Beilinson's regulator map 
up to a constant multiple.
\end{thm}

\vskip 1pc
The $0$-th arithmetic $K$-group has already been given 
by Gillet and Soul\'{e} in a different way \cite{gilsoul}.
They indeed defined it as the factor group of the free abelian group 
generated by $(\overline{E}, \omega )$, where $\overline{E}$ is 
a hermitian vector bundle on $X$ and 
$\omega \in \widetilde{A}_1(X)$, by the subgroup generated by 
$$
(\overline{E^{\prime }}, \omega )+
(\overline{E^{\prime \prime }}, \omega ^{\prime })-
(\overline{E}, \omega +\omega ^{\prime }-\CCGS _1(\mathscr E))
$$
for all short exact sequences $\mathscr E:0\to \overline{E^{\prime }}
\to \overline{E}\to \overline{E^{\prime \prime }}\to 0$.
In this paper, we denote this factor group by 
$\widehat{\mathscr K}_0(X)$.
We denote by $[(\overline{E}, \omega )]$ the element of 
$\widehat{\mathscr K}_0(X)$ determined by a pair $(\overline{E}, \omega )$.

\vskip 1pc
\begin{thm} \ \ 
There is a canonical isomorphism 
$$
\widehat{\alpha }:\widehat{\mathscr K}_0(X)\simeq \widehat{\pi }_1
(\vert \widehat{S}(X)\vert , \CCGS )=\widehat{K}_0(X).
$$
\end{thm}

{\it Proof}: \ 
This is an analogue of the fact that the Grothendieck 
group of a small exact category is isomorphic to the fundamental 
group of its $S$-construction.
Since $\widehat{S}_1(X)$ is the set of all hermitian vector bundles 
on $X$ and $\widehat{S}_0(X)=\{*\}$, a pointed simplicial loop 
$l_{\overline{E}}:S^1\to \vert \widehat{S}(X)\vert $ can be associated 
with a hermitian vector bundle $\overline{E}$ on $X$.
Moreover, any short exact sequence $\mathscr E:0\to 
\overline{E^{\prime }}\to \overline{E}\to \overline{E^{\prime \prime }}
\to 0$ gives a $2$-simplex $\Delta _{\mathscr E}$ in $\widehat{S}(X)$ 
whose faces are 
$\partial _0\Delta _{\mathscr E}=\overline{E^{\prime \prime }}$, 
$\partial _1\Delta _{\mathscr E}=\overline{E}$ and 
$\partial _2\Delta _{\mathscr E}=\overline{E^{\prime }}$.
If we regard $\Delta _{\mathscr E}$ as a cellular homotopy from 
$l_{\overline{E^{\prime }}}*l_{\overline{E^{\prime \prime }}}$ 
to $l_{\overline{E}}$, then the Bott-Chern form of this homotopy is 
$\CCGS _1(\mathscr E)$.
Hence we have 
$$
[(l_{\overline{E^{\prime }}}*l_{\overline{E^{\prime \prime }}}, 0)]=
[(l_{\overline{E}}, \CCGS _1(\mathscr E))]
$$
in $\widehat{\pi }_1(\vert \widehat{S}(X)\vert , \CCGS )$.
This tells that $(\overline{E}, \omega )\mapsto 
(l_{\overline{E}}, -\omega )$ gives rise to the homomorphism of groups 
$$
\widehat{\alpha }:\widehat{\mathscr K}_0(X)\to \widehat{\pi }_1
(\vert \widehat{S}(X)\vert , \CCGS )=\widehat{K}_0(X).
$$

Consider the following commutative diagram:
$$
\begin{CD}
K_1(X) @>{\rho }>> \widetilde{\mathscr A}_1(X) @>{a}>> 
\widehat{\mathscr K}_0(X) @>>> K_0(X) @>>> 0  \\
@VV{-\text{id}}V @VV{-\text{id}}V @VV{\widehat{\alpha }}V @VV{\text{id}}V \\
K_1(X) @>{\rho }>> \widetilde{\mathscr A}_1(X) @>{a}>> \widehat{K}_0(X) @>>> 
K_0(X) @>>> 0.
\end{CD}
$$
The upper sequence is exact by Thm.6.2 in \cite{gilsoul} and the lower 
one is exact by Thm.3.6.
Hence $\widehat{\alpha }$ is bijective by the five lemma.
\qed

\vskip 1pc
The Bott-Chern form of a pointed cellular map 
$f:S^{n+1}\to \vert \widehat{S}(X)\vert $ is defined as 
$$
\CCGS _n(f)=\CCGS (f_*([S^{n+1}]))\in \mathscr A_n(X),
$$
where $[S^{n+1}]\in C_{n+1}(S^{n+1})$ is the fundamental chain of 
$S^{n+1}$.
The Bott-Chern form $\CCGS _{n+1}(H)$ of a cellular homotopy 
$H:S^{n+1}\times I/\{*\}\times I\to \vert \widehat{S}(X)\vert $ 
can be defined in the same way.

Let 
$$
\CCGS _n:\widehat{K}_n(X)\to \mathscr A_n(X)
$$
be the map defined in \S 2.2.
That is to say, $\CCGS _n$ is defined as 
$$
\CCGS _n([(f, \omega )])=\CCGS _n(f)-d_{\mathscr A}\omega .
$$
Let us call it the {\it Chern form map}.
Applying Cor.2.5 and Cor.2.6 to the present situation, 
we have the following corollaries:

\vskip 1pc
\begin{cor} \ \ 
There is an exact sequence 
$$
K_{n+1}(X)\overset{\rho }{\to }H_{n+1}(\mathscr A_*(X), d_{\mathscr A})
\to \widehat{K}_n(X)\overset{(\zeta , \CCGS _n)}{\longrightarrow }
K_n(X)\oplus \KER d_{\mathscr A}\overset{cl}{\to }
H_n(\mathscr A_*(X), d_{\mathscr A})\to 0,
$$
where $cl$ is defined as $cl(x, \omega )=\rho (x)-[\omega ]$.
\end{cor}

\vskip 1pc
\begin{cor} \ \ 
We define a subgroup $KM_n(X)\subset \widehat{K}_n(X)$ by the kernel of 
the Chern form map, that is, 
$$
KM_n(X)=\KER \left(\CCGS _n:\widehat{K}_n(X)\to \mathscr A_n(X)\right).
$$
Then there is a long exact sequence 
$$
\cdots \to K_{n+1}(X)\overset{\rho }{\to }\underset{p}{\oplus }
H^{2p-n-1}(\mathscr A^*(X, p), d_{\mathscr A})\to KM_n(X)\to K_n(X)
\to \cdots . 
$$
\end{cor}

\vskip 1pc
Cor.3.9 tells that the group $KM_*(X)$, which is supposed to be 
the homotopy groups of the homotopy fiber of the regulator map, 
can be defined in terms of modified homotopy groups 
developed in the last section.

We finish this subsection by defining a pull back morphism.
Let $\varphi :X\to Y$ be a morphism of proper arithmetic varieties.
Then there is a commutative diagram 
$$
\begin{CD}
C_*(\vert \widehat{S}(Y)\vert ) @>{\CCGS }>> \mathscr A_*(Y)[1]  \\
@VV{\varphi ^*}V @VV{\varphi ^*}V  \\
C_*(\vert \widehat{S}(X)\vert ) @>{\CCGS }>> \mathscr A_*(X)[1].
\end{CD}
$$
Hence a homomorphism 
$$
\widehat{\varphi }^*:\widehat{K}_n(Y)\to \widehat{K}_n(X)
$$
is defined by 
$\widehat{\varphi }^*([(f, \omega )])=[(\varphi ^*f, \varphi ^*\omega )]$. 
In the case of $n=0$, the isomorphism $\widehat{\alpha }$ in Thm.3.7 
identifies the above $\widehat{\varphi }^*$ with the pull back morphism 
on $\widehat{\mathscr K}_0(X)$ defined in \cite{gilsoul}.
It is obvious that the pull back morphism $\widehat{\varphi }^*$ 
commutes with the Chern form map $\CCGS _n$.

\vskip 1pc
\subsection{Arakelov $K$-theory} \ \ 
Let $M$ be a compact algebraic K\"{a}hler manifold and $h_M$ 
a K\"{a}hler metric on $M$. 
Let $\mathscr H^n_{\mathbb R}(M)$ be the space of real harmonic forms 
on $M$ with respect to $h_M$ and $\mathscr H^{p,q}(M)$ the space of 
harmonic forms of type $(p, q)$. 
We set 
$$
\mathscr H^n_{\mathscr D}(M, p)=\begin{cases}
\mathscr H^{n-1}_{\mathbb R}(M)(p-1)\cap 
\underset{p^{\prime }<p, q^{\prime }<p}
{\underset{p^{\prime }+q^{\prime }=n-1}{\oplus }}
\mathscr H^{p^{\prime }, q^{\prime }}(M), \ & n<2p, \\
\mathscr H^{2p}_{\mathbb R}(M)(p)\cap \mathscr H^{p,p}(M), \ & n=2p.
\end{cases}
$$
The short exact sequence 
$$
0\to F^pH^{n-1}(M, \mathbb C)\to H^{n-1}(M, \mathbb R(p-1))\to 
H^n_{\mathscr D}(M, \mathbb R(p))\to 0
$$
for $n<2p$ and the short exact sequence 
$$
0\to H_{\mathscr D}^{2p}(M, \mathbb R(p))\to F^pH^{2p}(M, \mathbb C)
\to H^{2p}(M, \mathbb R(p-1))\to 0
$$
yield the isomorphism  
$$
H^n_{\mathscr D}(M, \mathbb R(p))\simeq 
\mathscr H^n_{\mathscr D}(M, p)
$$
for $n\leq 2p$.
Let $\mathscr H^n_{\mathscr A}(M, p)$ denote the image of 
$\mathscr H^n_{\mathscr D}(M, p)$ by the renormalization operator $\Theta $.
In other words, $\mathscr H^n_{\mathscr D}(M, p)$ 
(resp.~$\mathscr H^n_{\mathscr A}(M, p)$) is the space of harmonic forms 
in $\mathscr D^n(M, p)$ (resp.~$\mathscr A^n(M, p)$). 

Let us return to the arithmetic situation.
An {\it Arakelov variety} is a pair $\overline{X}=(X, h_X)$ 
of an arithmetic variety $X$ and an $F_{\infty }$-invariant 
K\"{a}hler metric $h_X$ on $X(\mathbb C)$. 
We now assume that $X$ is proper over $\mathbb Z$.
Let $\mathscr H_n(X)$ denote the space of harmonic 
forms with respect to $h_X$ in $\mathscr A_n(X)$, that is, 
$$
\mathscr H_n(X)=\begin{cases} \underset{p}{\oplus }
\mathscr H^{2p-n}_{\mathscr A}(X(\mathbb C), p)^{F_{\infty }^*=(-1)^{p-1}}, 
&n\geq 1, \\
\underset{p}{\oplus }
\mathscr H^{2p}_{\mathscr A}(X(\mathbb C), p)^{F_{\infty }^*=(-1)^p},
&n=0. \end{cases}
$$
Then there is an isomorphism $H_n(\mathscr A_*(X), d_{\mathscr A})
\simeq \mathscr H_n(X)$ and this yields the following:

\vskip 1pc
\begin{prop} \ \ 
there is an orthogonal decomposition 
$$
\KER d_{\mathscr A}=\IIm d_{\mathscr A}\oplus \mathscr H_n(X)
$$
in $\mathscr A_n(X)$.
\end{prop}

\vskip 1pc
\begin{defn} \ \ 
The subgroup $(\CCGS _n)^{-1}\left(\mathscr H_n(X)\right)$ 
of $\widehat{K}_n(X)$ is denoted by $K_n(\overline{X})$ and 
called the {\it $n$-th Arakelov $K$-group} of 
$\overline{X}=(X, h_X)$.
\end{defn}

\vskip 1pc
\begin{thm} \ \ 
There is an exact sequence 
$$
K_{n+1}(X)\to \underset{p}{\oplus }
H^{2p-n-1}(\mathscr A^*(X, p), d_{\mathscr A})\to K_n(\overline{X})
\to K_n(X)\to 0.
$$
\end{thm}

{\it Proof}: \ 
This exact sequence can be deduced from the fact that 
$[(0, \omega )]\in K_n(\overline{X})$ if and only if 
$d_{\mathscr A}\omega =0$, which follows from Prop.3.10. 
\qed

\vskip 1pc
K\"{u}nnemann has constructed a section of the inclusion 
from the Arakelov Chow group to the arithmetic Chow group 
in \cite{kunn}.
We now adapt his method to the inclusion \linebreak
$K_n(\overline{X})\hookrightarrow \widehat{K}_n(X)$ to get 
a section of it.
Let $\mathscr H:\mathscr A_n(X)\to \mathscr H_n(X)$ be the orthogonal 
projection with respect to the $L_2$-inner product.
Let $(f, \omega )$ be a pair of a pointed cellular map 
$f:S^{n+1}\to \vert \widehat{S}(X)\vert $ 
and $\omega \in \widetilde{A}_{n+1}(X)$. 
Then we can take $\omega _{\sharp }\in \widetilde{A}_{n+1}(X)$ 
such that $\CCGS _n(f)-d_{\mathscr A}\omega _{\sharp }$ is 
harmonic and $\mathscr H(\omega _{\sharp })=\mathscr H(\omega )$.
Existence and uniqueness of $\omega _{\sharp }$ follow from 
Prop.3.10.

If $(f, \omega )$ is homotopy equivalent to 
$(f^{\prime }, \omega ^{\prime })$, then we have 
$$
\CCGS _n(f^{\prime })-d_{\mathscr A}(\omega ^{\prime }
-\omega _{\sharp }-\omega )=\CCGS _n(f)-d_{\mathscr A}\omega _{\sharp }
$$
and $\mathscr H(\omega ^{\prime }-\omega _{\sharp }-\omega )=
\mathscr H(\omega ^{\prime })$.
Hence $\omega ^{\prime }_{\sharp }=\omega ^{\prime }-\omega _{\sharp }
-\omega $, therefore $(f, \omega _{\sharp })$ is homotopy equivalent to 
$(f^{\prime }, \omega ^{\prime }_{\sharp })$.
Hence we can define a section of the inclusion 
$$
\sigma :\widehat{K}_n(X)\to K_n(\overline{X})
$$
by $\sigma ([(f, \omega )])=[(f, \omega _{\sharp })]$. 
The map $\sigma $ is called the {\it harmonic projection} 
of $\widehat{K}_n(X)$.

\vskip 2pc
\section{A product formula for higher Bott-Chern forms}
\vskip 1pc

\subsection{A product formula} \ \ 
We begin this section by recalling the multiplicative structure on 
$\mathscr D^n(X, p)$ \cite{burgos1}.
Let $M$ be a compact complex algebraic manifold. 
A homomorphism  
$$
\bullet :\mathscr D^n(M, p)\otimes \mathscr D^m(M, q)
\to \mathscr D^{m+n}(M, p+q)
$$
is defined as 
$$
x\bullet y=(-1)^n(\partial x^{(p-1,n-p)}-\overline{\partial }x^{(n-p,p-1)})
\wedge y+x\wedge (\partial y^{(q-1,m-q)}-\overline{\partial }y^{(m-q,q-1)})
$$
if $n<2p$ and $m<2q$ and $x\bullet y=x\wedge y$ if $n=2p$ or $m=2q$.
Here $x^{(\alpha ,\beta )}$ is the $(\alpha , \beta )$-\linebreak 
component of 
the form $x$.
Then it satisfies 
$d_{\mathscr D}(x\bullet y)=d_{\mathscr D}x\bullet y+(-1)^nx\bullet 
d_{\mathscr D}y$ and $x\bullet y=(-1)^{nm}y\bullet x$. 
Moreover, it induces the product in the real Deligne cohomology, 
which is defined in \cite{esnault}.

The higher Bott-Chern forms are not compatible 
with products, that is, $\CC _{n+m}(\mathcal F\otimes \mathcal G)$ 
is not equal to $\CC _n(\mathcal F)\bullet \CC _m(\mathcal G)$ for 
an exact metrized $n$-cube $\mathcal F$ and an exact metrized $m$-cube 
$\mathcal G$.
But since the higher Chern character  
$K_n(M)\to H^{2p-n}_{\mathscr D}(M, p)$ respects the products, 
it is quite natural to expect that the difference 
$\CC _{n+m}(\mathcal F\otimes \mathcal G)-\CC _n(\mathcal F)\bullet 
\CC _m(\mathcal G)$ is expressed in terms of exact forms.

Let us define another operation on $\mathscr D^n(M, p)$.
For integers $i$ and $j$ satisfying $1\leq i\leq n$ and $1\leq j\leq m$, 
we set
$$
a_{i,j}^{n,m}=1-2\sbinom{n+m}{n}^{-1}\sum_{\alpha =0}^{i-1}
\sbinom{n+m-i-j+1}{n-\alpha }\sbinom{i+j-1}{\alpha },
$$
where $\sbinom{a}{b}=\frac{(a+b)!}{a!b!}$.
When $b<0$ or $a<b$, $\sbinom{a}{b}$ is assumed to be zero.

\vskip 1pc
\begin{lem} \ \ 
We have $a^{n,m}_{i,j}=-a^{n,m}_{n-i+1,m-j+1}$ and 
$a^{n,m}_{i,j}=-a^{m,n}_{j,i}$.
\end{lem}

{\it Proof}: \ 
Let us recall the following formula on binomial coefficients: 
$$
\sum_{\alpha =0}^a\sbinom{b}{a-\alpha }\sbinom{c}{\alpha }=\sbinom{b+c}{a}.
$$
Then we have 
{\allowdisplaybreaks
\begin{align*}
a^{n,m}_{i,j}+a^{n,m}_{n-i+1,m-j+1}&=2-2\sbinom{n+m}{n}^{-1}
\sum_{\alpha =0}^{i-1}
\sbinom{n+m-i-j+1}{n-\alpha }\sbinom{i+j-1}{\alpha }   \\*
&\hskip 1pc -2\sbinom{n+m}{n}^{-1}\sum_{\alpha =0}^{n-i}
\sbinom{i+j-1}{n-\alpha }\sbinom{n+m-i-j+1}{\alpha }   \\
&=2-2\sbinom{n+m}{n}^{-1}\sum_{\alpha =0}^n
\sbinom{n+m-i-j+1}{n-\alpha }\sbinom{i+j-1}{\alpha }  \\
&=0.
\end{align*}}
As for the second equality, we have 
{\allowdisplaybreaks
\begin{align*}
a^{m,n}_{j,i}&=1-2\sbinom{n+m}{m}^{-1}\sum_{\alpha =0}^{j-1}
\sbinom{m+n-j-i+1}{m-\alpha }\sbinom{j+i-1}{\alpha }  \\
&=-1+2\sbinom{n+m}{m}^{-1}\sum_{\alpha =j}^m
\sbinom{m+n-j-i+1}{m-\alpha }\sbinom{j+i-1}{\alpha }  \\
&=-1+2\sbinom{n+m}{m}^{-1}\sum_{\alpha =j}^m
\sbinom{m+n-j-i+1}{n-j-i+1+\alpha }\sbinom{j+i-1}{j+i-1-\alpha }.  \\
\intertext{If we put $\beta =i+j-1-\alpha $, then}
a^{m,n}_{j,i}&=-1+2\sbinom{n+m}{n}^{-1}\sum_{\beta =0}^{i-1}
\sbinom{m+n-i-j+1}{n-\beta }\sbinom{i+j-1}{\beta }  \\
&=-a^{n,m}_{i,j}, 
\end{align*}}
which completes the proof.
\qed

\vskip 1pc
For $x\in \mathscr D^{2p-n}(M, p)$ and 
$y\in \mathscr D^{2q-m}(M, q)$ with $n, m\geq 1$, $x\vartriangle y$ 
is defined as 
$$
x\vartriangle y=\underset{1\leq j\leq m}{\underset{1\leq i\leq n}{\sum }}
a_{i,j}^{n,m}x^{(p-n+i-1, p-i)}\wedge y^{(q-m+j-1, q-j)}.
$$
If $n=0$ or $m=0$, then $x\vartriangle y$ is defined to be zero.
The first claim of Lem.4.1 implies that $x\vartriangle y$ is in 
$\mathscr D^{2(p+q)-n-m-1}(M, p+q)$ and the second claim implies  
$x\vartriangle y=(-1)^{nm+n+m}y\vartriangle x$. 

\vskip 1pc
\begin{thm} \ \ 
If $\mathcal F$ is an exact metrized $n$-cube on $M$ and $\mathcal G$ 
is an exact metrized $m$-cube on $M$, then we have 
\begin{align*}
\CC _{n+m}(\mathcal F\otimes \mathcal G)&-\CC _n(\mathcal F)
\bullet \CC _m(\mathcal G)=(-1)^{n+1}d_{\mathscr D}
(\CC _n(\mathcal F)\vartriangle \CC _m(\mathcal G))  \\
&+(-1)^n\CC _{n-1}(\partial \mathcal F)\vartriangle \CC _m(\mathcal G)
-\CC _n(\mathcal F)\vartriangle \CC _{m-1}(\partial \mathcal G).
\end{align*}
\end{thm}

\vskip 1pc
\subsection{Proof of Thm.4.2} \ \ 
Let us first prepare some notations.
For differential forms $u_1, \cdots , u_n$ on $M$, let 
$(u_1, \cdots , u_n)^{(\alpha , \beta )}$ be the 
$(\alpha , \beta )$-component of $du_1\wedge \cdots \wedge du_n$.
When $u_i$ is a $(p_i, p_i)$-form, let 
$$
(u_1, \cdots , u_n)^{(i)}=\sum_p(u_1, \cdots , u_n)^{(p+i,p+n-i)}
$$
and 
$$
S_n^i(u_1, \cdots , u_n)=(i-1)!(n-i)!\sum_{\alpha =1}^n(-1)^{\alpha +1}
u_{\alpha }(u_1, \cdots , \widehat{u_{\alpha }}, \cdots , u_n)^{(i-1)}.
$$
Then $S_n^i(u_1, \cdots , u_n)\in \mathscr D_n(M)$ if 
$u_i\in \mathscr D_1(M)$.
We note that $S_n^i(\log |z_1|^2, \cdots , \log |z_n|^2)$ is the same form 
as $S_n^i$ defined in \S 1.5.

\vskip 1pc
\begin{lem} \ \ 
If $u_i$ is a $(p_i, p_i)$-form on $M$, then we have 
\begin{align*}
\partial S_n^i(u_1, \cdots , u_n)&=i!(n-i)!(u_1, \cdots , u_n)^{(i)} \\
&\hskip 1pc +(n-i)\sum_{\alpha =1}^n(-1)^{\alpha }
\partial \overline{\partial }
u_{\alpha }S_{n-1}^i(u_1, \cdots , \widehat{u_{\alpha }}, \cdots , u_n)
\end{align*}
and 
\begin{align*}
\overline{\partial }S_n^i(u_1, \cdots , u_n)
&=(i-1)!(n-i+1)!(u_1, \cdots , u_n)^{(i-1)} \\
&\hskip 1pc -(i-1)\sum_{\alpha =1}^n(-1)^{\alpha }
\partial \overline{\partial }
u_{\alpha }S_{n-1}^{i-1}(u_1, \cdots , \widehat{u_{\alpha }}, \cdots , u_n).
\end{align*}
\end{lem} 

{\it Proof}: \ 
We have 
\begin{small}
{\allowdisplaybreaks
\begin{align*}
\partial &S^i_n(u_1, \cdots , u_n)=(i-1)!(n-i)!\sum_{\alpha =1}^n
(-1)^{\alpha +1}\partial \left(u_{\alpha }(u_1, \cdots , 
\widehat{u_{\alpha }}, \cdots u_n)^{(i-1)}\right)  \\
&=(i-1)!(n-i)!\sum_{\alpha =1}^n(-1)^{\alpha +1}\partial u_{\alpha }
(\cdots \widehat{u_{\alpha }} \cdots )^{(i-1)}  \\
&\hskip 1pc +(i-1)!(n-i)!\sum_{\alpha =1}^n(-1)^{\alpha +1}u_{\alpha }
\left(\sum_{\beta <\alpha }(-1)^{\beta -1}\partial \overline{\partial }
u_{\beta }(\cdots \widehat{u_{\beta }}\cdots \widehat{u_{\alpha }}
\cdots )^{(i-1)}\right.  \\*
&\hskip 12pc +\left.\sum_{\alpha <\beta }(-1)^{\beta -1}\partial 
\overline{\partial }u_{\beta }(\cdots \widehat{u_{\alpha }}\cdots 
\widehat{u_{\beta }}\cdots )^{(i-1)}\right)  \\
&=i!(n-i)!(u_1, \cdots , u_n)^{(i)}+(i-1)!(n-i)!\sum_{\beta =1}^n
(-1)^{\beta }\partial \overline{\partial }u_{\beta }\left(
\sum_{\alpha <\beta }(-1)^{\alpha +1}u_{\alpha }(\cdots 
\widehat{u_{\alpha }}\cdots \widehat{u_{\beta }}\cdots )^{(i-1)}\right. \\*
&\hskip 12pc +\left.\sum_{\beta <\alpha }(-1)^{\alpha }u_{\alpha }
(\cdots \widehat{u_{\beta }}\cdots \widehat{u_{\alpha }}\cdots )^{(i-1)}
\right)  \\*
&=i!(n-i)!(u_1, \cdots , u_n)^{(i)}+(n-i)\sum_{\beta =1}^n(-1)^{\beta }
\partial \overline{\partial }u_{\beta }S_{n-1}^i
(u_1, \cdots \widehat{u_{\beta }}\cdots , u_n).
\end{align*}}
\end{small}
The second identity can be proved in a similar way.
\qed

\vskip 1pc
\begin{lem} \ \ 
For $(p_i, p_i)$-forms $u_i$ and $(q_j, q_j)$-forms $v_j$, we have 
\begin{align*}
S_{n+m}^k&(u_1, \cdots , u_n, v_1, \cdots , v_m)=
\sum_{i=1}^k\frac{(k-1)!(n+m-k)!}{(n-i)!(i-1)!}
S_n^i(u_1, \cdots , u_n)\wedge (v_1, \cdots , v_m)^{(k-i)}  \\
&+(-1)^n\sum_{j=1}^k\frac{(k-1)!(n+m-k)!}{(m-j)!(j-1)!}
(u_1, \cdots , u_n)^{(k-j)}\wedge S_m^j(v_1, \cdots , v_m).
\end{align*}
Hence we have 
\begin{align*}
\sum_{k=1}^{n+m}&(-1)^kS_{n+m}^k(u_1, \cdots , u_n, v_1, \cdots , v_m) \\
&=\underset{0\leq 1\leq m}{\sum_{1\leq i\leq n}}(-1)^{i+j}
\frac{(n+m-i-j)!(i+j-1)!}{(n-i)!(i-1)!}S_n^i(u_1, \cdots , u_n)\wedge 
(v_1, \cdots , v_m)^{(j)} \\
&\hskip 1pc +\underset{1\leq j\leq m}{\sum_{0\leq i\leq n}}(-1)^{n+i+j}
\frac{(n+m-i-j)!(i+j-1)!}{(m-j)!(j-1)!}(v_1, \cdots v_n)^{(i)}\wedge 
S_m^j(v_1, \cdots , v_m).
\end{align*}
\end{lem}

{\it Proof}: \ 
We have 
\begin{small}
{\allowdisplaybreaks
\begin{align*}
S&_{n+m}^k(u_1, \cdots , u_n, v_1, \cdots , v_m)  \\
&=(k-1)!(n+m-k)!\sum_{\alpha =1}^n(-1)^{\alpha +1}u_{\alpha }
(u_1, \cdots \widehat{u_{\alpha }} \cdots u_n, v_1, \cdots , v_m)^{(k-1)} \\
&\hskip 1pc +(-1)^n(k-1)!(n+m-k)!\sum_{\beta =1}^n(-1)^{\beta +1}
v_{\beta }(u_1, \cdots, u_n, v_1, \cdots \widehat{v_{\beta }} \cdots , 
v_m)^{(k-1)} \\
&=(k-1)!(n+m-k)!\sum_{\alpha =1}^n(-1)^{\alpha +1}u_{\alpha }
\left(\sum_{i=1}^k(u_1, \cdots \widehat{u_{\alpha }} \cdots u_n)^{(i-1)}
\wedge (v_1, \cdots , v_m)^{(k-i)}\right) \\
&\hskip 1pc +(-1)^n(k-1)!(n+m-k)!\sum_{\beta =1}^n(-1)^{\beta +1}
v_{\beta }\left(\sum_{j=1}^k(u_1, \cdots, u_n)^{(k-j)}\wedge 
(v_1, \cdots \widehat{v_{\beta }} \cdots , v_m)^{(j-1)}\right) \\
&=\sum_{i=1}^k\sfrac{(k-1)!(n+m-k)!}{(i-1)!(n-i)!}
S_n^i(u_1, \cdots , u_n)\wedge (v_1, \cdots , v_m)^{(k-i)} \\
&\hskip 1pc +(-1)^n\sum_{j=1}^k\sfrac{(k-1)!(n+m-k)!}{(j-1)!(m-j)!}
(u_1, \cdots, u_n)^{(k-j)}\wedge S_m^j(v_1, \cdots , v_m),
\end{align*}}
\end{small}
which completes the proof.
\qed

\vskip 1pc
If we assume that $u_i$ and $v_j$ are in $\mathscr D_1(M)$, then 
by Lem.4.3 we have 
{\allowdisplaybreaks
\begin{align*}
d_{\mathscr D}&\left(\sum_{i=1}^n(-1)^iS_n^i(u_1, \cdots , u_n)
\vartriangle \sum_{j=1}^m(-1)^jS_m^j(v_1, \cdots v_m)\right)  \\
&=-\underset{1\leq j\leq m}{\sum_{1\leq i\leq n}}(-1)^{i+j}a^{n,m}_{i,j}
dS_n^i(u_1, \cdots , u_n)\wedge S_m^j(v_1, \cdots v_m)  \\*
&\hskip 1pc +(-1)^n\underset{1\leq j\leq m}{\sum_{1\leq i\leq n}}(-1)^{i+j}
a^{n,m}_{i,j}S_n^i(u_1, \cdots , u_n)\wedge dS_m^j(v_1, \cdots v_m)  \\
&=-\underset{1\leq j\leq m}{\sum_{0\leq i\leq n}}(-1)^{i+j}i!(n-i)!
(a^{n,m}_{i,j}-a^{n,m}_{i+1,j})(u_1, \cdots , u_n)^{(i)}\wedge 
S_m^j(v_1, \cdots v_m)  \\*
&\hskip 1pc -\underset{1\leq j\leq m}{\sum_{1\leq i\leq n-1}}
(-1)^{i+j}((n-i)a^{n,m}_{i,j}+ia^{n,m}_{i+1,j})  \\*
&\hskip 5pc \times \left(\sum_{\alpha =1}^n(-1)^{\alpha }
\partial \overline{\partial }u_{\alpha }
S_{n-1}^i(u_1, \cdots , \widehat{u_{\alpha }}, \cdots , u_n)
\wedge S_m^j(v_1, \cdots , v_m)\right)  \\*
&\hskip 1pc +(-1)^n\underset{0\leq j\leq m}{\sum_{1\leq i\leq n}}
(-1)^{i+j}j!(m-j)!(a^{n,m}_{i,j}-a^{n,m}_{i,j+1})
S_n^i(u_1, \cdots , u_n)\wedge (v_1, \cdots v_m)^{(j)}  \\*
&\hskip 1pc +(-1)^n\underset{1\leq j\leq m-1}{\sum_{1\leq i\leq n}}
(-1)^{i+j}((m-j)a^{n,m}_{i,j}+ja^{n,m}_{i,j+1})  \\*
&\hskip 5pc \times \left(S_n^i(u_1, \cdots , u_n)\wedge \sum_{\beta =1}^m
(-1)^{\beta }\partial \overline{\partial }v_{\beta }
S_{m-1}^j(v_1, \cdots , \widehat{v_{\beta }}, \cdots , v_m)\right).
\end{align*}}
Let us compute the coefficients of the above expression.
Since $\sbinom{n+m-i-j+1}{n-\alpha }=\sbinom{n+m-i-j}{n-\alpha }+
\sbinom{n+m-i-j}{n-1-\alpha }$ and $\sbinom{i+j}{\alpha }=
\sbinom{i+j-1}{\alpha }+\sbinom{i+j-1}{\alpha -1}$, we have 
{\allowdisplaybreaks
\begin{align*}
a^{n,m}_{i,j}-a^{n,m}_{i+1,j}&=1-2\sbinom{n+m}{n}^{-1}
\sum_{\alpha =0}^{i-1}\sbinom{n+m-i-j+1}{n-\alpha }\sbinom{i+j-1}{\alpha } \\*
&\hskip 1pc -1+2\sbinom{n+m}{n}^{-1}
\sum_{\alpha =0}^{i}\sbinom{n+m-i-j}{n-\alpha }\sbinom{i+j}{\alpha }  \\
&=2\sbinom{n+m}{n}^{-1}\sbinom{n+m-i-j}{n-i}\sbinom{i+j-1}{i}
\end{align*}}
for $1\leq i\leq n-1$ and 
{\allowdisplaybreaks
\begin{align*}
a^{n,m}_{i,j}-a^{n,m}_{i,j+1}&=1-2\sbinom{n+m}{n}^{-1}
\sum_{\alpha =0}^{i-1}\sbinom{n+m-i-j+1}{n-\alpha }\sbinom{i+j-1}{\alpha } \\*
&\hskip 1pc -1+2\sbinom{n+m}{n}^{-1}
\sum_{\alpha =0}^{i-1}\sbinom{n+m-i-j}{n-\alpha }\sbinom{i+j}{\alpha }  \\
&=-2\sbinom{n+m}{n}^{-1}\sbinom{n+m-i-j}{n-i}\sbinom{i+j-1}{i-1}  \\
&=-2\sbinom{n+m}{n}^{-1}\sbinom{n+m-i-j}{m-j}\sbinom{i+j-1}{j}  \\
\end{align*}}
for $1\leq j\leq m-1$. 
It follows from the definition of $a^{n,m}_{i,j}$ and Lem.4.1 that 
{\allowdisplaybreaks
\begin{align*}
a^{n,m}_{1,j}&=1-2\sbinom{n+m}{n}^{-1}\sbinom{n+m-j}{n},  \\
a^{n,m}_{n,j}&=-1+2\sbinom{n+m}{n}^{-1}\sbinom{n+j-1}{n},  \\
a^{n,m}_{i,1}&=-1+2\sbinom{n+m}{m}^{-1}\sbinom{n+m-i}{m},  \\
a^{n,m}_{i,m}&=1-2\sbinom{n+m}{m}^{-1}\sbinom{m+i-1}{m}.
\end{align*}}
By Lem.A.1, we have 
{\allowdisplaybreaks
\begin{align*}
(n&-i)a^{n,m}_{i,j}+ia^{n,m}_{i+1,j}  \\
&=n-2\sbinom{n+m}{n}^{-1}\left((n-i)\sum_{\alpha =0}^{i-1}
\sbinom{n+m-i-j+1}{n-\alpha }\sbinom{i+j-1}{\alpha }+i\sum_{\alpha =0}^i
\sbinom{n+m-i-j}{n-\alpha }\sbinom{i+j}{\alpha }\right)  \\
&=n-2n\sbinom{n+m-1}{n-1}^{-1}\sum_{\alpha =0}^{i-1}
\sbinom{n+m-i-j}{n-1-\alpha }\sbinom{i+j-1}{\alpha }  \\
&=na^{n-1,m}_{i,j}
\end{align*}}
and 
{\allowdisplaybreaks
\begin{align*}
(m&-j)a^{n,m}_{i,j}+ja^{n,m}_{i,j+1}  \\
&=m-2\sbinom{n+m}{n}^{-1}\left((m-j)\sum_{\alpha =0}^{i-1}
\sbinom{n+m-i-j+1}{n-\alpha }\sbinom{i+j-1}{\alpha }+j\sum_{\alpha =0}^{i-1}
\sbinom{n+m-i-j}{n-\alpha }\sbinom{i+j}{\alpha }\right)  \\
&=m-2m\sbinom{n+m-1}{n}^{-1}\sum_{\alpha =0}^{i-1}\sbinom{n+m-i-j}{n-\alpha }
\sbinom{i+j-1}{\alpha }  \\
&=ma^{n,m-1}_{i,j}.
\end{align*}}
These computations imply that 
\begin{small}
{\allowdisplaybreaks 
\begin{align*}
d_{\mathscr D}&\left(\sum_{i=1}^n(-1)^iS_n^i(u_1, \cdots , u_n)\vartriangle 
\sum_{j=1}^n(-1)^jS_m^j(v_1, \cdots , v_m)\right)  \\
&=-2\sbinom{n+m}{n}^{-1}
\underset{1\leq j\leq m}{\sum_{0\leq i\leq n}}(-1)^{i+j}i!(n-i)!
\sbinom{n+m-i-j}{n-i}\sbinom{i+j-1}{i}(u_1, \cdots , u_n)^{(i)}\wedge 
S_m^j(v_1, \cdots , v_m)  \\*
&\hskip 1pc -n\underset{1\leq j\leq m}{\sum_{1\leq i\leq n-1}}(-1)^{i+j}
a^{n-1,m}_{i,j}\sum_{\alpha =1}^n(-1)^{\alpha }
\partial \overline{\partial }u_{\alpha }
S_{n-1}^i(u_1, \cdots , \widehat{u_{\alpha }}, \cdots , u_n)\wedge 
S_m^j(v_1, \cdots , v_m)  \\*
&\hskip 1pc -2(-1)^n\sbinom{n+m}{n}^{-1}
\underset{0\leq j\leq m}{\sum_{1\leq i\leq n}}(-1)^{i+j}j!(m-j)!
\sbinom{n+m-i-j}{m-j}\sbinom{i+j-1}{j}S_n^i(u_1, \cdots , u_n)\wedge 
(v_1, \cdots , v_m)^{(j)}  \\* 
&\hskip 1pc +(-1)^nm\underset{1\leq j\leq m-1}{\sum_{1\leq i\leq n}}
(-1)^{i+j}a^{n,m-1}_{i,j}S_n^i(u_1, \cdots , u_n)\wedge 
\sum_{\beta =1}^m(-1)^{\beta }\partial \overline{\partial }v_{\beta }
S_{m-1}^j(v_1, \cdots , \widehat{v_{\beta }}, \cdots , v_m) \\
&\hskip 1pc +\sum_{1\leq  j\leq m}(-1)^jn!(u_1, \cdots , u_n)^{(0)}\wedge 
S_m^j(v_1, \cdots , v_m)  \\*
&\hskip 1pc +\sum_{1\leq  j\leq m}(-1)^{n+j}n!(u_1, \cdots , u_n)^{(n)}
\wedge S_m^j(v_1, \cdots , v_m)  \\*
&\hskip 1pc +\sum_{1\leq i\leq n}(-1)^{n+i}m!S_n^i(u_1, \cdots , u_n)
\wedge (v_1, \cdots , v_m)^{(0)}  \\*
&\hskip 1pc +\sum_{1\leq i\leq n}(-1)^{n+m+i}m!S_n^i(u_1, \cdots , u_n)
\wedge (v_1, \cdots , v_m)^{(m)}  \\
&=-2\sbinom{n+m}{n}^{-1}\underset{1\leq j\leq m}{\sum_{0\leq i\leq n}}
(-1)^{i+j}\sfrac{(n+m-i-j)!(i+j-1)!}{(m-j)!(j-1)!}
(u_1, \cdots , u_n)^{(i)}\wedge S_m^j(v_1, \cdots , v_m)  \\*
&\hskip 1pc -n\underset{1\leq j\leq m}{\sum_{1\leq i\leq n-1}}(-1)^{i+j}
a^{n-1,m}_{i,j}\sum_{\alpha =1}^n(-1)^{\alpha }
\partial \overline{\partial }u_{\alpha }
S_{n-1}^i(u_1, \cdots , \widehat{u_{\alpha }}, \cdots , u_n)\wedge 
S_m^j(v_1, \cdots , v_m)  \\*
&\hskip 1pc -2(-1)^n\sbinom{n+m}{n}^{-1}
\underset{0\leq j\leq m}{\sum_{1\leq i\leq n}}(-1)^{i+j}
\sfrac{(n+m-i-j)!(i+j-1)!}{(n-i)!(i-1)!}S_n^i(u_1, \cdots , u_n)\wedge 
(v_1, \cdots , v_m)^{(j)}  \\* 
&\hskip 1pc +(-1)^nm\underset{1\leq j\leq m-1}{\sum_{1\leq i\leq n}}
(-1)^{i+j}a^{n,m-1}_{i,j}S_n^i(u_1, \cdots , u_n)\wedge 
\sum_{\beta =1}^m(-1)^{\beta }\partial \overline{\partial }v_{\beta }
S_{m-1}^j(v_1, \cdots , \widehat{v_{\beta }}, \cdots , v_m)  \\
&\hskip 1pc +\sum_{1\leq  j\leq m}(-1)^j
\left(\overline{\partial }S_n^1(u_1, \cdots , u_n)+(-1)^n\partial 
S_n^n(u_1, \cdots , u_n)\right)\wedge S_m^j(v_1, \cdots , v_m)  \\*
&\hskip 1pc +\sum_{1\leq i\leq n}(-1)^{n+i}S_n^i(u_1, \cdots , u_n)
\wedge \left(\overline{\partial }S_m^1(v_1, \cdots , v_m)
+(-1)^m\partial S_m^m(v_1, \cdots , v_m)\right).
\end{align*}}
\end{small}
Applying Lem.4.4 to the above, we have the following.

\vskip 1pc
\begin{prop} \ \ 
For $u_i\in \mathscr D_1(M)$ and $v_j\in \mathscr D_1(M)$, we have 
{\allowdisplaybreaks
\begin{small}
\begin{align*}
&d_{\mathscr D}\left(\sum_{i=1}^n(-1)^iS_n^i(u_1, \cdots , u_n)\vartriangle 
\sum_{j=1}^m(-1)^jS_m^j(v_1, \cdots , v_m)\right)  \\
&=2(-1)^{n+1}\sbinom{n+m}{n}^{-1}\sum_{k=1}^{n+m}(-1)^k
S_{n+m}^k(u_1, \cdots , u_n, v_1, \cdots , v_m)  \\
&\hskip 1pc -n\underset{1\leq j\leq m}{\sum_{1\leq i\leq n-1}}(-1)^{i+j}
a^{n-1,m}_{i,j}\sum_{\alpha =1}^n(-1)^{\alpha }
\partial \overline{\partial }u_{\alpha }
S_{n-1}^i(u_1, \cdots , \widehat{u_{\alpha }}, \cdots , u_n)\wedge 
S_m^j(v_1, \cdots , v_m)  \\*
&\hskip 1pc +(-1)^nm\underset{1\leq j\leq m-1}{\sum_{1\leq i\leq n}}
(-1)^{i+j}a^{n,m-1}_{i,j}S_n^i(u_1, \cdots , u_n)\wedge 
\sum_{\beta =1}^m(-1)^{\beta }\partial \overline{\partial }v_{\beta }
S_{m-1}^j(v_1, \cdots , \widehat{v_{\beta }}, \cdots , v_m)  \\
&\hskip 1pc +(-1)^n\left(\sum_{i=1}^n(-1)^iS_n^i(u_1, \cdots , u_n)\right)
\bullet \left(\sum_{j=1}^m(-1)^jS_m^j(v_1, \cdots , v_m)\right).
\end{align*}
\end{small}}
\end{prop}

\vskip 1pc
Let us return to the proof of Thm.4.2.
We may assume that $\mathcal F$ and $\mathcal G$ are emi-cubes.
For $s<t$, let $\pi _1:(\mathbb P^1)^t\to (\mathbb P^1)^s$ denote 
the projection defined by $(x_1, \cdots , x_t)\mapsto (x_1, \cdots , x_s)$ 
and let $\pi _2:(\mathbb P^1)^t\to (\mathbb P^1)^s$ denote the projection 
defined by $(x_1, \cdots , x_t)\mapsto (x_{t-s+1}, \cdots , x_t)$.
Then by Prop.4.5 we have 
\begin{small}
{\allowdisplaybreaks
\begin{align*}
&d_{\mathscr D}(\CC _n(\mathcal F)\vartriangle \CC _m(\mathcal G))  \\
&=\sfrac{(-1)^{n+m}}{4n!m!(2\pi \sqrt{-1})^{n+m}}
\int_{(\mathbb P^1)^{n+m}}\pi _1^*\CC _0(\TTR \mathcal F)\wedge 
\pi _2^*\CC _0(\TTR \mathcal G)\wedge d_{\mathscr D}\left(
(\sum_{i=1}^n(-1)^i\pi _1^*S_n^i)\vartriangle 
(\sum_{j=1}^m(-1)^j\pi _2^*S_m^j)\right)  \\
&=\sfrac{(-1)^{m+1}}{2(n+m)!(2\pi \sqrt{-1})^{n+m}}\int_{(\mathbb P^1)^{n+m}}
\CC _0(\TTR (\mathcal F\otimes \mathcal G))
\wedge \sum_{k=1}^{n+m}(-1)^kS_{n+m}^k  \\*
&\hskip 1pc +\sfrac{(-1)^{n+m+1}}{4(n-1)!m!(2\pi \sqrt{-1})^{n+m-1}}
\int_{(\mathbb P^1)^{n+m-1}}
\CC _0(\TTR (\partial \mathcal F\otimes \mathcal G))\wedge 
\underset{1\leq j\leq m}{\underset{1\leq i\leq n-1}{\sum }}
(-1)^{i+j}a^{n-1,m}_{i,j}\pi _1^*S_{n-1}^i\wedge \pi _2^*S_m^j \\*
&\hskip 1pc +\sfrac{(-1)^m}{4n!(m-1)!(2\pi \sqrt{-1})^{n+m-1}}
\int_{(\mathbb P^1)^{n+m-1}}
\CC _0(\TTR (\mathcal F\otimes \partial \mathcal G))\wedge 
\underset{1\leq j\leq m-1}{\underset{1\leq i\leq n}{\sum }}
(-1)^{i+j}a^{n,m-1}_{i,j}\pi _1^*S_n^i\wedge \pi _2^*S_{m-1}^j \\*
&\hskip 1pc +\sfrac{(-1)^m}{4n!m!(2\pi \sqrt{-1})^{n+m}}
\int_{(\mathbb P^1)^{n+m}}\CC _0(\TTR (\mathcal F\otimes \mathcal G))\wedge 
\left((\sum_{i=1}^n(-1)^i\pi _1^*S_n^i)\bullet 
(\sum_{j=1}^m(-1)^j\pi _2^*S_m^j)\right) \\
&=(-1)^{n+1}\CC _{n+m}(\mathcal F\otimes \mathcal G)
+\CC _{n-1}(\partial \mathcal F)\vartriangle \CC _m(\mathcal G)  \\*
&\hskip 1pc +(-1)^{n+1}\CC _n(\mathcal F)\vartriangle 
\CC _m(\partial \mathcal G)
+(-1)^n\CC _n(\mathcal F)\bullet \CC _m(\mathcal G),
\end{align*}}
\end{small}
which completes the proof.
\qed

\vskip 1pc
\subsection{Renormalization} \ \ 
In this subsection we renormalize the identity in Thm.4.2.
Let us define 
$$
\bullet :\mathscr A^n(M, p)\times \mathscr A^m(M, q)\to 
\mathscr A^{n+m}(M, p+q)
$$
as 
$$
x\bullet y=\frac{(-1)^n}{4\pi \sqrt{-1}}(\partial x^{(p-1,n-p)}
-\overline{\partial }x^{(n-p,p-1)})\wedge y+\frac{1}{4\pi \sqrt{-1}}
x\wedge (\partial y^{(q-1,m-q)}-\overline{\partial }y^{(m-q,q-1)})
$$
if $n<2p, m<2q$ and $x\bullet y=x\wedge y$ if $n=2p$ or $m=2q$. 
Then it follows that 
$$
\Theta (x\bullet y)=\Theta (x)\bullet \Theta (y) 
$$
for $x\in \mathscr D^n(M, p)$ and $y\in \mathscr D^m(M, q)$, 
where $\Theta :\mathscr D^*(M, *)\to \mathscr A^*(M, *)$ is the 
renormalization operator defined in \S 1.6.

For $n, m\geq 1$ let us define 
$$
\vartriangle :\mathscr A^{2p-n}(M, p)\times \mathscr A^{2q-m}(M, q)\to 
\mathscr A^{2p+2q-n-m-1}(M, p+q)
$$
by 
$$
x\vartriangle y=\frac{1}{4\pi \sqrt{-1}}\underset{1\leq j\leq m}
{\underset{1\leq i\leq n}{\sum }}a^{n,m}_{i,j}x^{(p-n+i-1,p-i)}
\wedge y^{(q-m+j-1,q-j)}.
$$
If $n=0$ or $m=0$, then $x\vartriangle y$ is defined to be zero.
Then it follows that 
$$
\Theta (x\vartriangle y)=\Theta (x)\vartriangle \Theta (y)
$$
for $x\in \mathscr D^{2p-n}(M, p)$ and $y\in \mathscr D^{2q-m}(M ,q)$.
Then we can rewrite the product formula as follows:

\vskip 1pc
\begin{prop} \ \ 
If $\mathcal F$ is an exact metrized $n$-cube on $M$ and $\mathcal G$ 
is an exact metrized $m$-cube on $M$, then we have 
\begin{align*}
\CCGS _{n+m}(\mathcal F\otimes \mathcal G)&-\CCGS _n(\mathcal F)
\bullet \CCGS _m(\mathcal G)=(-1)^{n+1}d_{\mathscr A}
(\CCGS _n(\mathcal F)\vartriangle \CCGS _m(\mathcal G))  \\
&+(-1)^n\CCGS _{n-1}(\partial \mathcal F)\vartriangle \CCGS _m(\mathcal G)
-\CCGS _n(\mathcal F)\vartriangle \CCGS _{m-1}(\partial \mathcal G).
\end{align*}
\end{prop}

\vskip 2pc
\section{Product}
\vskip 1pc

\subsection{Notations on bisimplicial sets} \ \ 
A {\it bisimplicial set} is a contravariant functor from 
the category of pairs of finite ordered sets to the category of sets.
The product $S\times T$ and the reduced product $S\wedge T$ 
of two simplicial sets $S, T$ are examples of bisimplicial sets.
For a bisimplicial set $S$, let $S_{n,m}=S([n], [m])$.
The topological realization of a bisimplicial set 
$S$ is defined as 
$$
\vert S\vert =\underset{n,m}{\coprod }S_{n,m}\times \Delta ^n\times 
\Delta ^m/\sim ,
$$
where the relation $\sim $ is given in a similar way 
to the case of a simplicial set.
Then $\vert S\vert $ is a CW-complex such that the set of 
$k$-cells of $\vert S\vert $ can be identified with the set of 
nondegenerate elements of $\underset{n+m=k}{\coprod }S_{n,m}$.

For a bisimplicial set $S$, $[n]\mapsto S_{n,n}$  becomes a simplicial 
set.
It is denoted by $\Delta (S)$ and called the {\it diagonal simplicial set} 
of $S$.
The topological realization $\vert \Delta (S)\vert $ is homeomorphic 
to $\vert S\vert $ and $\vert \Delta (S)\vert $ 
becomes a subdivision of $\vert S\vert $ as CW-complexes.
Hence the identity map $\vert S\vert \to \vert \Delta (S)\vert $ 
is cellular, although the inverse is not.

\vskip 1pc
\subsection{Product in higher $K$-theory} \ \ 
In this subsection we review the product in higher algebraic 
$K$-theory by means of $S$-construction \cite{wal}.
For a small exact category $\mathfrak A$, let $S_nS_m\mathfrak A$ be 
the set of functors 
$$
E:\AAR [n]\times \AAR [m]\to \mathfrak A, \ (i\leq j, \alpha \leq \beta )
\mapsto E_{(i,j)\times (\alpha ,\beta )}
$$
satisfying the following conditions:
\begin{enumerate}
\item \ $E_{(i,i)\times (\alpha ,\beta )}=0$ and 
$E_{(i,j)\times (\alpha ,\alpha )}=0$. 
\item \ For any $i\leq j\leq k$ and $\alpha \leq \beta $, 
$E_{(i,j)\times (\alpha ,\beta )}\to E_{(i,k)\times (\alpha ,\beta )}
\to E_{(j,k)\times (\alpha ,\beta )}$ is a short exact sequence 
of $\mathfrak A$.
\item \ For any $i\leq j$ and $\alpha \leq \beta \leq \gamma $, 
$E_{(i,j)\times (\alpha ,\beta )}\to E_{(i,j)\times (\alpha ,\gamma )}
\to E_{(i,j)\times (\beta ,\gamma )}$ is a short exact sequence 
of $\mathfrak A$.
\end{enumerate}
Then $([n], [m])\mapsto S_nS_m\mathfrak A$ is a bisimplicial set, 
which is denoted by $S^{(2)}\mathfrak A$. 
The natural identification $S_1S_m\mathfrak A=S_m\mathfrak A$ yields 
a morphism of bisimplicial sets 
$$
S^1\wedge S\mathfrak A\to S^{(2)}\mathfrak A, 
$$
and its adjoint map 
$\vert S\mathfrak A\vert \to \Omega \vert S^{(2)}\mathfrak A\vert $ 
turns out to be a homotopy equivalence.

When $\mathfrak A$ is equipped with a tensor product, 
a morphism of bisimplicial sets 
$$
m:S\mathfrak A\wedge S\mathfrak A\to S^{(2)}\mathfrak A 
$$
is defined by $m(E, F)_{(i,j)\times (\alpha , \beta )}=E_{i,j}
\otimes F_{\alpha , \beta }$.
This induces a pairing 
$$
m_*:\pi _{n+1}(\vert S\mathfrak A\vert )\times 
\pi _{m+1}(\vert S\mathfrak A\vert )\to 
\pi _{n+m+2}(\vert S^{(2)}\mathfrak A\vert ).
$$
Combining this pairing with the isomorphisms $K_n(\mathfrak A)\simeq 
\pi _{n+1}(\vert S\mathfrak A\vert )\simeq 
\pi _{n+2}(\vert S^{(2)}\mathfrak A\vert )$ yields the product 
in higher algebraic $K$-theory $K_*(\mathfrak A)$.

\vskip 1pc
\subsection{$G$-construction} \ \ 
In \cite{gilgray}, Gillet and Grayson constructed a simplicial set 
$G\mathfrak A$ associated with a small exact category $\mathfrak A$ 
that is homotopy equivalent to the loop space of the $S$-construction 
$S\mathfrak A$.
In this subsection we recall their construction.

Let $G_n\mathfrak A$ be the set of pairs $(E^+, E^-)$ of 
$E^+, E^-\in S_{n+1}\mathfrak A$ with 
$\partial _0E^+=\partial _0E^-$.
Let us define the boundary maps $\partial _k:G_n\mathfrak A\to 
G_{n-1}\mathfrak A$ by $(E^+, E^-)\mapsto (\partial _{k+1}E^+, 
\partial _{k+1}E^-)$ and the degeneracy maps $s_k:G_n\mathfrak A\to 
G_{n+1}\mathfrak A$ by $(E^+, E^-)\mapsto (s_{k+1}E^+, s_{k+1}E^-)$.
Then $(G_*\mathfrak A, \partial _k, s_k)$ becomes a simplicial set.
We fix $0=(0, 0)\in G_0\mathfrak A$ as the base point of $G\mathfrak A$.

Let $\iota _k$ denote the element of $\Delta [1]_n=\HOM ([n], [1])$ 
defined as
$$
\iota _k(i)=\begin{cases}0, &i<k,  \\ 1, &i\geq k. \end{cases}
$$
Then $\Delta [1]_n=\{\iota _0, \iota _1, \cdots , \iota _{n+1}\}$.
Let 
$$
\chi _n^{\pm }:\Delta [1]_n\times G_n\mathfrak A\to S_n\mathfrak A
$$
be maps defined as 
$$
\chi _n^{\pm }(\iota _k, (E^+, E^-))=\begin{cases}
\partial _0E^{\pm }, &k=0,   \\
(s_0)^{k-1}(\partial _1)^kE^{\pm },   &k\geq 1.
\end{cases}
$$
Then $\chi ^{\pm }=\{\chi ^{\pm }_n\}:\Delta (\Delta [1]\times 
G\mathfrak A)\to S\mathfrak A$ are morphisms of simplicial sets 
such that \linebreak 
$\chi ^{\pm }(\{0\}\times G\mathfrak A)=*$ and 
$\chi ^+|_{\{1\}\times G\mathfrak A}=\chi ^-|_{\{1\}\times G\mathfrak A}$.

Let $T^1$ be a simplicial set given by the following cocartesian square:
\begin{center}
\setlength{\unitlength}{1mm}
\begin{picture}(55,18)
   \put(2,7){$\{0\}\cup \{1\}$}
   \put(31,13){$\Delta [1]$}
   \put(31,1){$\Delta [1]$}
   \put(51,7){$T^1$.}
   \put(20,8.5){\vector(2,1){10}}
   \put(20,7.5){\vector(2,-1){10}}
   \put(40,13.5){\vector(2,-1){10}}
   \put(40,2.5){\vector(2,1){10}}
\end{picture}
\end{center}
We fix the image of $\{0\}\in \Delta [1]$ as the base point of $T^1$.
Since the topological realization of $\Delta [1]$ is the unit interval 
$I=[0, 1]$, the topological realization of $T^1$ is the barycentric 
subdivision of the circle $S^1=I/\partial I$.

\vskip 1pc
\begin{thm}\cite{gilgray} \ \ 
By gluing the morphisms $\chi ^{\pm }$, we obtain a morphism of 
simplicial sets
$$
\chi :\Delta (T^1\wedge G\mathfrak A)\to S\mathfrak A.
$$
The adjoint map $\vert G\mathfrak A\vert \to \Omega 
\vert S\mathfrak A\vert$ to $\vert \chi \vert $ is homotopy equivalent.
Therefore we have an isomorphism 
$$
\pi _i(\vert G\mathfrak A\vert ,0)\simeq K_i(\mathfrak A).
$$
\end{thm}

\vskip 1pc
We next introduce a description of the product in $K$-theory 
by means of $G$-construction, which is also given in \cite{gilgray}.
Let 
$$
G_nG_m\mathfrak A=\{(E^{++}, E^{+-}, E^{-+}, E^{--}); 
E^{\pm \pm }\in S_{n+1}S_{m+1}\mathfrak A, 
\partial _0E^{+\pm }=\partial _0E^{-\pm }, 
\partial ^{\prime }_0E^{\pm +}=\partial ^{\prime }_0E^{\pm -}\},
$$
where $\partial _0$ is the boundary map on the first factor of the 
bisimplicial set $S^{(2)}\mathfrak A$ and $\partial ^{\prime }_0$ 
is the boundary map on the second factor of $S^{(2)}\mathfrak A$.
Then $([n], [m])\mapsto G_nG_m\mathfrak A$ becomes a bisimplicial set, 
which is denoted by $G^{(2)}\mathfrak A$.
Let $R:G_n\mathfrak A\to G_0G_n\mathfrak A$ be the morphism defined by 
$R(E^+, E^-)=(E^+, E^-, 0, 0)$.
Then it is shown in \cite{gilgray} that $R$ induces a homotopy 
equivalent morphism $R:G\mathfrak A\to G^{(2)}\mathfrak A$.
Therefore we have an isomorphism 
$$
\pi _n(\vert G^{(2)}\mathfrak A\vert )\simeq K_n(\mathfrak A).
$$
A morphism of bisimplicial sets 
$$
m^G:G\mathfrak A\wedge G\mathfrak A\to G^{(2)}\mathfrak A
$$
is defined by $m^G(E, F)^{\pm \pm }=E^{\pm }\otimes F^{\pm }$ for 
$E=(E^+, E^-)\in G_n\mathfrak A$ and $F=(F^+, F^-)\in G_m\mathfrak A$.
Then the pairing 
$$
m^G_*:\pi _n(\vert G\mathfrak A\vert )\times \pi _m(\vert G\mathfrak A\vert )
\to \pi _{n+m}(\vert G^{(2)}\mathfrak A\vert )
$$
induces a product in $K_*(\mathfrak A)$, which agrees with 
the one defined by using the $S$-construction.

Finally, let us define a cube associated with an element of 
$G\mathfrak A$ or $G^{(2)}\mathfrak A$.
The morphism $\chi :\Delta (T^1\wedge G\mathfrak A)\to S\mathfrak A$ 
in Thm.5.1 yields a morphism of complexes 
$$
\CUB :C_*(\vert G\mathfrak A\vert )\overset{\chi _*}{\longrightarrow }
C_*(\vert S\mathfrak A\vert )[-1]\overset{\CUB }{\longrightarrow }
\CUB _*(\mathfrak A).
$$
The cube $\CUB (E)$ associated with $E=(E^+, E^-)\in G_n\mathfrak A$ 
is defined as the image by $\CUB $ of the element of 
$C_*(\vert G\mathfrak A\vert )$ represented by $E$. 
Since $\chi _*([E])=[E^+]-[E^-]$ in 
$C_{n+1}(\vert S\mathfrak A\vert )$, we have 
$$
\CUB (E)=\CUB (E^+)-\CUB (E^-).
$$
Let us define the cube associated with 
$E=(E^{\pm \pm })\in G_nG_m\mathfrak A$ as 
$$
\CUB (E)=\CUB (E^{++})-\CUB (E^{+-})-\CUB (E^{-+})+\CUB (E^{--}), 
$$
where $\CUB (E^{\pm \pm })$ is the image of the element of 
$E^{\pm \pm }\in S_{n+1}S_{m+1}\mathfrak A$ by the homomorphism
$$
S_{n+1}S_{m+1}\mathfrak A\to \CUB _n(S_{m+1}\mathfrak A)\to 
\CUB _n(\CUB _m(\mathfrak A))=\CUB _{n+m}(\mathfrak A).
$$
If $E=(E^{\pm \pm })$ is degenerate, then we can show in 
the same way as the proof of Lem.3.1 that the cube $\CUB (E^{\pm \pm })$ 
associated with $E$ is zero in $\CUB _{n+m}(\mathfrak A)$.
Hence $E=(E^{\pm \pm })\mapsto \CUB (E)$ induces a homomorphism 
of complexes 
$$
\CUB :C_*(\vert G^{(2)}\mathfrak A\vert )\to \CUB _*(\mathfrak A).
$$

\vskip 1pc
\begin{prop} \ \ 
We have the following commutative diagram: 
$$
\setlength{\unitlength}{1mm}
\begin{picture}(100,41)
  \put(1,35){$C_*(\vert G\mathfrak A\vert )\otimes 
          C_*(\vert G\mathfrak A\vert )$}
  \put(10,18){$C_*(\vert G^{(2)}\mathfrak A\vert )$}
  \put(11,2){$C_*(\vert G\mathfrak A\vert )$}
  \put(61,35){$\CUB _*(\mathfrak A)\otimes \CUB _*(\mathfrak A)$}
  \put(71,18){$\CUB _*(\mathfrak A)$}
  \put(71,2){$\CUB _*(\mathfrak A)$.}
  \put(20,33){\vector(0,-1){9}}
  \put(20,7){\vector(0,1){9}}
  \put(78,33){\vector(0,-1){9}}
  \put(78,7){\vector(0,1){9}}
  \put(41,36){\vector(1,0){17}}
  \put(41,19){\vector(1,0){17}}
  \put(41,3){\vector(1,0){17}}
  \put(42,37){$\scriptstyle{\CUB \otimes \CUB }$}
  \put(47,20){$\scriptstyle{\CUB }$}
  \put(47,4){$\scriptstyle{\CUB }$}
  \put(21,28){$\scriptstyle{m^G_*}$}
  \put(79,28){$\scriptstyle{\otimes }$}
  \put(21,11){$\scriptstyle{R_*}$}
  \put(79,11){$\scriptstyle{\IId }$}
\end{picture}
$$
\end{prop}

\vskip 1pc
\subsection{The pairing $\widehat{\mathscr K}_0\times \widehat{K}_n\to 
\widehat{K}_n$} \ \ 
For a proper arithmetic variety $X$, let $\widehat{\mathscr P}(X)$ 
be the category of hermitian vector bundles on $X$ and 
$\widehat{G}(X)=G(\widehat{\mathscr P}(X))$ its $G$-construction.
Then there is a homomorphism of complexes 
$$
\CCGS :C_*(\vert \widehat{G}(X)\vert )\overset{\CUB }{\longrightarrow }
\widehat{\CUB }_*(X)\overset{\CCGS }{\longrightarrow }\mathscr A_*(X).
$$

\vskip 1pc
\begin{prop} \ \ 
The morphism of simplicial sets 
$\chi :\Delta (T^1\wedge \widehat{G}(X))\to \widehat{S}(X)$ yields 
an isomorphism 
$$
\widehat{\chi }_*:\widehat{\pi }_n(\vert \widehat{G}(X)\vert , \CCGS )
\simeq \widehat{\pi }_{n+1}(\vert \widehat{S}(X)\vert , \CCGS )
$$
by $[(f, \omega )]\mapsto [(\chi (1\wedge f), -\omega )]$.
Hence there is a canonical isomorphism 
$$
\widehat{K}_n(X)\simeq \widehat{\pi }_n(\vert \widehat{G}(X)\vert , \CCGS )
$$
if $n\geq 1$.
\end{prop}

{\it Proof}: \ 
It is obvious that the map $(f, \omega )\mapsto 
(\chi (1\wedge f), -\omega )$ gives rise to a homomorphism 
of the modified homotopy groups.
Consider the following commutative diagram:
$$
\begin{CD}
\pi _{n+1}(\vert \widehat{G}(X)\vert ) @>{\rho }>> 
\widetilde{\mathscr A}_{n+1}(X) @>>> 
\widehat{\pi }_n(\vert \widehat{G}(X)\vert , \CCGS ) @>>> 
\pi _n(\vert \widehat{G}(X)\vert ) @>>> 0  \\
@VV{-\chi _*}V @VV{-\text{id}}V @VV{\widehat{\chi }_*}V 
@VV{\chi _*}V  @.  \\
\pi _{n+2}(\vert \widehat{S}(X)\vert ) @>{\rho }>> 
\widetilde{\mathscr A}_{n+1}(X) @>>> 
\widehat{\pi }_{n+1}(\vert \widehat{S}(X)\vert , \CCGS ) @>>> 
\pi _{n+1}(\vert \widehat{S}(X)\vert ) @>>> 0,  \\
\end{CD}
$$
where the upper and lower sequences are exact by Thm.2.3 and $\chi _*$ 
is an isomorphism by Thm.5.1.
Hence the proposition follows from the five lemma.
\qed 

\vskip 1pc
If we set $\widehat{G}^{(2)}(X)=G^{(2)}(\widehat{\mathscr P}(X))$, 
then there is a homomorphism of complexes 
$$
\CCGS :C_*(\vert \widehat{G}^{(2)}(X)\vert )\overset{\CUB }{\longrightarrow }
\widehat{\CUB }_*(X)\overset{\CCGS }{\longrightarrow }\mathscr A_*(X)
$$
and the following square is commutative by Prop.5.2:
$$
\begin{CD}
C_*(\vert \widehat{G}(X)\vert ) @>{\CCGS }>> \mathscr A_*(X)  \\
@VV{R_*}V   @VV{\text{id}}V  \\
C_*(\vert \widehat{G}^{(2)}(X)\vert ) @>{\CCGS }>> \mathscr A_*(X).
\end{CD}
$$
Hence $R$ induces an isomorphism 
$$
\widehat{R}_*:\widehat{\pi }_n(\vert \widehat{G}(X)\vert , \CCGS )
\simeq \widehat{\pi }_n(\vert \widehat{G}^{(2)}(X)\vert , \CCGS )
$$
by Thm.2.3.

Before discussing product in $\widehat{K}_*(X)$, let us recall the 
product in $\widehat{\mathscr K}_0(X)$ defined by Gillet and Soul\'{e} 
in \cite{gilsoul}.
The pairing 
$$
\widehat{\mathscr K}_0(X)\times \widehat{\mathscr K}_0(X)\to 
\widehat{\mathscr K}_0(X)
$$ 
is given by 
$$
[(\overline{E}, \omega )]\times [(\overline{F}, \tau )]=
[(\overline{E}\otimes \overline{F}, \CCGS _0(\overline{E})\bullet \tau 
+\omega \bullet \CCGS _0(\overline{F})+\omega \bullet d_{\mathscr A}\tau )].
$$
Then $\widehat{\mathscr K}_0(X)$ becomes a commutative associative 
algebra.

To construct a product in higher arithmetic $K$-theory, we will 
use the $G$-construction.
However, we have not had any expression of $\widehat{K}_0(X)$ by 
means of the $G$-construction.
Hence we have to distinguish the cases including the $K_0$-group 
from the general case.

Let $(\overline{E}, \eta )$ be a pair of a hermitian vector bundle 
on $X$ and $\eta \in \widetilde{A}_1(X)$ and let $(f, \omega )$ be 
a pair of a pointed cellular map $f: S^n\to \vert \widehat{G}(X)\vert $ 
and $\omega \in \widetilde{A}_{n+1}(X)$.
Let us define the product of these pairs as 
$$
(\overline{E}, \eta )\times (f, \omega )=(\overline{E}\otimes f, 
\CCGS _0(\overline{E})\bullet \omega +\eta \bullet \CCGS _n(f)+
\eta \bullet d_{\mathscr A}\omega ),
$$
where $\overline{E}\otimes f:S^n\overset{f}{\to }
\vert \widehat{G}(X)\vert \overset{\overline{E}\otimes }{\longrightarrow }
\vert \widehat{G}(X)\vert $.

\vskip 1pc
\begin{thm} \ \ 
The above product gives rise to a pairing 
$$
\times :\widehat{\mathscr K}_0(X)\times \widehat{K}_n(X)\to \widehat{K}_n(X).
$$
\end{thm}

{\it Proof}: \ 
To prove the theorem, we have to show that 
$(\overline{E}, \eta )\times (f, \omega )$ is compatible with 
the equivalence relations of $\widehat{\mathscr K}_0(X)$ and 
$\widehat{K}_n(X)$.
Let us first show the compatibility with the relation of $\widehat{K}_n(X)$.

Let $H:S^n\times I/\{*\}\times I\to \vert \widehat{G}(X)\vert $ be 
a cellular homotopy from $(f, \omega )$ to $(f^{\prime }, \omega ^{\prime })$.
Then $\omega ^{\prime }-\omega =(-1)^{n+1}\CCGS _{n+1}(H)$ and the map 
$$
\overline{E}\otimes H:S^n\times I/\{*\}\times I\overset{H}{\to }
\vert \widehat{G}(X)\vert \overset{\overline{E}\otimes }{\longrightarrow }
\vert \widehat{G}(X)\vert 
$$
is a cellular homotopy from $\overline{E}\otimes f$ to 
$\overline{E}\otimes f^{\prime }$.
Furthermore, by Prop.4.6 we have 
\begin{align*}
\CCGS _{n+1}(\overline{E}\otimes H)&=\CCGS _0(\overline{E})\bullet 
\CCGS _{n+1}(H)  \\
&=(-1)^{n+1}\CCGS _0(\overline{E})\bullet (\omega ^{\prime }-\omega ).
\end{align*}
This tells that $\overline{E}\otimes H$ is a cellular homotopy 
from $(\overline{E}, \eta )\times (f, \omega )$ to 
$(\overline{E}, \eta )\times (f^{\prime }, \omega ^{\prime })$.

Next we show the compatibility with the relation of 
$\widehat{\mathscr K}_0(X)$.
Let 
$$
\mathscr E:0\to \overline{E}\to \overline{F}\to \overline{G}\to 0
$$
be a short exact sequence of hermitian vector bundles on $X$.
Consider the following $1$-dimensional subcomplex of 
$\vert \widehat{G}(X)\vert $:
\begin{center}
\setlength{\unitlength}{1mm}
\begin{picture}(52,11)
 \put(10,6){\vector(1,0){18}}
 \put(46,6){\vector(-1,0){18}}
 \put(1,2){$\scriptstyle{(\overline{E}\oplus \overline{G}, 0)}$}
 \put(22,2){$\scriptstyle{(\overline{F}\oplus \overline{G}, \overline{G})}$}
 \put(44,2){$\scriptstyle{(\overline{F}, 0)}$,}
 \put(17,7.5){$e_1$}
 \put(35,7.5){$e_2$}
 \put(10,6){\circle*{1}}
 \put(28,6){\circle*{1}}
 \put(46,6){\circle*{1}}
\end{picture}
\end{center}
where 
$$
e_1=\left(
\setlength{\unitlength}{1mm}
\begin{picture}(63,10)
  \put(15,7){\vector(1,0){9}}
  \put(32,5){\vector(0,-1){7}}
  \put(2,6){$\overline{E}\oplus \overline{G}$}
  \put(26,6){$\overline{F}\oplus \overline{G}$}
  \put(30,-6){$\overline{G}$}
  \put(47,7){\vector(1,0){9}}
  \put(59,5){\vector(0,-1){7}}
  \put(43,6){$0$}
  \put(57,6){$\overline{G}$}
  \put(57,-6){$\overline{G}$}
\end{picture}
\right), 
e_2=\left(
\setlength{\unitlength}{1mm}
\begin{picture}(54,10)
  \put(7,7){\vector(1,0){9}}
  \put(23,5){\vector(0,-1){7}}
  \put(2,6){$\overline{F}$}
  \put(17,6){$\overline{F}\oplus \overline{G}$}
  \put(21,-6){$\overline{G}$}
  \put(38,7){\vector(1,0){9}}
  \put(50,5){\vector(0,-1){7}}
  \put(34,6){$0$}
  \put(48,6){$\overline{G}$}
  \put(48,-6){$\overline{G}$}
\end{picture}
\right).
$$
We denote by $\iota _{\mathscr E}:I\to \vert \widehat{G}(X)\vert $ 
the cellular map such that $\iota _{\mathscr E}(I)=e_1e_2^{-1}$.

For a pointed cellular map $f:S^n\to \vert \widehat{G}(X)\vert $, let 
$$
H:S^n\times I/\{*\}\times I\overset{T}{\to }I\times S^n/I\times \{*\}
\overset{\iota _{\mathscr E}\wedge f}{\longrightarrow }
\vert \widehat{G}(X)\vert \wedge \vert \widehat{G}(X)\vert 
\overset{m^G}{\longrightarrow }\vert \widehat{G}^{(2)}(X)\vert ,
$$
where $T$ is the switching map defined by $T(s, t)=(t, s)$.
If $H_0(s)=H(s, 0)$, then $H_0$ is written as 
$$
S^n\overset{f}{\to }\vert \widehat{G}(X)\vert 
\overset{\iota _{\overline{E}\oplus \overline{G}}\wedge \IId}
{\longrightarrow }\vert \widehat{G}(X)\vert \wedge 
\vert \widehat{G}(X)\vert \overset{m^G}{\longrightarrow }
\vert \widehat{G}^{(2)}(X)\vert ,
$$
where $\iota _{\overline{E}\oplus \overline{G}}: S^0\to 
\vert \widehat{G}(X)\vert $ is the pointed map determined by 
$(\overline{E}\oplus \overline{G}, 0)\in \widehat{G}_0(X)$.
Since the diagram 
$$
\begin{CD} 
\vert \widehat{G}(X)\vert  
@>{\iota _{\overline{E}\oplus \overline{G}}\wedge \IId}>> 
\vert \widehat{G}(X)\vert \wedge \vert \widehat{G}(X)\vert   \\
@VV{(\overline{E}\oplus \overline{G})\otimes }V @VV{m^G}V  \\
\vert \widehat{G}(X)\vert  @>{R}>> \vert \widehat{G}^{(2)}(X)\vert 
\end{CD}
$$
is commutative, it follows that 
$H_0=R((\overline{E}\oplus \overline{G})\otimes f)$. 
If $H_1(s)=H(s, 1)$, then we can show $H_1=R(\overline{F}\otimes f)$ 
in the same way.
Moreover, Prop.4.6 implies 
\begin{align*}
\CCGS _{n+1}(H)&=(-1)^n\CCGS _{n+1}(m^G_*(\iota _{\mathscr E}
\wedge f)_*([I\times S^n]))  \\
&\equiv (-1)^n\CCGS _1(\iota _{\mathscr E})\bullet \CCGS _n(f)  \\
&=(-1)^n\CCGS _1(\mathscr E)\bullet \CCGS _n(f)
\end{align*}
modulo $\IIm d_{\mathscr A}$.
Hence $H$ is a cellular homotopy from $(R((\overline{E}\oplus \overline{G})
\otimes f), \CCGS _1(\mathscr E)\bullet \CCGS _n(f))$ to 
$(R(\overline{F}\otimes f), 0)$.
Since $\widehat{R}_*:\widehat{\pi }_n(\vert \widehat{G}(X)\vert , \CCGS )
\to \widehat{\pi }_n(\vert \widehat{G}^{(2)}(X)\vert , \CCGS )$ is 
bijective, we have 
$$
[((\overline{E}\oplus \overline{G})\otimes f, \CCGS _1(\mathscr E)
\bullet \CCGS _n(f))]=[(\overline{F}\otimes f, 0)]
$$
in $\widehat{\pi }_n(\vert \widehat{G}(X)\vert , \CCGS )$.

The short exact sequence $\mathscr E$ gives a relation 
$$
[(\overline{E}, 0)]+[(\overline{G}, 0)]=
[(\overline{F}, -\CCGS _1(\mathscr E))]
$$
in $\widehat{\mathscr K}_0(X)$.
We have 
$$
[(\overline{E}, 0)\times (f, \omega )]*[(\overline{G}, 0)\times 
(f, \omega )]=[((\overline{E}\otimes f)*(\overline{G}\otimes f), 
(\CCGS _0(\overline{E})+\CCGS _0(\overline{G}))\bullet \omega )]
$$
and 
\begin{align*}
[(\overline{F}, &-\CCGS _1(\mathscr E))\times (f, \omega )]  \\
&=[(\overline{F}\otimes f, \CCGS _0(\overline{F})\bullet \omega 
-\CCGS _1(\mathscr E)\bullet \CCGS _n(f)-d_{\mathscr A}
\CCGS _1(\mathscr E)\bullet \omega )]  \\
&=[((\overline{E}\oplus \overline{G})\otimes f, (\CCGS _0(\overline{E})
+\CCGS _0(\overline{G}))\bullet \omega )]
\end{align*}
in $\widehat{\pi }_n(\vert \widehat{G}(X)\vert , \CCGS )$.
Hence Thm.5.4 follows from Lem.5.5 and Lem.5.6 below.
\qed

\vskip 1pc
\begin{lem} \ \ 
For a pointed cellular map $f:S^n\to \vert \widehat{G}(X)\vert $ and 
two hermitian vector bundles $\overline{E}, \overline{G}$ on $X$, 
we have 
$$
[((\overline{E}\otimes f)\oplus (\overline{G}\otimes f), 0)]=
[((\overline{E}\otimes f)*(\overline{G}\otimes f), 0)]
$$
in $\widehat{\pi }_n(\vert \widehat{G}(X)\vert , \CCGS )$.
\end{lem}

{\it Proof}: \ 
Let us first describe the map $(\overline{E}\otimes f)\oplus (\overline{G}
\otimes f)$ explicitly.
Since $f$ is a pointed cellular map, the map 
$$
S^n\overset{\Delta }{\hookrightarrow }S^n\times S^n 
\overset{(\overline{E}\otimes f)\times 
(\overline{G}\otimes f)}{\longrightarrow }\vert \Delta (\widehat{G}(X)
\times \widehat{G}(X))\vert 
$$
is also a pointed cellular map.
Moreover, the direct sum of hermitian vector bundles induces a 
morphism of simplicial sets 
$$
\oplus :\Delta (\widehat{G}(X)\times \widehat{G}(X))\to \widehat{G}(X).
$$
Then the map $(\overline{E}\otimes f)\oplus (\overline{G}\otimes f)$ 
is expressed as the composition of these two cellular maps, that is, 
$$
(\overline{E}\otimes f)\oplus (\overline{G}\otimes f):
S^n\overset{\Delta }{\hookrightarrow }S^n\times S^n 
\overset{(\overline{E}\otimes f)\times 
(\overline{G}\otimes f)}{\longrightarrow }\vert \Delta (\widehat{G}(X)
\times \widehat{G}(X))\vert \overset{\oplus }{\to }\vert 
\widehat{G}(X)\vert .
$$

Consider the homomorphism of complexes 
$$
\CCGS \oplus \CCGS :
C_*(\vert \Delta (\widehat{G}(X)\times \widehat{G}(X))\vert )\to 
\mathscr A_*(X)\oplus \mathscr A_*(X)
$$
defined by $(E, F)\mapsto (\CCGS _n(E), \CCGS _n(F))$ for 
$E, F\in \widehat{G}_n(X)$ and the inclusion 
$$
in_1 \ (\text{resp.~$in_2$}):\widehat{G}(X)\hookrightarrow 
\Delta (\widehat{G}(X)\times \widehat{G}(X))
$$
defined by $in_1(t)=(t, *)$ (resp.~$in_2(t)=(*, t)$). 
Then we have the following commutative diagram:
$$
\begin{CD}
C_*(\vert \widehat{G}(X)\vert ) @>{\CCGS }>> \mathscr A_*(X)  \\
@VV{\text{${in_1}_*$ (resp.~${in_2}_*$)}}V  
@VV{\text{$in_1$ (resp.~$in_2$)}}V \\
C_*(\vert \Delta (\widehat{G}(X)\times \widehat{G}(X))\vert ) 
@>{\CCGS \oplus \CCGS }>> \mathscr A_*(X)\oplus \mathscr A_*(X),
\end{CD}
$$
where the right vertical arrow is defined by 
$in_1(\omega )=(\omega , 0)$ (resp.~$in_2(\omega )=(0, \omega )$).
On the other hand, the projection 
$$
pr_1 \ (\text{resp.~$pr_2$}):\Delta (\widehat{G}(X)\times \widehat{G}(X))
\to \widehat{G}(X)
$$ 
defined by $pr_1(x, y)=x$ (resp.~$pr_2(x, y)=y$) is also a morphism 
of simplicial sets and it satisfies the following 
commutative diagram:
$$
\begin{CD}
C_*(\vert \Delta (\widehat{G}(X)\times \widehat{G}(X))\vert ) 
@>{\CCGS \oplus \CCGS }>> \mathscr A_*(X)\oplus \mathscr A_*(X)  \\
@VV{\text{${pr_1}_*$ (resp.~${pr_2}_*$)}}V  
@VV{\text{$pr_1$ (resp.~$pr_2$)}}V \\
C_*(\vert \widehat{G}(X)\vert ) @>{\CCGS }>> \mathscr A_*(X),
\end{CD}
$$
where the right vertical arrow is defined by 
$pr_1(\omega , \tau )=\omega $ (resp.~$pr_2(\omega , \tau )=\tau $).
Hence we have four homomorphisms 
\begin{center}
\setlength{\unitlength}{1mm}
\begin{picture}(110,10)
  \put(1,4.5){$\widehat{\pi }_n(\vert \Delta (\widehat{G}(X)\times 
             \widehat{G}(X))\vert , \CCGS \oplus \CCGS )$}
  \put(83,4.5){$\widehat{\pi }_n(\vert \widehat{G}(X)\vert , \CCGS )$}
  \put(81,5){\vector(-1,0){12}}
  \put(69,6.5){\vector(1,0){12}}
  \put(73,8){$\scriptstyle{\widehat{pr_j}_*}$}
  \put(73,1){$\scriptstyle{\widehat{in_j}_*}$}
\end{picture}
\end{center}
that form an isomorphism 
$$
\widehat{\pi }_n(\vert \widehat{G}(X)\vert , \CCGS )
\oplus \widehat{\pi }_n(\vert \widehat{G}(X)\vert , \CCGS )\simeq 
\widehat{\pi }_n(\vert \Delta (\widehat{G}(X)\times \widehat{G}(X))\vert , 
\CCGS \oplus \CCGS )
$$
by $(x, y)\mapsto \widehat{in_1}_*(x)*\widehat{in_2}_*(y)$.
The inverse of it is $\widehat{pr_1}_*\oplus \widehat{pr_2}_*$.
Since 
\begin{gather*}
\widehat{pr_1}_*([(((\overline{E}\otimes f)\times (\overline{G}\otimes f))
\Delta , 0)])=[(\overline{E}\otimes f, 0)]  \\
\intertext{and}
\widehat{pr_2}_*([(((\overline{E}\otimes f)\times (\overline{G}\otimes f))
\Delta , 0)])=[(\overline{G}\otimes f, 0)],
\end{gather*}
we have 
$$
[(((\overline{E}\otimes f)\times (\overline{G}\otimes f))\Delta , 0)]=
\widehat{in_1}_*([(\overline{E}\otimes f, 0)])*\widehat{in_2}_*
([(\overline{G}\otimes f, 0)])
$$ 
in $\widehat{\pi }_*(\vert \Delta (\widehat{G}(X)\times \widehat{G}(X))
\vert , \CCGS \oplus \CCGS )$.

The following commutative diagram 
$$
\begin{CD}
C_*(\vert \Delta (\widehat{G}(X)\times \widehat{G}(X))\vert ) 
@>{\CCGS \oplus \CCGS }>> \mathscr A_*(X)\oplus \mathscr A_*(X)  \\
@VV{\oplus _*}V   @VV{+}V  \\
C_*(\vert \widehat{G}(X)\vert )  @>{\CCGS }>> \mathscr A_*(X)
\end{CD}
$$
yields a homomorphism 
$$
\widehat{\oplus }_*:\widehat{\pi }_n(\vert \Delta 
(\widehat{G}(X)\times \widehat{G}(X))\vert , \CCGS \oplus \CCGS )\to 
\widehat{\pi }_n(\vert \widehat{G}(X)\vert , \CCGS ).
$$
Since $\widehat{\oplus }_*\widehat{in_j}_*$ is the identity homomorphism, 
we have 
\begin{align*}
[((\overline{E}\otimes f)\oplus (\overline{G}\otimes f), 0)]&=
\widehat{\oplus }_*\left([(((\overline{E}\otimes f)\times 
(\overline{G}\otimes f))\Delta , 0)]\right)  \\
&=[(\overline{E}\otimes f, 0)]*[(\overline{G}\otimes f, 0)]  \\
&=[((\overline{E}\otimes f)*(\overline{G}\otimes f), 0)],
\end{align*}
which completes the proof.
\qed

\vskip 1pc
\begin{lem} \ \ 
In the same notations as Lem.5.5, we have 
$$
[((\overline{E}\oplus \overline{G})\otimes f, 0)]=
[((\overline{E}\otimes f)\oplus (\overline{G}\otimes f), 0)]
$$
in $\widehat{\pi }_n(\vert \widehat{G}(X)\vert , \CCGS )$. 
\end{lem}

{\it Proof}: \ 
Consider the following diagram:
$$
\begin{CD}
\Delta (\widehat{G}(X)\times \widehat{G}(X)) 
@>{(\overline{E}\otimes \ \ )\times (\overline{G}\otimes \ \ )}>> 
\Delta (\widehat{G}(X)\times \widehat{G}(X))  @.  \\
@AA{\Delta }A  @VV{\oplus }V  @.  \\
\widehat{G}(X) @>{(\overline{E}\oplus \overline{G})\otimes }>>
\widehat{G}(X) @>{R}>> \widehat{G}^{(2)}(X). 
\end{CD}
$$
This square is commutative up to homotopy.
Let $\alpha _0$ be the upper map of the diagram and $\alpha _1$ 
the lower map.
Then for $P=(P^{\pm })\in \widehat{G}_n(X)$, $\alpha _0(P)$ and 
$\alpha _1(P)$ are elements of $\widehat{G}_0\widehat{G}_n(X)$ 
described as 
\begin{align*}
\alpha _0(P)&=((\overline{E}\otimes P^+)\oplus (\overline{G}\otimes P^+), 
(\overline{E}\otimes P^-)\oplus (\overline{G}\otimes P^-), 0, 0)  \\
\intertext{and}
\alpha _1(P)&=((\overline{E}\oplus \overline{G})\otimes P^+, 
(\overline{E}\oplus \overline{G})\otimes P^-, 0, 0).
\end{align*}
The canonical isometries $(\overline{E}\otimes P^{\pm })\oplus 
(\overline{G}\otimes P^{\pm })\simeq (\overline{E}\oplus \overline{G})
\otimes P^{\pm }$ give an element of $\widehat{G}_1\widehat{G}_n(X)$ 
whose Bott-Chern form is zero.
Collecting these elements for all $P=(P^{\pm })$ provides 
a morphism of bisimplicial sets 
$$
\Psi :\Delta [1]\times \widehat{G}(X)\to \widehat{G}^{(2)}(X)
$$
such that $\Psi (0, s)=\alpha _0(s)$ and $\Psi (1, s)=\alpha _1(s)$.
Therefore for any pointed celular map 
$f:S^n\to \vert \widehat{G}(X)\vert $, 
$$
H:S^n\times I/\{*\}\times I\overset{T}{\to }I\times S^n/I\times \{*\}
\overset{\IId \times f}{\longrightarrow }I\times 
\vert \widehat{G}(X)\vert /I\times \{*\}\overset{\vert \Psi \vert }
{\longrightarrow }\vert \widehat{G}^{(2)}(X)\vert 
$$
is a cellular homotopy from $R((\overline{E}\otimes f)\oplus 
(\overline{G}\otimes f))$ to $R((\overline{E}\oplus \overline{G})\otimes f)$ 
such that $\CCGS _{n+1}(H)=0$.
Since $\widehat{R}_*:\widehat{\pi }_n(\vert \widehat{G}(X)\vert , \CCGS )
\to \widehat{\pi }_n(\vert \widehat{G}^{(2)}(X)\vert , \CCGS )$ is 
bijective, we have 
$$
[((\overline{E}\otimes f)\oplus (\overline{G}\otimes f), 0)]=
[((\overline{E}\oplus \overline{G})\otimes f, 0)]
$$
in $\widehat{\pi }_n(\vert \widehat{G}(X)\vert , \CCGS )$.
\qed 

\vskip 1pc
In the same way as above, we can define a pairing 
$\widehat{K}_n(X)\times \widehat{\mathscr K}_0(X)\to \widehat{K}_n(X)$ by 
$$
[(f, \omega )]\times [(\overline{E}, \eta )]=[(f\otimes \overline{E}, 
(-1)^n\CCGS _n(f)\bullet \eta +\omega \bullet \CCGS _0(\overline{E})
+\omega \bullet d_{\mathscr A}\eta )],
$$
where $f\otimes \overline{E}:S^n\overset{f}{\to }\vert \widehat{G}(X)\vert 
\overset{\otimes \overline{E}}{\longrightarrow }\vert \widehat{G}(X)\vert $.
Combining these pairings with the isomorphism 
$\widehat{\alpha }:\widehat{\mathscr K}_0(X)\simeq \widehat{K}_0(X)$, 
we can obtain a pairing 
$$
\times :\widehat{K}_n(X)\times \widehat{K}_m(X)\to \widehat{K}_{n+m}(X)
$$
when $n=0$ or $m=0$.

\vskip 1pc
\subsection{The pairing of higher arithmetic $K$-theory} \ \ 
In this subsection we define a pairing $\widehat{K}_n(X)\times 
\widehat{K}_m(X)\to \widehat{K}_{n+m}(X)$ in the case of $n, m\geq 1$.
For two pointed cellular maps $f:S^n\to \vert \widehat{G}(X)\vert $ 
and $g:S^m\to \vert \widehat{G}(X)\vert $, we set 
$$
f\times g:S^{n+m}=S^n\wedge S^m\overset{f\wedge g}{\longrightarrow}
\vert \widehat{G}(X)\vert \wedge \vert \widehat{G}(X)\vert 
\overset{m^G}{\longrightarrow }\vert \widehat{G}^{(2)}(X)\vert .
$$
For $\omega \in \widetilde{\mathscr A}_{n+1}(X)$ and 
$\tau \in \widetilde{\mathscr A}_{m+1}(X)$, 
the product $(f, \omega )\times (g, \tau )$ is defined as 
\begin{align*}
(f, &\omega )\times (g, \tau )  \\
&=(f\times g, (-1)^n\CCGS _n(f)\bullet \tau+\omega \bullet \CCGS _m(g)
+\omega \bullet d_{\mathscr A}\tau
+(-1)^n\CCGS _n(f)\vartriangle \CCGS _m(g)).
\end{align*}
 
\vskip 1pc
\begin{prop} \ \ 
The above product gives rise to a pairing 
$$
\widehat{m}^G_*:\widehat{\pi }_n(\vert \widehat{G}(X)\vert , \CCGS )
\times \widehat{\pi }_m(\vert \widehat{G}(X)\vert , \CCGS )\to 
\widehat{\pi }_{n+m}(\vert \widehat{G}^{(2)}(X)\vert , \CCGS ).
$$
\end{prop}  

{\it Proof}: \ 
For a cellular homotopy $H$ from $(f, \omega )$ to 
$(f^{\prime }, \omega ^{\prime })$, let $\widetilde{H}$ be 
a cellular map given as 
$$
\widetilde{H}:S^{n+m}\times I/\{*\}\times I
\longrightarrow \vert \widehat{G}(X)\vert \wedge \vert \widehat{G}(X)\vert 
\overset{m^G}{\longrightarrow }\vert \widehat{G}^{(2)}(X)\vert ,
$$
where the first map is defined by $(s_1, s_2, t)\mapsto (H(s_1, t), g(s_2))$ 
for $s_1\in S^n, s_2\in S^m$ and $t\in I$.
Then $\widetilde{H}$ is a homotopy from $f\times g$ to 
$f^{\prime }\times g$ and Prop.4.6 and Prop.5.2 imply 
{\allowdisplaybreaks
\begin{align*}
\CCGS _{n+m+1}(\widetilde{H})&=(-1)^m\CCGS _{n+m+1}(m^G_*
(H\times g)_*([S^n\times I\times S^m]))  \\
&\equiv (-1)^m\CCGS _{n+1}(H)\bullet \CCGS _m(g)+(-1)^{n+m+1}
\CCGS _n(\partial H_*([S^n\times I]))\vartriangle \CCGS _m(g)  \\
&=(-1)^{n+m+1}(\omega ^{\prime }-\omega )\bullet \CCGS _m(g)+(-1)^{m+1}
(\CCGS _n(f^{\prime })-\CCGS _n(f))\vartriangle \CCGS _m(g)  \\
&=(-1)^{n+m+1}(\omega ^{\prime }\bullet \CCGS _m(g)+(-1)^n
\CCGS _n(f^{\prime })\vartriangle \CCGS _m(g))  \\*
&\hskip 1pc -(-1)^{n+m+1}(\omega \bullet \CCGS _m(g)
+(-1)^n\CCGS _n(f)\vartriangle \CCGS _m(g))
\end{align*}}
modulo $\IIm d_{\mathscr A}$.
This tells that the map $\widetilde{H}$ is a cellular homotopy from 
$(f, \omega )\times (g, \tau )$ to $(f^{\prime }, \omega ^{\prime })
\times (g, \tau )$.

If $H^{\prime }$ is a cellular homotopy from $(g, \tau )$ to 
$(g^{\prime }, \tau ^{\prime })$, we can verify in the same way 
that the map 
$$
S^{n+m}\times I/\{*\}\times I\overset{f\wedge H^{\prime }}{\longrightarrow }
\vert \widehat{G}(X)\vert \wedge \vert \widehat{G}(X)\vert 
\overset{m^G}{\longrightarrow }\vert \widehat{G}^{(2)}(X)\vert 
$$
is a cellular homotopy from $(f, \omega )\times (g, \tau )$ to 
$(f, \omega )\times (g^{\prime }, \tau ^{\prime })$.
\qed

\vskip 1pc
\begin{defn} \ \ 
For $n, m\geq 1$, we define a product in higher arithmetic $K$-theory
$$
\times :\widehat{K}_n(X)\times \widehat{K}_m(X)\to \widehat{K}_{n+m}(X)
$$
by the following homomorphism: 
$$
\widehat{\pi }_n(\vert \widehat{G}(X)\vert , \CCGS )\times 
\widehat{\pi }_m(\vert \widehat{G}(X)\vert , \CCGS )
\overset{\widehat{m}^G_*}{\longrightarrow }
\widehat{\pi }_{n+m}(\vert \widehat{G}^{(2)}(X)\vert , \CCGS )
\overset{\widehat{R}_*^{-1}}{\longrightarrow }
\widehat{\pi }_{n+m}(\vert \widehat{G}(X)\vert , \CCGS ).
$$
\end{defn}

\vskip 1pc
\begin{prop} \ \ 
The Chern form map respects the products, that is, we have 
$$
\CCGS _{n+m}(x\times y)=\CCGS _n(x)\bullet \CCGS _m(y)
$$
for $x\in \widehat{K}_n(X)$ and $y\in \widehat{K}_m(X)$.
\end{prop}

{\it Proof}: \ 
We assume $n, m\geq 1$.
Let us define the Chern form map on $\widehat{\pi }_{n+m}(\vert 
\widehat{G}^{(2)}(X)\vert , \CCGS )$ as 
$$
\CCGS _{n+m}([(f, \omega )])=\CCGS (f_*([S^{n+m}]))+
d_{\mathscr A}\omega \in \mathscr A_{n+m}(X).
$$
Then it follows that $\CCGS _{n+m}(\widehat{R}_*(x))=\CCGS _{n+m}(x)$ 
for any $x\in \widehat{\pi }_{n+m}(\vert \widehat{G}(X)\vert , \CCGS )$.
Hence it is sufficient to show that 
$$
\CCGS _{n+m}(\widehat{m}^G_*(x, y))=\CCGS _n(x)\bullet \CCGS _m(y).
$$

For $x=[(f, \omega )]$ and $y=[(g, \tau )]$, Prop.4.6 implies 
{\allowdisplaybreaks
\begin{align*}
\CCGS _{n+m}(\widehat{m}^G_*(x, y))&=\CCGS _{n+m}(f\times g)+d_{\mathscr A}
\left((-1)^n\CCGS _n(f)\bullet \tau +\omega \bullet \CCGS _m(g)\right. \\
&\hskip 10pc \left.+\omega \bullet d_{\mathscr A}\tau +(-1)^n
\CCGS _n(f)\vartriangle \CCGS _m(g)\right)   \\
&=(\CCGS _n(f)+d_{\mathscr A}\omega )\bullet (\CCGS _m(g)
+d_{\mathscr A}\tau )  \\
&=\CCGS _n(x)\bullet \CCGS _m(y). 
\end{align*}}
The case where $n=0$ or $m=0$ is trivial.
\qed

\vskip 1pc
{\it Remark}: \ 
As far as the author knows, there does not exist a homomorphism 
$$
\CUB :C_*(\vert \widehat{S}^{(2)}(X)\vert )[-2]\to \widehat{\CUB }_*(X)
$$
that makes the following diagram commutative:
$$
\begin{CD}
C_*(\vert \widehat{S}(X)\vert )[-1]\otimes 
C_*(\vert \widehat{S}(X)\vert )[-1] @>{\CUB \otimes \CUB }>> 
\widehat{\CUB }_*(X)\otimes \widehat{\CUB }_*(X)  \\
@VV{m_*}V  @VV{\otimes }V  \\
C_*(\vert \widehat{S}^{(2)}(X)\vert )[-2] @>{\CUB }>> \widehat{\CUB }_*(X).
\end{CD}
$$
Hence it seems impossible to define a product in 
$\widehat{K}_*(X)$ using the $S$-construction.

\vskip 1pc
\subsection{The commutativity of the product} \ \ 
In this subsection we discuss the commutativity of the product 
in $\widehat{K}_*(X)$.
When $n=0$ or $m=0$, it is easy to prove that the product 
$\widehat{K}_n(X)\times \widehat{K}_m(X)\to \widehat{K}_{n+m}(X)$ is 
commutative.
So we consider only the case of $n, m\geq 1$.

For a small exact category $\mathfrak A$, let 
$L:G_n\mathfrak A\to G_nG_0\mathfrak A$ be the morphism defined by \linebreak
$L(E^+, E^-)=(E^+, 0, E^-, 0)$. 
Then it induces a homotopy equivalent morphism 
$$
L:G\mathfrak A\to G^{(2)}\mathfrak A
$$
and it is homotopy equivalent to the morphism $R$.

Let us define a new product in $\widehat{K}_*(X)$ using the morphism $L$.
For a proper arithmetic variety $X$, the diagram 
$$
\begin{CD}
C_*(\vert \widehat{G}(X)\vert ) @>{\CCGS }>> \mathscr A_*(X)  \\
@VV{L_*}V  @VV{\text{id}}V  \\
C_*(\vert \widehat{G}^{(2)}(X)\vert ) @>{\CCGS }>> \mathscr A_*(X)
\end{CD}
$$
is commutative.
Therefore we have an isomorphism 
$$
\widehat{L}_*:\widehat{\pi }_n(\vert \widehat{G}(X)\vert , \CCGS )\simeq 
\widehat{\pi }_n(\vert \widehat{G}^{(2)}(X)\vert , \CCGS ).
$$

\vskip 1pc
\begin{defn} \ \ 
For $n, m\geq 1$, we define a new product 
$$
\underset{L}{\times }:\widehat{K}_n(X)\times \widehat{K}_m(X)
\to \widehat{K}_{n+m}(X)
$$
by the pairing  
$$
\widehat{\pi }_n(\vert \widehat{G}(X)\vert , \CCGS )\times 
\widehat{\pi }_m(\vert \widehat{G}(X)\vert , \CCGS )
\overset{\widehat{m}^G_*}{\to }
\widehat{\pi }_m(\vert \widehat{G}^{(2)}(X)\vert , \CCGS )
\overset{\widehat{L}_*^{-1}}{\to }
\widehat{\pi }_m(\vert \widehat{G}(X)\vert , \CCGS ).
$$
\end{defn}

\vskip 1pc
Let us compare this new product $\underset{L}{\times }$ with 
the previously defined one. 
Let $T:\widehat{S}_n\widehat{S}_m(X)\to \widehat{S}_m\widehat{S}_n(X)$ 
be the switching map defined by $T(E)_{(i,j)\times (\alpha ,\beta )}=
E_{(\alpha ,\beta )\times (i,j)}$.
Then the map 
$$
\underset{n, m}{\coprod }\widehat{S}_n\widehat{S}_m(X)\times 
\Delta ^n\times \Delta ^m\to \underset{n, m}{\coprod }\widehat{S}_m
\widehat{S}_n(X)\times \Delta ^m\times \Delta ^n
$$
defined by $(x, t_1, t_2)\mapsto (T(x), t_2, t_1)$ induces 
an involution of the topological realization 
$$
\mathscr T:\vert \widehat{S}^{(2)}(X)\vert \to 
\vert \widehat{S}^{(2)}(X)\vert .
$$
Similarly we can define the switching map $T:\widehat{G}_n\widehat{G}_m(X)
\to \widehat{G}_m\widehat{G}_n(X)$ and the involution 
$$
\mathscr T:\vert \widehat{G}^{(2)}(X)\vert \to 
\vert \widehat{G}^{(2)}(X)\vert .
$$

\vskip 1pc
\begin{lem} \ \ 
The diagram 
$$
\begin{CD}
C_*(\vert \widehat{G}^{(2)}(X)\vert ) @>{\CCGS }>> \mathscr A_*(X)  \\
@VV{\mathscr T_*}V  @VV{\IId }V  \\
C_*(\vert \widehat{G}^{(2)}(X)\vert ) @>{\CCGS }>> \mathscr A_*(X)
\end{CD}
$$
is commutative.
Hence we have an isomorphism 
$$
\widehat{\mathscr T}_*:\widehat{\pi }_n(\vert \widehat{G}^{(2)}(X)\vert , 
\CCGS )\simeq \widehat{\pi }_n(\vert \widehat{G}^{(2)}(X)\vert , \CCGS )
$$
by $[(f, \omega )]\mapsto [(\mathscr Tf, \omega )]$.
\end{lem} 

{\it Proof}: \ 
If we denote by $[E]$ the element of $C_*(\vert \widehat{G}^{(2)}(X)\vert )$ 
determined by $E\in \widehat{G}_n\widehat{G}_m(X)$, then we have 
$\mathscr T_*([E])=(-1)^{nm}[T(E)]$.
Hence we have 
\begin{align*}
\CCGS _{n+m}(\mathscr T_*([E]))&=(-1)^{nm}\CCGS _{n+m}(\CUB (T(E))) \\
&=(-1)^{nm}\CCGS _{n+m}(T_{n,m}(\CUB (E))). 
\end{align*}
In the above, $T_{n,m}(\mathcal F)$ for an exact metrized 
$(n+m)$-cube $\mathcal F$ is defined as 
$$
T_{n,m}(\mathcal F)_{\alpha _1, \cdots , \alpha _{n+m}}=
\mathcal F_{\alpha _{n+1}, \cdots , \alpha _{n+m}, \alpha _1, \cdots , 
\alpha _n}.
$$
Then it is easy to see that $\CCGS _{n+m}(T_{n,m}(\mathcal F))=(-1)^{nm}
\CCGS _{n+m}(\mathcal F)$, therefore 
$\CCGS _{n+m}(\mathscr T_*([E]))=$ \linebreak 
$\CCGS _{n+m}(E)$.
\qed

\vskip 1pc
\begin{prop} \ \ 
Let $x\in \widehat{K}_n(X)$ and $y\in \widehat{K}_m(X)$ with 
$n, m\geq 1$. 
Then we have 
$$
x\times y=(-1)^{nm}y\underset{L}{\times }x.
$$
\end{prop}

{\it Proof}: \ 
For two pointed spaces $S_1$ and $S_2$, let 
$T:S_1\wedge S_2\to S_2\wedge S_1$ denote the switching map given by 
$T(s_1, s_2)=(s_2, s_1)$.
For two pointed cellular maps $f:S^n\to \vert \widehat{G}(X)\vert $ and 
$g:S^m\to \vert \widehat{G}(X)\vert $, we consider the following 
diagram:
$$
\begin{CD}
S^n\wedge S^m @>{f\wedge g}>> \vert \widehat{G}(X)\vert \wedge 
\vert \widehat{G}(X)\vert @>{m^G}>> 
\vert \widehat{G}^{(2)}(X)\vert   \\
@VV{T}V @VV{T}V  @VV{\mathscr T}V \\
S^m\wedge S^n @>{g\wedge f}>> \vert \widehat{G}(X)\vert \wedge 
\vert \widehat{G}(X)\vert @>{m^G}>> 
\vert \widehat{G}^{(2)}(X)\vert .
\end{CD}
$$
The left square is obviously commutative, but the right one is not.
In fact, for $E\in \widehat{G}_n(X)$ and 
$F\in \widehat{G}_m(X)$ we have 
\begin{align*}
\mathscr Tm^G(E, F)_{(i,j)\times (\alpha , \beta )}&=
E_{\alpha , \beta }\otimes F_{i,j}  \\
\intertext{and}
m^GT(E, F)_{(i,j)\times (\alpha , \beta )}&=
F_{i,j}\otimes E_{\alpha , \beta }.
\end{align*}
Hence the homotopy from $\mathscr Tm^G$ to $m^GT$ can be constructed by 
means of the canonical isometry $\overline{P}\otimes \overline{Q}\simeq 
\overline{Q}\otimes \overline{P}$ of hermitian vector bundles.
More precisely, if $\widehat{G}^{(3)}(X)$ is the trisimplicial set 
obtained by taking $G$-construction three times and 
$R:\widehat{G}^{(2)}(X)\to \widehat{G}^{(3)}(X)$ is the homotopy 
equivalent map induced from the map $R:\widehat{G}_n\widehat{G}_m(X)\to 
\widehat{G}_0\widehat{G}_n\widehat{G}_m(X)$ given by 
$R(E)^{+ \pm \pm }=E^{\pm \pm }$ and $R(E)^{- \pm \pm }=0$, then 
a homotopy between the diagram 
$$
\begin{CD}
\vert \widehat{G}(X)\vert \wedge \vert \widehat{G}(X)\vert  
@>{m^G}>> \vert \widehat{G}^{(2)}(X)\vert @.  \\
@VV{T}V  @VV{\mathscr T}V  @.   \\
\vert \widehat{G}(X)\vert \wedge \vert \widehat{G}(X)\vert  
@>{m^G}>> \vert \widehat{G}^{(2)}(X)\vert @>{R}>> 
\vert \widehat{G}^{(3)}(X)\vert 
\end{CD}
$$
is given by means of $\overline{P}\otimes \overline{Q}\simeq 
\overline{Q}\otimes \overline{P}$ in the same way as $\Psi $ 
in the proof of Lem.5.6.
Hence we have 
$$
[(\mathscr T(f\times g), 0)]=[((g\times f)T, 0)]
$$
in $\widehat{\pi }_{n+m}(\vert \widehat{G}^{(2)}(X)\vert , \CCGS )$.

If $x=[(f, \omega )]$ and $y=[(g, \tau )]$, then we have 
\begin{align*}
\widehat{\mathscr T}_*&\widehat{m}^G_*([(f, \omega )], [(g, \tau )])  \\
&=[(\mathscr T(f\times g), (-1)^n\CCGS _n(f)\bullet \tau 
+\omega \bullet \CCGS _m(g)+\omega \bullet d_{\mathscr A}\tau 
+(-1)^n\CCGS _n(f)\vartriangle \CCGS _m(g))]  \\
&=[((g\times f)T, (-1)^{nm}\tau \bullet \CCGS _n(f)
+(-1)^{(n+1)m}\CCGS _m(g)\bullet \omega +(-1)^{nm}\tau \bullet 
d_{\mathscr A}\omega   \\
&\hskip 10pc +(-1)^{(n+1)m}\CCGS _m(g)\vartriangle \CCGS _n(f))].
\end{align*}
Since $T:S^{n+m}\to S^{n+m}$ is homotopic to $(-1)^{nm}\IId _{S^{n+m}}$, 
we have 
$$
\widehat{\mathscr T}_*\widehat{m}^G_*([(f, \omega )], [(g, \tau )])=
(-1)^{nm}\widehat{m}^G_*([(g, \tau )], [(f, \omega )])
$$
in $\widehat{\pi }_{n+m}(\vert \widehat{G}^{(2)}(X)\vert , \CCGS )$.
The commutativity of the square 
$$
\begin{CD}
\vert \widehat{G}^{(2)}(X)\vert @<{R}<< \vert \widehat{G}(X)\vert   \\
@VV{\mathscr T}V  @VV{\IId }V  \\
\vert \widehat{G}^{(2)}(X)\vert @<{L}<< \vert \widehat{G}(X)\vert 
\end{CD}
$$
implies 
\begin{align*}
(-1)^{nm}\widehat{L}_*([(g, \tau )]\underset{L}{\times }[(f, \omega )])
&=(-1)^{nm}\widehat{m}^G_*([(g, \tau )], [(f, \omega )])  \\
&=\widehat{\mathscr T}_*\widehat{m}^G_*([(f, \omega )], [(g, \tau )])  \\
&=\widehat{\mathscr T}_*\widehat{R}_*([(f, \omega )]\times [(g, \tau )]) \\
&=\widehat{L}_*([(f, \omega )]\times [(g, \tau )]).
\end{align*}
Since $\widehat{L}_*$ is bijective, we have completed the proof.
\qed

\vskip 1pc
\begin{prop} \ \ 
For $x\in \widehat{\pi }_n(\vert \widehat{G}(X)\vert , \CCGS )$, 
the element 
$$
\widehat{R}_*(x)-\widehat{L}_*(x)\in 
\widehat{\pi }_n(\vert \widehat{G}^{(2)}(X)\vert , \CCGS )
$$
is contained in $\IIm (\widetilde{A}_{n+1}(X)\to 
\widehat{\pi }_n(\vert \widehat{G}^{(2)}(X)\vert , \CCGS ))$ and 
$2(\widehat{R}_*(x)-\widehat{L}_*(x))=0$.
In particular, for $x\in \widehat{K}_n(X)$ and $y\in \widehat{K}_m(X)$ 
with $n, m\geq 1$, the element $x\times y-x\underset{L}{\times }y$ 
is contained in 
$\IIm (\widetilde{\mathscr A}_{n+m+1}(X)\to \widehat{K}_{n+m}(X))$ 
and $2(x\times y-x\underset{L}{\times }y)=0$.
\end{prop} 

{\it Proof}: \ 
Since $R$ and $L$ are homotopy equivalent, it is clear that 
$\widehat{R}_*(x)-\widehat{L}_*(x)$ is contained in 
$\IIm (\widetilde{A}_{n+1}(X)\to 
\widehat{\pi }_n(\vert \widehat{G}^{(2)}(X)\vert , \CCGS ))$.
Let $f:S^n\to \vert \widehat{G}(X)\vert $ be a pointed cellular map. 
Then there is a pointed cellular map 
$$
H:S^n\times I/\{*\}\times I\to \vert \widehat{G}^{(2)}(X)\vert 
$$
such that $H(s, 0)=Rf(s)$ and $H(s, 1)=Lf(s)$.
Let 
$$
H^{\prime }=\mathscr TH:S^n\times I/\{*\}\times I\to 
\vert \widehat{G}^{(2)}(X)\vert ,
$$
then we have $H^{\prime }(s, 0)=Lf(s)$ and $H^{\prime }(s, 1)=Rf(s)$.
The commutative square in Lem.5.11 implies 
$\CCGS _{n+1}(H^{\prime })=\CCGS _{n+1}(H)$.
Gluing the maps $H$ and $H^{\prime }$ on the boundaries, we obtain 
a cellular map 
$$
H\cup H^{\prime }:S^n\times T^1/\{*\}\times T^1\to 
\vert \widehat{G}^{(2)}(X)\vert ,
$$
where $T^1$ is the barycentric subdivision of $S^1$.

\vskip 1pc
\begin{lem} \ \ 
If $n\geq 1$, there is a surjection 
$$
p:S^{n+1}\to S^n\times S^1/\{*\}\times S^1
$$
such that $p^{-1}((S^n-\{*\})\times S^1)\to 
(S^n-\{*\})\times S^1$ is a homeomorphism.
\end{lem} 

{\it Proof}: \ 
We describe the space $S^{n+1}$ as follows:
$$
S^{n+1}=\left\{(z, t_1, \cdots , t_n)\in \mathbb C\times \mathbb R^n; 
|z|^2+t_1^2+\cdots +t_n^2=1\right\}.
$$
Let $S^{n-1}=\{(0, t_1, \cdots , t_n)\in S^{n+1}\}$. 
Then a homeomorphism $S^{n+1}-S^{n-1}\to B^n\times S^1$
is defined by 
$$
(z, t_1, \cdots , t_n)\mapsto \left((t_1, \cdots , t_n), 
\frac{z}{\sqrt{1-t_1^2-\cdots -t_n^2}}\right),
$$
where $B^n=\{(t_1, \cdots , t_n)\in \mathbb R^n; t_1^2+\cdots +t_n^2<1\}$.
Since $S^n\times S^1/\{*\}\times S^1$ is the one-point compactification 
of $B^n\times S^1$, this homeomorphism can be extended to 
$$
p:S^{n+1}\to S^n\times S^1/\{*\}\times S^1
$$
which satisfies the condition as mentioned above.
\qed

\vskip 1pc
Let us return to the proof of Prop.5.13.
Since $T^1$ is the barycentric subdivision of $S^1$, 
we have a pointed cellular map 
$$
F:S^{n+1}\overset{p}{\to }S^n\times T^1/\{*\}\times T^1
\overset{H\cup H^{\prime }}{\longrightarrow }
\vert \widehat{G}^{(2)}(X)\vert .
$$
Then we have 
$$
\CCGS _{n+1}(F)=\CCGS _{n+1}(H)+\CCGS _{n+1}(H^{\prime })
=2\CCGS _{n+1}(H).
$$
In other words, $2\CCGS _{n+1}(H)$ is contained in the image of 
$\pi _{n+1}(\vert \widehat{G}^{(2)}(X)\vert )\to 
\widetilde{\mathscr A}_{n+1}(X)$.
Hence $2[(0, \CCGS _{n+1}(H))]=0$ in 
$\widehat{\pi }_n(\vert \widehat{G}^{(2)}(X)\vert , \CCGS )$ by 
Thm.2.3.
For any $x=[(f, \omega )]\in $ \linebreak
$\widehat{\pi }_n(\vert \widehat{G}(X)\vert , 
\CCGS )$, we have 
\begin{align*}
\widehat{R}_*(x)-\widehat{L}_*(x)&=[(Rf, 0)]-[(Lf, 0)]  \\
&=(-1)^{n+1}[(0, \CCGS _{n+1}(H))],
\end{align*}
therefore $2(\widehat{R}_*(x)-\widehat{L}_*(x))=0$. 
\qed

\vskip 1pc
Prop.5.12 and Prop.5.13 yield the following theorem:

\vskip 1pc
\begin{thm} \ \ 
Let $x\in \widehat{K}_n(X)$ and $y\in \widehat{K}_m(X)$.
Then $x\times y-(-1)^{nm}y\times x$ is a $2$-torsion element 
contained in $\IIm (\widetilde{\mathscr A}_{n+m+1}(X)\to 
\widehat{K}_{n+m}(X))$.
Hence the product in $\widehat{K}_*(X)$ is graded commutative up to 
$2$-torsion.
\end{thm}

\vskip 1pc
\subsection{The lack of the associativity} \ \ 
In this subsection we discuss the associativity of the product in 
$\widehat{K}_*(X)$.
Let $\widehat{G}^{(3)}(X)$ be the trisimplicial set given by taking 
$\widehat{G}$ three times. 
Let $m^G$ denote the morphisms of trisimplicial sets 
$$
\widehat{G}^{(2)}(X)\wedge \widehat{G}(X)\to \widehat{G}^{(3)}(X)
$$
and
$$
\widehat{G}(X)\wedge \widehat{G}^{(2)}(X)\to \widehat{G}^{(3)}(X)
$$
given by the tensor product of hermitian vector bundles.
A homotopy equivalent map 
$$
R:\widehat{G}(X)\to \widehat{G}^{(3)}(X)
$$
is given by $R(E)^{+ + \pm }=E^{\pm }$ and $R(E)^{+ - \pm }=R(E)^{- + \pm }
=R(E)^{- - \pm }=0$ for $E=(E^{\pm })\in \widehat{G}_n(X)$.
Under the above notations, the following diagram 
$$
\begin{CD}
\widehat{G}(X)\wedge \widehat{G}(X)\wedge \widehat{G}(X) 
@>{m^G\wedge 1}>> \widehat{G}^{(2)}(X)\wedge \widehat{G}(X) \\
@VV{1\wedge m^G}V  @VV{m^G}V   \\
\widehat{G}(X)\wedge \widehat{G}^{(2)}(X) @>{m^G}>> \widehat{G}^{(3)}(X) 
\end{CD}
$$
is commutative up to the homotopy arising from the natural isometry 
$(\overline{E}\otimes \overline{F})\otimes \overline{G}\simeq 
\overline{E}\otimes (\overline{F}\otimes \overline{G})$ of 
hermitian vector bundles.
This commutative diagram is the source of the associativity of the product 
in usual algebraic $K$-theory $K_*(X)$.

For two pointed cellular maps $f:S^n\to \vert \widehat{G}^{(2)}(X)\vert $ 
and $g:S^m\to \vert \widehat{G}(X)\vert $, let 
$$
f\times g:S^{n+m}\overset{f\wedge g}{\longrightarrow }
\vert \widehat{G}^{(2)}(X)\vert \wedge \vert \widehat{G}(X)\vert 
\overset{\widehat{m}^G}{\longrightarrow }\vert \widehat{G}^{(3)}(X)\vert .
$$
Let us define a pairing 
$$
\widehat{m}^G_*:\widehat{\pi }_n(\vert \widehat{G}^{(2)}(X)\vert , \CCGS )
\times \widehat{\pi }_m(\vert \widehat{G}(X)\vert , \CCGS )\to 
\widehat{\pi }_{n+m}(\vert \widehat{G}^{(3)}(X)\vert , \CCGS )
$$
by 
\begin{align*}
([(f&, \omega )], [(g, \tau )])\mapsto  \\
&[(f\times g, (-1)^n\CCGS _n(f)\bullet \tau +\omega \bullet \CCGS _m(g)+
\omega \bullet d_{\mathscr A}\tau +(-1)^n\CCGS _n(f)\vartriangle 
\CCGS _m(g))].
\end{align*}
The well-definedness of the pairing can be shown in the same way 
as the proof of Thm.5.7.
We can also define a pairing 
$$
\widehat{m}^G_*:\widehat{\pi }_n(\vert \widehat{G}(X)\vert , \CCGS )\times 
\widehat{\pi }_m(\vert \widehat{G}^{(2)}(X)\vert , \CCGS )\to 
\widehat{\pi }_{n+m}(\vert \widehat{G}^{(3)}(X)\vert , \CCGS )
$$
by the same expression as above.
Then the associativity of the product in $\widehat{K}_*(X)$ is equivalent 
to the commutativity of the following diagram: 
$$
\setlength{\unitlength}{1mm}
\begin{picture}(130,71)
  \put(35,65){$\scriptstyle{\widehat{\pi }_n
          (\vert \widehat{G}(X)\vert , \CCGS )\times 
          \widehat{\pi }_m(\vert \widehat{G}(X)\vert , \CCGS )\times 
          \widehat{\pi }_l(\vert \widehat{G}(X)\vert , \CCGS )}$}
  \put(2,43){$\scriptstyle{\widehat{\pi }_{n+m}
       (\vert \widehat{G}^{(2)}(X)\vert , \CCGS )
     \times \widehat{\pi }_l(\vert \widehat{G}(X)\vert , \CCGS )}$}
  \put(83,43){$\scriptstyle{\widehat{\pi }_n
       (\vert \widehat{G}(X)\vert , \CCGS )\times 
       \widehat{\pi }_{m+l}(\vert \widehat{G}^{(2)}(X)\vert , \CCGS )}$}
  \put(54,20){$\scriptstyle{\widehat{\pi }_{n+m+l}
              (\vert \widehat{G}^{(3)}(X)\vert , \CCGS )}$}
  \put(56,2){$\scriptstyle{\widehat{\pi }_{n+m+l}
            (\vert \widehat{G}(X)\vert , \CCGS )}$.}
  \put(45,63){\vector(-1,-1){15}}
  \put(95,63){\vector(1,-1){15}}
  \put(25,40){\vector(2,-1){32}}
  \put(115,40){\vector(-2,-1){32}}
  \put(70,7){\vector(0,1){10}}
  \put(40,55){$\scriptstyle{\widehat{m}^G_*\times \text{id}}$}
  \put(104,55){$\scriptstyle{\text{id}\times \widehat{m}^G_*}$}
  \put(34,30){$\scriptstyle{\widehat{m}^G_*}$}
  \put(101,30){$\scriptstyle{\widehat{m}^G_*}$}
  \put(71,11){$\scriptstyle{\widehat{R}_*}$}
\end{picture}
$$
However, the diagram is not commutative.
Let $[(f, \omega )]\in 
\widehat{\pi }_n(\vert \widehat{G}(X)\vert , \CCGS )$, 
$[(g, \tau )]\in \widehat{\pi }_m(\vert \widehat{G}(X)\vert , \CCGS )$ 
and $[(h, \eta )]\in \widehat{\pi }_l(\vert \widehat{G}(X)\vert , \CCGS )$.
Then in the same way as the proof of Lem.5.6, the identity 
$$
[((f\times g)\times h, 0)]=[(f\times (g\times h), 0)]
$$
in $\widehat{\pi }_{n+m+l}(\vert \widehat{G}^{(3)}(X)\vert , \CCGS )$
can be shown.
The quadrusimplicial set $\widehat{G}^{(4)}(X)$ is necessary to show 
this identity.
Hence an easy calculation leads us to the following:

\vskip 1pc
\begin{prop} \ \ 
We have 
$$
\widehat{m}^G_*\left(\widehat{m}^G_*([(f, \omega )], [(g, \tau )]), 
[(h, \eta )]\right)-\widehat{m}^G_*\left([(f, \omega )], 
\widehat{m}^G_*([(g, \tau )], [(h, \eta )])\right)=
[(0, r(f, g, h, \omega , \tau , \eta ))]
$$
in $\widehat{\pi }_{n+m+l}(\vert \widehat{G}^{(3)}(X)\vert , \CCGS )$, 
where 
{\allowdisplaybreaks
\begin{align*}
r&(f, g, h, \omega , \tau , \eta ) \\
&=(-1)^n((\CCGS _n(f)+d_{\mathscr A}\omega )\bullet \tau )\bullet 
(\CCGS _l(h)+d_{\mathscr A}\eta )  \\*
&\hskip 4pc -(-1)^n(\CCGS _n(f)+d_{\mathscr A}\omega )
\bullet (\tau \bullet (\CCGS _l(h)+d_{\mathscr A}\eta ))  \\
&\hskip 1pc +(-1)^{n+m}\CCGS _{n+m}(f\times g)\bullet \eta 
+(-1)^n(\CCGS _n(f)\vartriangle \CCGS _m(g))\bullet d_{\mathscr A}\eta \\*
&\hskip 4pc -(-1)^{n+m}\CCGS _n(f)\bullet (\CCGS _m(g)\bullet \eta )  \\
&\hskip 1pc +(\omega \bullet \CCGS _m(g))\bullet \CCGS _l(h)
-(-1)^{n+m}d_{\mathscr A}\omega \bullet (\CCGS _m(g)\vartriangle 
\CCGS _l(h))-\omega \bullet \CCGS _{m+l}(g\times h)  \\
&\hskip 1pc +(\omega \bullet \CCGS _m(g))\bullet d_{\mathscr A}\eta 
-(-1)^{n+m}d_{\mathscr A}\omega \bullet (\CCGS _m(g)\bullet \eta )  \\
&\hskip 1pc +(-1)^n(\CCGS _n(f)\vartriangle \CCGS _m(g))\bullet \CCGS _l(h)
+(-1)^{n+m}\CCGS _{n+m}(f\times g)\vartriangle \CCGS _l(h)  \\
&\hskip 1pc -(-1)^{n+m}\CCGS _n(f)\bullet (\CCGS _m(g)\vartriangle 
\CCGS _l(h))-(-1)^n\CCGS _n(f)\vartriangle \CCGS _{m+l}(g\times h).
\end{align*}}
\end{prop}

\vskip 1pc
For $\omega \in \mathscr A_n(X)$ or $\mathscr D_n(X)$ 
and for an integer $i$ satisfying $1\leq i\leq n$, we set 
$$
\omega ^{(-i,-n+i-1)}=\underset{p}{\sum }\omega ^{(p-i,p-n+i-1)}, 
$$
where $\omega ^{(p-i,p-n+i-1)}$ is the $(p-i,p-n+i-1)$-component of 
the differential form $\omega $.
Then for $\omega \in \mathscr A_n(X)$ and $\tau \in \mathscr A_m(X)$, 
it follows that 
$$
\omega \vartriangle \tau =\frac{1}{4\pi \sqrt{-1}}\sum a^{n,m}_{i,j}
\omega ^{(-i,-n+i-1)}\wedge \tau ^{(-j,-m+j-1)}.
$$

\vskip 1pc
\begin{lem} \ 
{\rm 1)} \ We assume $n, m, l\geq 1$. 
For $\alpha \in \mathscr A_n(X)$, $\beta \in \mathscr A_m(X)$ and 
$\gamma \in \mathscr A_l(X)$, we have 
\begin{align*}
(\alpha &\bullet \beta )\bullet \gamma -\alpha \bullet 
(\beta \bullet \gamma )  \\
&=\frac{(-1)^{n+m}}{4(2\pi \sqrt{-1})^2}(\partial \alpha ^{(-1,-n)}
+\overline{\partial }\alpha ^{(-n,-1)})\wedge (\partial \beta ^{(-1,-m)}
+\overline{\partial }\beta ^{(-m,-1)})\wedge \gamma  \\
&\hskip 1pc -\frac{1}{4(2\pi \sqrt{-1})^2}\alpha \wedge 
(\partial \beta ^{(-1,-m)}+\overline{\partial }\beta ^{(-m,-1)})\wedge 
(\partial \gamma ^{(-1,-l)}+\overline{\partial }\gamma ^{(-l,-1)}). 
\end{align*}

\noindent
{\rm 2)} \ If $d_{\mathscr A}\alpha =0$ and $d_{\mathscr A}\beta =0$, then 
we have 
$$
(\alpha \bullet \beta )\bullet \gamma -\alpha \bullet 
(\beta \bullet \gamma )\equiv 
\begin{cases}
\frac{-1}{4(2\pi \sqrt{-1})^2}\alpha \wedge d\beta \wedge d_{\mathscr A}
\gamma, \ &l\geq 2,  \\
0, &l=1 
\end{cases} 
$$
modulo $\IIm d_{\mathscr A}$. 

\noindent
{\rm 3)} \ If $d_{\mathscr A}\beta =0$ and $d_{\mathscr A}\gamma =0$, then 
we have 
$$
(\alpha \bullet \beta )\bullet \gamma -\alpha \bullet 
(\beta \bullet \gamma )\equiv 
\begin{cases}
\frac{(-1)^{n+m}}{4(2\pi \sqrt{-1})^2}d_{\mathscr A}\alpha \wedge 
d\beta \wedge \gamma , \ &n\geq 2, \\
0, &n=1
\end{cases}
$$
modulo $\IIm d_{\mathscr A}$. 

\noindent
{\rm 4)} \ If $d_{\mathscr A}\alpha =0$ and $d_{\mathscr A}\gamma =0$, then 
we have 
$$
(\alpha \bullet \beta )\bullet \gamma -\alpha \bullet 
(\beta \bullet \gamma )\equiv 
\begin{cases}
\frac{(-1)^{m}}{4(2\pi \sqrt{-1})^2}\alpha \wedge dd_{\mathscr A}\beta 
\wedge \gamma , \ &m\geq 2, \\
0, &m=1
\end{cases}
$$
modulo $\IIm d_{\mathscr A}$. 
\end{lem}

{\it Proof}: \ 
The identity in 1) follows from an easy calculation. 
If $l\geq 2$ and $d_{\mathscr A}\alpha =d_{\mathscr A}\beta =0$, 
then we have 
\begin{align*}
&(\alpha \bullet \beta )\bullet \gamma -\alpha \bullet 
(\beta \bullet \gamma )  \\
&\hskip 1pc =\frac{(-1)^{n+m}}{4(2\pi \sqrt{-1})^2}
d\alpha \wedge d\beta \wedge \gamma -\frac{1}{4(2\pi \sqrt{-1})^2}
\alpha \wedge d\beta \wedge (d\gamma +d_{\mathscr A}\gamma )  \\
&\hskip 1pc =\frac{(-1)^{n+m}}{4(2\pi \sqrt{-1})^2}
d(\alpha \wedge d\beta \wedge \gamma )-\frac{1}{4(2\pi \sqrt{-1})^2}
\alpha \wedge d\beta \wedge d_{\mathscr A}\gamma .
\end{align*}
The form $\frac{1}{4(2\pi \sqrt{-1})^2}\alpha \wedge d\beta \wedge \gamma $ 
is contained in $\mathscr A^{2(p+q+r)-n-m-l-1}(X, p+q+r)$ and 
$$
d_{\mathscr A}\left(\frac{1}{4(2\pi \sqrt{-1})^2}(\alpha \wedge d\beta 
\wedge \gamma )\right)=\frac{-1}{4(2\pi \sqrt{-1})^2}d(\alpha \wedge d\beta 
\wedge \gamma ).
$$
Hence we have 
$$
(\alpha \bullet \beta )\bullet \gamma -\alpha \bullet 
(\beta \bullet \gamma )\equiv \frac{-1}{4(2\pi \sqrt{-1})^2}
\alpha \wedge d\beta \wedge d_{\mathscr A}\gamma  
$$
modulo $\IIm d_{\mathscr A}$. 
When $l=1$, we have 
$$
(\alpha \bullet \beta )\bullet \gamma -\alpha \bullet 
(\beta \bullet \gamma )=d_{\mathscr A}\left(
\frac{(-1)^{n+m+1}}{4(2\pi \sqrt{-1})^2}(\alpha \wedge d\beta 
\wedge \gamma )\right).
$$
Hence 2) holds.
The identities in 3) and 4) can be proved in the same way. 
\qed

\vskip 1pc
We calculate $r(f, g, h, \omega , \tau , \eta )$ by using Lem.5.17. 
If $n, m, l\geq 1$, then we have 
{\allowdisplaybreaks
\begin{align*}
(-1)^n&((\CCGS _n(f)+d_{\mathscr A}\omega )\bullet \tau )
\bullet (\CCGS _l(h)+d_{\mathscr A}\eta )  \\*
&-(-1)^n(\CCGS _n(f)+d_{\mathscr A}\omega )
\bullet (\tau \bullet (\CCGS _l(h)+d_{\mathscr A}\eta ))  \\
&\hskip 1pc \equiv \frac{(-1)^{n+m+1}}{4(2\pi \sqrt{-1})^2}
(\CCGS _n(f)+d_{\mathscr A}\omega )\wedge dd_{\mathscr A}
\tau \wedge (\CCGS _l(h)+d_{\mathscr A}\eta )  \\
\end{align*}
and
\begin{align*}
(\omega &\bullet \CCGS _m(g))\bullet d_{\mathscr A}\eta 
-(-1)^{n+m}d_{\mathscr A}\omega \bullet (\CCGS _m(g)\bullet \eta )  \\*
&\equiv \frac{(-1)^{n+m+1}}{4(2\pi \sqrt{-1})^2}
d_{\mathscr A}\omega \wedge d\CCGS _m(g)\wedge 
d_{\mathscr A}\eta 
\end{align*}
modulo $\IIm d_{\mathscr A}$. 
Since 
$$
\CCGS _{n+m}(f\times g)=\CCGS _n(f)\bullet \CCGS _m(g)+(-1)^{n+1}
d_{\mathscr A}(\CCGS _n(f)\vartriangle \CCGS _m(g))
$$
by Prop.4.6, we have 
{\allowdisplaybreaks
\begin{align*}
(-1)^{n+m}&\CCGS _{n+m}(f\times g)\bullet \eta 
+(-1)^n(\CCGS _n(f)\vartriangle \CCGS _m(g))\bullet d_{\mathscr A}\eta \\*
&-(-1)^{n+m}\CCGS _{n}(f)\bullet (\CCGS _m(g)\bullet \eta )  \\
&\hskip 1pc =(-1)^{n+m}(\CCGS _n(f)\bullet \CCGS _m(g))\bullet \eta 
-(-1)^{n+m}\CCGS _{n}(f)\bullet (\CCGS _m(g)\bullet \eta )  \\
&\hskip 1pc \equiv \frac{(-1)^{n+m+1}}{4(2\pi \sqrt{-1})^2}
\CCGS _n(f)\wedge d\CCGS _m(g)\wedge d_{\mathscr A}\eta 
\end{align*}}
modulo $\IIm d_{\mathscr A}$.
In the same way as above we have 
\begin{align*}
(\omega &\bullet \CCGS _m(g))\bullet \CCGS _l(h)
-\omega \bullet \CCGS _{m+l}(g\times h)
-(-1)^{n+m}d_{\mathscr A}\omega \bullet 
(\CCGS _m(g)\vartriangle \CCGS _l(h)) \\
&\equiv \frac{(-1)^{n+m+1}}{4(2\pi \sqrt{-1})^2}
d_{\mathscr A}\omega \wedge d\CCGS _m(g)\wedge \CCGS _l(h) 
\end{align*}
modulo $\IIm d_{\mathscr A}$.
As for the last four terms of the identity in Prop.5.16, we have 
the following:

\vskip 1pc
\begin{prop} \ \ 
If $n ,m, l\geq 1$, then we have 
\begin{multline*}
(-1)^n(\CCGS _n(f)\vartriangle \CCGS _m(g))\bullet \CCGS _l(h)
-(-1)^{n+m}\CCGS _n(f)\bullet (\CCGS _m(g)\vartriangle \CCGS _l(h)) \\*
+(-1)^{n+m}\CCGS _{n+m}(f\times g)\vartriangle \CCGS _l(h)
-(-1)^n\CCGS _n(f)\vartriangle \CCGS _{m+l}(g\times h)  \\*
\equiv (-1)^{n+m+1}\frac{1}{4(2\pi \sqrt{-1})^2}
\CCGS _n(f)\wedge d\CCGS _m(g)\wedge \CCGS _l(h)
\end{multline*}
modulo $\IIm d_{\mathscr A}$.
\end{prop}

\vskip 1pc
We will prove this identity in \S 5.9.
Substituting these identities into that in Prop.5.16 yields 
\begin{align*}
&r(f, g, h, \omega , \tau , \eta )  \\
&\hskip 1pc \equiv \frac{(-1)^{n+m+1}}{4(2\pi \sqrt{-1})^2}
(\CCGS _n(f)+d_{\mathscr A}\omega )\wedge 
d(\CCGS _m(g)+d_{\mathscr A}\tau )\wedge (\CCGS _l(h)+d_{\mathscr A}\eta ) 
\end{align*}
modulo $\IIm d_{\mathscr A}$.

\vskip 1pc
\begin{thm} \ \ 
The product in higher arithmetic K-theory does not satisfy the 
associativity. 
In fact, if $x\in \widehat{K}_n(X), y\in \widehat{K}_m(X)$ and 
$z\in \widehat{K}_l(X)$ for $n, m, l\geq 1$, we have 
$$
(x\times y)\times z-x\times (y\times z)=[(0, 
\sfrac{(-1)^{n+m+1}}{4(2\pi \sqrt{-1})^2}\CCGS _n(x)\wedge 
d\CCGS _m(y)\wedge \CCGS _l(z))]
$$
in $\widehat{K}_{n+m+l}(X)$. 
Hence $(x\times y)\times z=x\times (y\times z)$ holds when 
$nml=0$ or $y\in K_m(\overline{X})$ or $x=y=z$.
\end{thm}

{\it Proof}: \ 
When $n,m,l\geq 1$, we have already proved this identity. 
The identity $(x\times y)\times z=x\times (y\times z)$ in the case of 
$nml=0$ follows from the definition of the product and 
Lem.5.17.
\qed

\vskip 1pc
\subsection{The product in Arakelov $K$-theory} \ \ 
For a proper Arakelov variety $\overline{X}=(X, h_X)$, let us define 
a pairing 
$$
K_n(\overline{X})\times K_m(\overline{X})\to 
K_{n+m}(\overline{X})
$$
by $(x, y)\mapsto \sigma (x\times y)$, where $\sigma $ is the 
harmonic projection defined in \S 3.3.

\vskip 1pc
\begin{thm} \ \ 
The above pairing makes $K_*(\overline{X})$ a graded associative 
algebra.
That is to say, it follows that 
$$
\sigma (\sigma (x\times y)\times z)=\sigma (x\times \sigma (y\times z))
$$
for $x, y, z\in K_*(\overline{X})$.
\end{thm}

{\it Proof}: \ 
This identity is obvious when $nml=0$, so we may assume $nml\geq 1$.
We first prove the identity 
$$
\sigma (\sigma (x\times y)\times z)=\sigma ((x\times y)\times z)
$$
for $x\in K_n(\overline{X}), y\in K_m(\overline{X})$ and 
$z\in K_l(\overline{X})$.
It follows from the definition of $\sigma $ that 
$$
\sigma (x\times y)=x\times y+[(0, \alpha )]
$$
for $\alpha \in \mathscr A_{n+m+1}(X)$ with $\mathscr H(\alpha )=0$.
Then we have 
$$
\sigma (x\times y)\times z=(x\times y)\times z+
[(0, \alpha \bullet \CCGS _l(z))],
$$
therefore 
$$
\sigma (\sigma (x\times y)\times z)=(x\times y)\times z+
[(0, \alpha \bullet \CCGS _l(z)+\beta )]
$$
for $\beta \in \mathscr A_{n+m+l+1}(X)$ with $\mathscr H(\beta )=0$.
Since $\CCGS _l(z)$ is harmonic, we have 
$$
\alpha \bullet \CCGS _l(z)=(-1)^{n+m+1}(\partial \alpha ^{(-1,-n-m)}
-\overline{\partial }\alpha ^{(-n-m,-1)})\wedge \CCGS _l(z).
$$
Since $\partial \alpha ^{(-1,-n-m)}\wedge \CCGS _l(z)$ is 
$\partial $-exact and 
$\overline{\partial }\alpha ^{(-n-m,-1)}\wedge \CCGS _l(z)$ 
is $\overline{\partial }$-exact, we have 
$\mathscr H(\alpha \bullet \CCGS _l(z))=0$, so 
$\mathscr H(\alpha \bullet \CCGS _l(z)+\beta )=0$. 
Therefore we have 
$\sigma (\sigma (x\times y)\times z)=\sigma ((x\times y)\times z)$. 

In the same way, we can show 
$\sigma (x\times \sigma (y\times z))=\sigma (x\times (y\times z))$.
Hence Thm.5.19 implies 
$$
\sigma (\sigma (x\times y)\times z)=\sigma (x\times \sigma (y\times z)),
$$
which completes the proof.
\qed

\vskip 1pc
\subsection{Proof of Prop.5.18} \ \ 
We will prove the following formula, which is not renormalized: 
\begin{multline*}
(-1)^n(\CC _n(f)\vartriangle \CC _m(g))\bullet \CC _l(h)
-(-1)^{n+m}\CC _n(f)\bullet (\CC _m(g)\vartriangle \CC _l(h)) \\
+(-1)^{n+m}\CC _{n+m}(f\times g)\vartriangle \CC _l(h)
-(-1)^n\CC _n(f)\vartriangle \CC _{m+l}(g\times h)  \\
\equiv (-1)^{n+m+1}\CC _n(f)\wedge d\CC _m(g)\wedge \CC _l(h)\end{multline*}
modulo $\IIm d_{\mathscr D}$.
Let $\Phi $ denote the left hand side of the above identity and let 
$\Phi (f)$ (resp.~ $\Phi (g)$ and $\Phi (h)$) be the part of $\Phi $ 
including derivatives of $\CC _n(f)$ (resp.~ $\CC _m(g)$ and 
$\CC _l(h)$). 
In other words, $\Phi =\Phi (f)+\Phi (g)+\Phi (h)$ such that 
{\allowdisplaybreaks
\begin{align*}
\Phi (f)&=(-1)^{m+1}(\partial \CC _n(f)^{(-1,-n)}-\overline{\partial }
\CC _n(f)^{(-n,-1)})\wedge (\CC _m(g)\vartriangle \CC _l(h))  \\
&\hskip 1pc +(-1)^m((\partial \CC _n(f)^{(-1,-n)}-\overline{\partial }
\CC _n(f)^{(-n,-1)})\wedge \CC _m(g))\vartriangle \CC _l(h)  \\
&\hskip 1pc +(-1)^m\left(\underset{1\leq j\leq m}
{\underset{1\leq i\leq n}{\sum }}a^{n,m}_{i,j}d\CC _n(f)^{(-n+i-1,-i)}
\wedge \CC _m(g)^{(-m+j-1,-j)}\right)\vartriangle \CC _l(h),   \\
\Phi (g)&=(-1)^{n+m}(\CC _n(f)\wedge (\partial \CC _m(g)^{(-1,-m)}
-\overline{\partial }\CC _m(g)^{(-m,-1)}))\vartriangle \CC _l(h)  \\
&\hskip 1pc +(-1)^{n+m+1}\left(\underset{1\leq j\leq m}
{\underset{1\leq i\leq n}{\sum }}a^{n,m}_{i,j}
\CC _n(f)^{(-n+i-1,-i)}\wedge d\CC _m(g)^{(-m+j-1,-j)}\right)
\vartriangle \CC _l(h)  \\
&\hskip 1pc +(-1)^{n+m+1}\CC _n(f)\vartriangle ((\partial \CC _m(g)^{(-1,-m)}
-\overline{\partial }\CC _m(g)^{(-m,-1)})\wedge \CC _l(h))  \\
&\hskip 1pc +(-1)^{n+m+1}\CC _n(f)\vartriangle \left(\underset{1\leq k\leq l}
{\underset{1\leq j\leq m}{\sum }}a^{m,l}_{j,k} 
d\CC _m(g)^{(-m+j-1,-j)}\wedge \CC _l(h)^{(-l+k-1,-k)}\right)  \\
\intertext{and} 
\Phi (h)&=(-1)^n(\CC _n(f)\vartriangle \CC _m(g))\wedge 
(\partial \CC _l(h)^{(-1,-l)}-\overline{\partial }\CC _l(h)^{(-l,-1)}) \\
&\hskip 1pc +(-1)^{n+1}\CC _n(f)\vartriangle (\CC _m(g)\wedge (\partial 
\CC _l(h)^{(-1,-l)}-\overline{\partial }\CC _l(h)^{(-l,-1)}))  \\
&\hskip 1pc +(-1)^n\CC _n(f)\vartriangle \left(\underset{1\leq k\leq l}
{\underset{1\leq j\leq m}{\sum }}a^{m,l}_{j,k}
\CC _m(g)^{(-m+j-1,-j)}\wedge d\CC _l(h)^{(-l+k-1,-k)}\right).
\end{align*}

Let us first calculate $\Phi (f)$. 
It follows from $d_{\mathscr D}(\CC _n(f))=0$ that 
$\partial \CC _n(f)^{(-n+i-1, -i)}=$\linebreak 
$-\overline{\partial }\CC _n(f)^{(-n+i, -i-1)}$ for $1\leq i\leq n-1$. 
Then $\Phi (f)$ is expressed as follows: 
\begin{align*}
\Phi (f)=&b^{n,m,l}_{0,j,k}\underset{1\leq k\leq l}
{\underset{1\leq j\leq m}{\sum }}
\overline{\partial }\CC _n(f)^{(-n,-1)}\wedge 
\CC _m(g)^{(-m+j-1,-j)}\wedge \CC _l(h)^{(-l+k-1,-k)} \\*
&+\underset{1\leq k\leq l}{\underset{1\leq j\leq m}{
\underset{1\leq i\leq n}{\sum }}}
b^{n,m,l}_{i,j,k}\partial \CC _n(f)^{(-n+i-1,-i)}\wedge 
\CC _m(g)^{(-m+j-1,-j)}\wedge \CC _l(h)^{(-l+k-1,-k)}
\end{align*}
where 
{\allowdisplaybreaks
\begin{align*}
b^{n,m,l}_{0,j,k}&=(-1)^ma^{m,l}_{j,k}+(-1)^{m+1}a^{n+m,l}_{j,k}
+(-1)^ma^{n+m,l}_{j,k}\times a^{n,m}_{1,j}  \\*
&=(-1)^ma^{m,l}_{j,k}+2(-1)^{m+1}\sbinom{n+m}{n}^{-1}\sbinom{n+m-j}{n}
a^{n+m,l}_{j,k}, \\
b^{n,m,l}_{n,m,k}&=(-1)^{m+1}a^{m,l}_{j,k}+(-1)^ma^{n+m,l}_{n+j,k}
+(-1)^ma^{n+m,l}_{n+j,k}\times a^{n,m}_{n,j}  \\*
&=(-1)^{m+1}a^{m,l}_{j,k}+2(-1)^m\sbinom{n+m}{n}^{-1}\sbinom{n-j+1}{n}
a^{n+m,l}_{n+j,k}, 
\intertext{and}
b^{n,m,l}_{i,j,k}&=(-1)^{m+1}a^{n+m,l}_{i+j,k}\times a^{n,m}_{i+1,j}
+(-1)^ma^{n+m,l}_{i+j,k}\times a^{n,m}_{i,j}  \\*
&=2(-1)^m\sbinom{n+m}{n}^{-1}\sbinom{n+m-i-j}{n-i}
\sbinom{i+j-1}{i}a^{n+m,l}_{i+j,k} 
\end{align*}
for $1\leq i\leq n-1$.

\vskip 1pc
\begin{lem} \ \ If $1\leq j\leq m$ and $1\leq k\leq l$, then 
$$
\sbinom{n+m}{n}^{-1}\sum_{i=0}^n\sbinom{n+m-i-j}{n-i}
\sbinom{i+j-1}{i}a^{n+m,l}_{i+j,k}=a^{m,l}_{j,k}.
$$
\end{lem}

{\it Proof}: \ 
By Lem.A.2 and Lem.A.3, we have 
\begin{align*}
\text{(LHS)}&=\sbinom{n+m}{n}^{-1}\sum_{i=0}^n\sbinom{n+m-i-j}{n-i}
\sbinom{i+j-1}{i}  \\*
&\hskip 1pc -2\sbinom{n+m}{n}^{-1}\sbinom{n+m+l}{n+m}^{-1}\sum_{i=0}^n
\sbinom{n+m-i-j}{n-i}\sbinom{i+j-1}{i}\sum_{\alpha =0}^{i+j-1}
\sbinom{n+m+l-i-j-k+1}{n+m-\alpha }\sbinom{i+j+k-1}{\alpha }  \\
&=1-2\sbinom{m+l}{m}^{-1}\sum_{\alpha =0}^{j-1}
\sbinom{m+l-j-k+1}{m-\alpha }\sbinom{j+k-1}{\alpha }   \\
&=a^{m,l}_{j,k},
\end{align*}
which completes the proof.
\qed

\vskip 1pc
Let $c^{n,m,l}_{i,j,k}$ be the rational number defined as follows: 
$$
c^{n,m,l}_{i,j,k}=(-1)^ma^{m,l}_{j,k}-2(-1)^m\sbinom{n+m}{n}^{-1}
\sum_{\alpha =0}^{i-1}\sbinom{n+m-\alpha -j}{n-\alpha }
\sbinom{\alpha +j-1}{\alpha }a^{n+m,l}_{\alpha +j,k}. 
$$
Then we have 
\begin{align*}
c^{n,m,l}_{1,j,k}&=b^{n,m,l}_{0,j,k},   \\
c^{n,m,l}_{i,j,k}-c^{n,m,l}_{i+1,j,k}&=b^{n,m,l}_{i,j,k}  \\
\intertext{for $1\leq i\leq n-1$ and Lem.5.21 implies}
c^{n,m,l}_{n,j,k}&=b^{n,m,l}_{n,j,k}. 
\end{align*}
We define a differential form $\Psi $ as 
$$
\Psi =\underset{1\leq k\leq l}{\underset{1\leq j\leq m}
{\underset{1\leq i\leq n}{\sum }}}c^{n,m,l}_{i,j,k}\CC _n(f)^{(-n+i-1,-i)}
\wedge \CC _m(g)^{(-m+j-1,-j)}\wedge \CC _l(h)^{(-l+k-1,-k)}.
$$

\vskip 1pc
\begin{lem} \ \ 
We have $c^{n,m,l}_{n-i+1,m-j+1,l-k+1}=c^{n,m,l}_{i,j,k}$.
Hence $\Psi $ is contained in \linebreak
$\mathscr D_{n+m+l+2}(X)$.
\end{lem}

{\it Proof}:
We have 
{\allowdisplaybreaks
\begin{align*}
c^{n,m,l}_{n-i+1,m-j+1,l-k+1}&=(-1)^ma^{m,l}_{m-j+1,l-k+1}  \\
&\hskip 1pc -2(-1)^m\sbinom{n+m}{n}^{-1}\sum_{\alpha =0}^{n-i}
\sbinom{n-\alpha +j-1}{n-\alpha }\sbinom{\alpha +m-j}{\alpha }
a^{n+m,l}_{\alpha +m-j+1,l-k+1}   \\
&=(-1)^{m+1}a^{m,l}_{j,k}+2(-1)^m\sbinom{n+m}{n}^{-1}
\sum_{\alpha =0}^{n-i}\sbinom{n-\alpha +j-1}{n-\alpha }
\sbinom{\alpha +m-j}{\alpha }a^{n+m,l}_{n-\alpha +j,k}   \\
&=(-1)^{m+1}a^{m,l}_{j,k}+2(-1)^m\sbinom{n+m}{n}^{-1}\sum_{\beta =i}^n
\sbinom{\beta +j-1}{\beta }\sbinom{n-\beta +m-j}{n-\beta }
a^{n+m,l}_{\beta +j,k}.  
\end{align*}}
Hence Lem.5.21 implies 
{\allowdisplaybreaks
\begin{align*}
c^{n,m,l}_{n-i+1,m-j+1,l-k+1}&=(-1)^{m+1}a^{m,l}_{j,k}+2(-1)^m
\left(a^{m,l}_{j,k}-\sbinom{n+m}{n}^{-1}\sum_{\beta =0}^{i-1}
\sbinom{\beta +j-1}{\beta }\sbinom{n-\beta +m-j}{n-\beta }
a^{n+m,l}_{\beta +j,k}\right)   \\
&=(-1)^ma^{m,l}_{j,k}-2(-1)^m\sbinom{n+m}{n}^{-1}\sum_{\beta =0}^{i-1}
\sbinom{n+m-\beta -j}{n-\beta }\sbinom{\beta +j-1}{\beta }
a^{n+m,l}_{\beta +j,k}   \\
&=c^{n,m,l}_{i,j,k},
\end{align*}}
which completes the proof.
\qed

\vskip 1pc
Let us denote the parts of $d\Psi $ including derivatives of 
$\CC _n(f), \CC _m(g)$ and $\CC _l(h)$ by $\Psi (f)$, 
$\Psi (g)$ and $\Psi (h)$ respectively. 
Then $d\Psi =\Psi (f)+\Psi (g)+\Psi (h)$ and 
{\allowdisplaybreaks
\begin{align*}
\Psi (f)&=\underset{1\leq k\leq l}{\underset{1\leq j\leq m}
{\underset{1\leq i\leq n}{\sum }}}c^{n,m,l}_{i,j,k}d\CC _n(f)^{(-n+i-1,-i)}
\wedge \CC _m(g)^{(-m+j-1,-j)}\wedge \CC _l(h)^{(-l+k-1,-k)}  \\
&=\underset{1\leq k\leq l}{\underset{1\leq j\leq m}{\sum }}
c^{n,m,l}_{1,j,k}\overline{\partial }\CC _n(f)^{(-n,-1)}
\wedge \CC _m(g)^{(-m+j-1,-j)}\wedge \CC _l(h)^{(-l+k-1,-k)}  \\*
&\hskip 1pc +\underset{1\leq k\leq l}{\underset{1\leq j\leq m}
{\underset{1\leq i\leq n-1}{\sum }}}(c^{n,m,l}_{i,j,k}-
c^{n,m,l}_{i+1,j,k})\partial \CC _n(f)^{(-n+i-1,-i)}
\wedge \CC _m(g)^{(-m+j-1,-j)}\wedge \CC _l(h)^{(-l+k-1,-k)}  \\*
&\hskip 1pc +\underset{1\leq k\leq l}{\underset{1\leq j\leq m}
{\sum }}c^{n,m,l}_{n,j,k}\partial \CC _n(f)^{(-1,-n)}
\wedge \CC _m(g)^{(-m+j-1,-j)}\wedge \CC _l(h)^{(-l+k-1,-k)}  \\
&=\Phi (f).
\end{align*}}
Let us express $\Phi (h)-\Psi (h)$ as follows:
\begin{align*}
&\Phi (h)-\Psi (h)  \\
&\hskip 1pc =d^{n,m,l}_{i,j,0}\underset{1\leq j\leq m}
{\underset{1\leq i\leq n}{\sum }}\CC _n(f)^{(-n+i-1,-i)}\wedge 
\CC _m(g)^{(-m+j-1,-j)}\wedge \overline{\partial }\CC _l(h)^{(-l,-1)} \\*
&\hskip 2pc +\underset{1\leq k\leq l}{\underset{1\leq j\leq m}
{\underset{1\leq i\leq n}{\sum }}}
d^{n,m,l}_{i,j,k}\CC _n(f)^{(-n+i-1,-i)}\wedge 
\CC _m(g)^{(-m+j-1,-j)}\wedge \partial \CC _l(h)^{(-l+k-1,-k)}. 
\end{align*}

\vskip 1pc
\begin{lem} \ \ 
It follows that $d^{n,m,l}_{i,j,k}=0$, therefore $\Phi (h)-\Psi (h)=0$. 
\end{lem}

{\it Proof}: \ 
When $1\leq k\leq l-1$, we have 
{\allowdisplaybreaks
\begin{align*}
d^{n,m,l}_{i,j,k}&=(-1)^na^{n,m+l}_{i,j+k}\times (a^{m,l}_{j,k}-
a^{m,l}_{j,k+1})-(-1)^{n+m}(c^{n,m,l}_{i,j,k}-c^{n,m,l}_{i,j,k+1})  \\
&=(-1)^{n+1}(1-a^{n,m+l}_{i,j+k})(a^{m,l}_{j,k}-a^{m,l}_{j,k+1})  \\*
&\hskip 1pc +2(-1)^n\sbinom{n+m}{n}^{-1}\sum_{\alpha =0}^{i-1}
\sbinom{n+m-\alpha -j}{n-\alpha }\sbinom{\alpha +j-1}{\alpha }
(a^{n+m,l}_{\alpha +j,k}-a^{n+m,l}_{\alpha +j,k+1})  \\
&=4(-1)^n\sbinom{n+m+l}{n}^{-1}\sbinom{m+l}{l}^{-1}
\sbinom{m+l-j-k}{m-j}\sbinom{j+k-1}{j-1}\sum_{\alpha =0}^{i-1}
\sbinom{n+m+l-i-j-k+1}{n-\alpha }\sbinom{i+j+k-1}{\alpha }  \\*
&\hskip 1pc -4(-1)^n\sbinom{n+m+l}{l}^{-1}\sbinom{n+m}{n}^{-1}
\sum_{\alpha =0}^{i-1}\sbinom{n+m+l-\alpha -j-k}{n+m-\alpha -j}
\sbinom{\alpha +j+k-1}{\alpha +j-1}\sbinom{n+m-\alpha -j}{n-\alpha }
\sbinom{\alpha +j-1}{\alpha }. 
\end{align*}}
Since
\begin{align*}
\sbinom{n+m+l}{l}\sbinom{n+m}{n}&=\sbinom{n+m+l}{n}\sbinom{m+l}{l},  \\
\sbinom{n+m+l-\alpha -j-k}{n+m-\alpha -j}\sbinom{n+m-\alpha -j}{n-\alpha }
&=\sbinom{n+m+l-\alpha -j-k}{m+l-j-k}\sbinom{m+l-j-k}{m-j},  \\
\sbinom{\alpha +j+k-1}{\alpha +j-1}\sbinom{\alpha +j-1}{\alpha }
&=\sbinom{\alpha +j+k-1}{\alpha }\sbinom{j+k-1}{j-1},
\end{align*} 
we have 
{\allowdisplaybreaks
\begin{align*}
d^{n,m,l}_{i,j,k}&=4(-1)^n\sbinom{n+m+l}{n}^{-1}\sbinom{m+l}{l}^{-1}
\sbinom{m+l-j-k}{m-j}\sbinom{j+k-1}{j-1}  \\*
&\hskip 2pc \times \left(\sum_{\alpha =0}^{i-1}
\sbinom{n+m+l-i-j-k+1}{n-\alpha }\sbinom{i+j+k-1}{\alpha }
-\sum_{\alpha =0}^{i-1}\sbinom{n+m+l-\alpha -j-k}{n-\alpha }
\sbinom{\alpha +j+k-1}{\alpha }\right)  \\
&=0
\end{align*}}
by Lem.A.3.

When $k=0$, we have 
{\allowdisplaybreaks
\begin{align*}
d^{n,m,l}_{i,j,0}&=(-1)^{n+1}a^{n,m}_{i,j}+(-1)^na^{n,m+l}_{i,j}\times 
(1+a^{m,l}_{j,1})-(-1)^{n+m}c^{n,m,l}_{i,j,1}  \\
&=(-1)^{n+1}\left((1-2\sbinom{n+m}{n}^{-1}\sum_{\alpha =0}^{i-1}
\sbinom{n+m-i-j+1}{n-\alpha }\sbinom{i+j-1}{\alpha }\right)  \\*
&\hskip 1pc +2(-1)^n\sbinom{m+l}{m}^{-1}\sbinom{m+l-j}{l}
\left(1-2\sbinom{n+m+l}{n}^{-1}\sum_{\alpha =0}^{i-1}
\sbinom{n+m+l-i-j+1}{n-\alpha }\sbinom{i+j-1}{\alpha }\right)  \\*
&\hskip 1pc -(-1)^n\left(-1+2\sbinom{m+l}{l}^{-1}\sbinom{m+l-j}{l}\right) \\*
&\hskip 1pc +2(-1)^n\sbinom{n+m}{n}^{-1}\sum_{\alpha =0}^{i-1}
\sbinom{n+m-\alpha -j}{n-\alpha }\sbinom{\alpha +j-1}{\alpha }
\left(-1+2\sbinom{n+m+l}{l}^{-1}\sbinom{n+m+l-\alpha -j}{l}\right).
\end{align*}} 
Since 
$$
\sbinom{n+m+l-\alpha -j}{l}\sbinom{n+m-\alpha -j}{n-\alpha }=
\sbinom{n+m+l-\alpha -j}{n-\alpha }\sbinom{m+l-j}{l},
$$
we have 
{\allowdisplaybreaks
\begin{align*}
d^{n,m,l}_{i,j,0}&=2(-1)^n\sbinom{n+m}{n}^{-1}\left(
\sum_{\alpha =0}^{i-1}\sbinom{n+m-i-j+1}{n-\alpha }\sbinom{i+j-1}{\alpha }
-\sum_{\alpha =0}^{i-1}\sbinom{n+m-\alpha -j}{n-\alpha }
\sbinom{\alpha +j-1}{\alpha }\right)  \\*
&\hskip 1pc +4(-1)^{n+1}\sbinom{n+m+l}{n}^{-1}\sbinom{m+l}{l}^{-1}
\sbinom{m+l-j}{l}  \\*
&\hskip 3pc \times \left(\sum_{\alpha =0}^{i-1}
\sbinom{n+m+l-i-j+1}{n-\alpha }\sbinom{i+j-1}{\alpha }-\sum_{\alpha =0}^{i-1}
\sbinom{n+m+l-\alpha -j}{n-\alpha }\sbinom{\alpha +j-1}{\alpha }\right)  \\
&=0
\end{align*}}
by Lem.A.3.
We can verify $d^{n,m,l}_{i,j,l}=0$ in the same way, so we omit the 
proof.
\qed

\vskip 1pc
We finally calculate $\Phi (g)-\Psi (g)$. 
Let us express it as follows:
\begin{align*}
&\Phi (g)-\Psi (g)  \\
&\hskip 1pc =e^{n,m,l}_{i,0,k}\underset{1\leq k\leq l}
{\underset{1\leq i\leq n}{\sum }}\CC _n(f)^{(-n+i-1,-i)}\wedge 
\overline{\partial }\CC _m(g)^{(-m,-1)}\wedge \CC _l(h)^{(-l+k-1,-k)} \\
&\hskip 2pc +\underset{1\leq k\leq l}{\underset{1\leq j\leq m}
{\underset{1\leq i\leq n}{\sum }}}
e^{n,m,l}_{i,j,k}\CC _n(f)^{(-n+i-1,-i)}\wedge 
\partial \CC _m(g)^{(-m+j-1,-j)}\wedge \CC _l(h)^{(-l+k-1,-k)}.
\end{align*}

\vskip 1pc
\begin{lem} \ \ 
It follows that $e^{n,m,l}_{i,j,k}=0$ if $1\leq j\leq m-1$ and 
$e^{n,m,l}_{i,0,k}=e^{n,m,l}_{i,m,k}=(-1)^{n+m+1}$.
\end{lem}

{\it Proof}: \ 
When $1\leq j\leq m-1$, we have
{\allowdisplaybreaks
\begin{align*}
e^{n,m,l}_{i,j,k}&=(-1)^{n+m+1}a^{n+m,l}_{i+j,k}\times (a^{n,m}_{i,j}
-a^{n,m}_{i,j+1})+(-1)^{n+m+1}a^{n,m+l}_{i,j+k}\times (a^{m,l}_{j,k}
-a^{m,l}_{j+1,k})  \\*
&\hskip 1pc +(-1)^n(c^{n,m,l}_{i,j,k}-c^{n,m,l}_{i,j+1,k})  \\
&=(-1)^{n+m+1}a^{n+m,l}_{i+j,k}\times (a^{n,m}_{i,j}
-a^{n,m}_{i,j+1})+(-1)^{n+m+1}(a^{n,m+l}_{i,j+k}-1)\times (a^{m,l}_{j,k}
-a^{m,l}_{j+1,k})  \\*
&\hskip 1pc +2(-1)^{n+m+1}\sbinom{n+m}{n}^{-1}\sum_{\alpha =0}^{i-1}
\sbinom{n+m-\alpha -j}{n-\alpha }\sbinom{\alpha +j-1}{\alpha }
a^{n+m,l}_{\alpha +j,k}  \\*
&\hskip 1pc -2(-1)^{n+m+1}\sbinom{n+m}{n}^{-1}\sum_{\alpha =0}^{i-1}
\sbinom{n+m-\alpha -j-1}{n-\alpha }\sbinom{\alpha +j}{\alpha }
a^{n+m,l}_{\alpha +j+1,k}.
\end{align*}}
Hence we have 
{\allowdisplaybreaks
\begin{align*}
&e^{n,m,l}_{i+1,j,k}-e^{n,m,l}_{i,j,k}  \\
&\hskip 1pc =(-1)^{n+m+1}a^{n+m,l}_{i+j+1,k}\times (a^{n,m}_{i+1,j}-
a^{n,m}_{i+1,j+1})-(-1)^{n+m+1}a^{n+m,l}_{i+j,k}\times (a^{n,m}_{i,j}
-a^{n,m}_{i,j+1})  \\*
&\hskip 2pc +(-1)^{n+m+1}(a^{n,m+l}_{i+1,j+k}-a^{n,m+l}_{i,j+k})
(a^{m,l}_{j,k}-a^{m,l}_{j+1,k})  \\*
&\hskip 2pc +2(-1)^{n+m+1}\sbinom{n+m}{n}^{-1}\left(\sbinom{n+m-i-j}{n-i}
\sbinom{i+j-1}{i}a^{n+m,l}_{i+j,k}-\sbinom{n+m-i-j-1}{n-i}\sbinom{i+j}{i}
a^{n+m,l}_{i+j+1,k}\right) \\
&\hskip 1pc =2(-1)^{n+m}\sbinom{n+m}{n}^{-1}\sbinom{n+m-i-j-1}{n-i-1}
\sbinom{i+j}{i}a^{n+m,l}_{i+j+1,k}-2(-1)^{n+m}\sbinom{n+m}{n}^{-1}
\sbinom{n+m-i-j}{n-i}\sbinom{i+j-1}{i-1}a^{n+m,l}_{i+j,k} \\*
&\hskip 2pc +(-1)^{n+m+1}(a^{n,m+l}_{i+1,j+k}-a^{n,m+l}_{i,j+k})
(a^{m,l}_{j,k}-a^{m,l}_{j+1,k})  \\*
&\hskip 2pc +2(-1)^{n+m+1}\sbinom{n+m}{n}^{-1}\left(\sbinom{n+m-i-j}{n-i}
\sbinom{i+j-1}{i}a^{n+m,l}_{i+j,k}-\sbinom{n+m-i-j-1}{n-i}\sbinom{i+j}{i}
a^{n+m,l}_{i+j+1,k}\right) \\
&\hskip 1pc =2(-1)^{n+m}\sbinom{n+m}{n}^{-1}\sbinom{n+m-i-j}{n-i}
\sbinom{i+j}{i}(a^{n+m,l}_{i+j+1,k}-a^{n+m,l}_{i+j,k}) \\*
&\hskip 2pc +(-1)^{n+m+1}(a^{n,m+l}_{i+1,j+k}-a^{n,m+l}_{i,j+k})
(a^{m,l}_{j,k}-a^{m,l}_{j+1,k})  \\
&\hskip 1pc =-4(-1)^{n+m}\sbinom{n+m}{n}^{-1}\sbinom{n+m-i-j}{n-i}
\sbinom{i+j}{i}\sbinom{n+m+l}{n+m}^{-1}\sbinom{n+m+l-i-j-k}{n+m-i-j}
\sbinom{i+j+k-1}{i+j} \\*
&\hskip 2pc +4(-1)^{n+m}\sbinom{n+m+l}{n}^{-1}\sbinom{n+m+l-i-j-k}{n-i}
\sbinom{i+j+k-1}{i}\sbinom{m+l}{m}^{-1}\sbinom{m+l-j-k}{m-j}
\sbinom{j+k-1}{j}  \\
&\hskip 1pc =0 
\end{align*}}
and 
{\allowdisplaybreaks
\begin{align*}
e^{n,m,l}_{1,j,k}&=(-1)^{n+m+1}a^{n+m,l}_{j+1,k}\times (a^{n,m}_{1,j}
-a^{n,m}_{1,j+1})+(-1)^{n+m+1}(a^{n,m+l}_{1,j+k}-1)(a^{m,l}_{j,k}
-a^{m,l}_{j+1,k}) \\*
&\hskip 1pc +2(-1)^{n+m+1}\sbinom{n+m}{n}^{-1}\left(\sbinom{n+m-j}{n}
a^{n+m,l}_{j,k}-\sbinom{n+m-j-1}{n}a^{n+m,l}_{j+1,k}\right)  \\
&=(-1)^{n+m+1}(a^{n,m+l}_{1,j+k}-1)(a^{m,l}_{j,k}-a^{m,l}_{j+1,k})
-2(-1)^{n+m}\sbinom{n+m}{n}^{-1}\sbinom{n+m-j}{n}(a^{n+m,l}_{j,k}-
a^{n+m,l}_{j+1,k})  \\
&=4(-1)^{n+m}\sbinom{n+m+l}{n}^{-1}\sbinom{n+m+l-j-k}{n}
\sbinom{m+l}{m}^{-1}\sbinom{m+l-j-k}{m-j}\sbinom{j+k-1}{j}  \\*
&\hskip 1pc -4(-1)^{n+m}\sbinom{n+m}{n}^{-1}\sbinom{n+m-j}{n}
\sbinom{n+m+l}{n+m}^{-1}\sbinom{n+m+l-j-k}{n+m-j}\sbinom{j+k-1}{j} \\
&=0.
\end{align*}}
Hence $e^{n,m,l}_{i,j,k}=0$ if $1\leq j\leq m-1$. 

When $j=0$, we have 
{\allowdisplaybreaks
\begin{align*}
e^{n,m,l}_{i,0,k}&=(-1)^{n+m+1}a^{n+m,l}_{i,k}\times (1+a^{n,m}_{i,1})
+(-1)^{n+m}a^{n,m+l}_{i,k}\times (1-a^{m,l}_{1,k})-(-1)^{n+1}
c^{n,m,l}_{i,1,k}  \\
&=2(-1)^{n+m+1}\sbinom{n+m}{n}^{-1}\sbinom{n+m-i}{m}a^{n+m,l}_{i,k}
+2(-1)^{n+m}\sbinom{m+l}{m}^{-1}\sbinom{m+l-k}{m}a^{n,m+l}_{i,k}  \\*
&\hskip 1pc +(-1)^{n+m}a^{m,l}_{1,k}-2(-1)^{n+m}\sbinom{n+m}{n}^{-1}
\sum_{\alpha =0}^{i-1}\sbinom{n+m-\alpha -1}{n-\alpha }
a^{n+m,l}_{\alpha +1,k}. 
\end{align*}}
Hence we have 
{\allowdisplaybreaks
\begin{align*}
&e^{n,m,l}_{i+1,0,k}-e^{n,m,l}_{i,0,k}  \\
&\hskip 1pc =-2(-1)^{n+m}\sbinom{n+m}{n}^{-1}\sbinom{n+m-i-1}{m}
a^{n+m,l}_{i+1,k}+2(-1)^{n+m}\sbinom{n+m}{n}^{-1}\sbinom{n+m-i}{m}
a^{n+m,l}_{i,k}  \\*
&\hskip 2pc +2(-1)^{n+m}\sbinom{m+l}{m}^{-1}\sbinom{m+l-k}{m}
(a^{n,m+l}_{i+1,k}-a^{n,m+l}_{i,k})-2(-1)^{n+m}\sbinom{n+m}{n}^{-1}
\sbinom{n+m-i-1}{n-i}a^{n+m,l}_{i+1,k}  \\
&\hskip 1pc =2(-1)^{n+m+1}\sbinom{n+m}{n}^{-1}\sbinom{n+m-i}{m}
(a^{n+m,l}_{i+1,k}-a^{n+m,l}_{i,k}) \\*
&\hskip 2pc +2(-1)^{n+m}\sbinom{m+l}{m}^{-1}\sbinom{m+l-k}{m}
(a^{n,m+l}_{i+1,k}-a^{n,m+l}_{i,k})  \\
&\hskip 1pc =4(-1)^{n+m}\sbinom{n+m}{n}^{-1}\sbinom{n+m-i}{m}
\sbinom{n+m+l}{n+m}^{-1}\sbinom{n+m+l-i-k}{n+m-i}\sbinom{i+k-1}{i}  \\*
&\hskip 2pc -4(-1)^{n+m}\sbinom{m+l}{m}^{-1}\sbinom{m+l-k}{m}
\sbinom{n+m+l}{n}^{-1}\sbinom{n+m+l-i-k}{n-i}\sbinom{i+k-1}{i} \\
&\hskip 1pc =0 
\end{align*}}
and 
{\allowdisplaybreaks
\begin{align*}
e^{n,m,l}_{1,0,k}&=2(-1)^{n+m+1}\sbinom{n+m}{n}^{-1}\sbinom{n+m-1}{m}
a^{n+m,l}_{1,k}+2(-1)^{n+m}\sbinom{m+l}{m}^{-1}\sbinom{m+l-k}{m}
a^{n,m+l}_{1,k}  \\*
&\hskip 1pc +(-1)^{n+m}a^{m,l}_{1,k}-2(-1)^{n+m}\sbinom{n+m}{n}^{-1}
\sbinom{n+m-1}{n}a^{n+m,l}_{1,k}  \\
&=2(-1)^{n+m+1}\left(1-2\sbinom{n+m+l}{n+m}^{-1}\sbinom{n+m+l-k}{n+m}
\right)+(-1)^{n+m}\left(1-2\sbinom{m+l}{m}^{-1}\sbinom{m+l-k}{m}\right) \\*
&\hskip 1pc +2(-1)^{n+m}\sbinom{m+l}{m}^{-1}\sbinom{m+l-k}{m}
\left(1-2\sbinom{n+m+l}{n}^{-1}\sbinom{n+m+l-k}{n}\right)  \\
&=(-1)^{n+m+1}+4(-1)^{n+m}\sbinom{n+m+l}{n+m}^{-1}\sbinom{n+m+l-k}{n+m} \\*
&\hskip 1pc -4(-1)^{n+m}\sbinom{m+l}{m}^{-1}\sbinom{m+l-k}{m}
\sbinom{n+m+l}{n}^{-1}\sbinom{n+m+l-k}{n} \\
&=(-1)^{n+m+1}. 
\end{align*}
Hence $e^{n,m,l}_{i,0,k}=(-1)^{n+m+1}$.
We can verify $e^{n,m,l}_{i,m,k}=(-1)^{n+m+1}$ in the same way, 
so we omit the proof. 
\qed

\vskip 1pc
Let us return to the proof of Prop.5.18.
By the above calculations, we have 
\begin{align*}
&\Phi -d\Psi   \\
&=(-1)^{n+m+1}\underset{1\leq k\leq l}{\underset{1\leq i\leq n}
{\sum }}\CC _n(f)^{(-n+i-1,-i)}\wedge \left(\overline{\partial }
\CC _m(g)^{(-m,-1)}+\partial \CC _m(g)^{(-1,-m)}\right)\wedge 
\CC _l(h)^{(-l+k-1,-k)}  \\
&=(-1)^{n+m+1}\CC _n(f)\wedge d\CC _m(g)\wedge \CC _l(h).
\end{align*}
Since $\Psi \in \mathscr D_{n+m+l+2}(X)$ and 
$d_{\mathscr D}\Psi =-d\Psi $, we have completed the proof.
\qed

\vskip 2pc
\section{Direct images}
\vskip 1pc

\subsection{Higher analytic torsion forms} \ \ 
We start this section by recalling the higher analytic torsion forms 
defined by Bismut and K\"{o}hler \cite{biskoh}.
First we fix some notations.

Let $\varphi :M\to N$ be a smooth projective morphism of compact 
complex algebraic manifolds. 
Let $T\varphi $ be the relative tangent bundle of $\varphi $ and 
we fix a metric $h_{\varphi }$ on $T\varphi $ that induces a K\"{a}hler 
metric on the fiber $\varphi ^{-1}(y)$ on each point $y\in N$. 
The pair $(\varphi , h_{\varphi })$ is called a 
{\it K\"ahler fibration}.
A real closed $(1,1)$-form $\Omega $ on $M$ is called a 
{\it K\"ahler form} with respect to $h_{\varphi }$ if 
the restriction of $\Omega $ to each fiber $\varphi ^{-1}(y)$ is 
a K\"{a}hler form with respect to $h_{\varphi }$.
Let $K_{\varphi }$ be the curvature form of 
$(T\varphi , h_{\varphi })$ and 
$$
\TD ^{GS}(\overline{T\varphi })=\TD \left(-
\frac{K_{\varphi }}{2\pi \sqrt{-1}}\right),
$$
where $\TD $ is the Todd polynomial. 

Let $\overline{E}$ be a $\varphi $-acyclic hermitian vector bundle 
on $M$, that is, $\overline{E}$ is a hermitian vector bundle on $M$ 
such that the higher direct image $R^i\varphi _*\overline{E}$ is 
trivial if $i>0$. 
Then the direct image $\varphi _*\overline{E}$ becomes a vector 
bundle and is equipped with the $L_2$-hermitian metric.
Let $\mathscr E^p$ be an infinite dimensional vector bundle on $N$ 
whose fiber on $y\in N$ is the vector space of smooth sections of 
$\Lambda ^pT^{* (1,0)}\varphi \otimes \overline{E}$ over 
$\varphi ^{-1}(y)$, where $T^{* (1,0)}\varphi $ is the holomorphic part 
of the complexified relative cotangent bundle 
$T^*\varphi \otimes \mathbb C$.
Then a hermitian metric on $\mathscr E^p$ can be defined by means of 
the K\"{a}hler metric $h_{\varphi }$.

Let us fix a K\"{a}hler form $\Omega $ with respect to $h_{\varphi }$.
Then the Bismut superconnection $B_u$ and the number operator $N_u$ 
for $u>0$ are defined on $\mathscr E=\underset{p}{\oplus }\mathscr E^p$.
The zeta function $\zeta (s)$ is defined as 
$$
\zeta (s)=\frac{-1}{2\Gamma (s)}\int_0^{\infty }u^{s-1}
\TTR _s(N_u\exp (-B_u^2))du, 
$$
where $\TTR_s$ is the supertrace.
The higher analytic torsion form of $\overline{E}$ 
is defined as 
$$
T(\overline{E}, \varphi , \Omega )=\Theta _1\left(\left.\frac{d}{ds}
\zeta (s)\right|_{s=0}\right),
$$
where $\Theta _1:\mathscr D_1(M)\to \mathscr A_1(M)$ is the 
renormalization operator defined in \S 1.6.
When the morphism $\varphi $ and the K\"{a}hler form $\Omega $ are 
specified, we sometimes abbreviate $T(\overline{E}, \varphi , \Omega )$ to 
$T(\overline{E})$.

The Grothendieck-Riemann-Roch theorem states that two closed differential 
forms \linebreak
$\int_{M/N}\CCGS (\overline{E})\TD ^{GS}(\overline{T\varphi })$ and 
$\CCGS (\varphi _*\overline{E})$ determine the same cohomology class.
The higher analytic torsion form of $T(\overline{E}, \varphi , \Omega )$ 
gives a homotopy between these two forms.
Namely, we have the following:

\vskip 1pc
\begin{thm}\cite{biskoh} \ \ 
Under the above notations, we have 
$$
d_{\mathscr A}T(\overline{E}, \varphi , \Omega )=
\CCGS _0(\varphi _*\overline{E})
-\int_{M/N}\TD ^{GS}(\overline{T\varphi })\CCGS _0(\overline{E}).
$$
\end{thm}

\vskip 1pc
\subsection{Independence of the analytic torsion forms} \ \ 
In this subsection we review independence of 
$T(\overline{E}, \varphi , \Omega )$ modulo $\IIm d_{\mathscr A}$ 
of the choice of the K\"{a}hler form $\Omega $, 
due also to Bismut and K\"{o}hler \cite{biskoh}.
Let $\Omega $ and $\Omega ^{\prime }$ be two K\"{a}hler forms with 
respect to $h_{\varphi }$.
For $0\leq l\leq 1$, let $\Omega _l=l\Omega +(1-l)\Omega ^{\prime }$.
Then $\Omega _l$ is also a  K\"{a}hler form with respect to 
$h_{\varphi }$.
The Bismut superconnection and the number operator determined by 
$\Omega _l$ are denoted by $B_{u,l}$ and $N_{u,l}$ respectively.
The holomorphic and anti-holomorphic parts of $B_{u,l}$ are denoted by 
$B_{u,l}^{\prime }$ and $B_{u,l}^{\prime \prime }$ respectively.
In addition, let $M_u$ be the section of 
$\Lambda ^*T^*N\otimes \END (\mathscr E)$ defined in 
\cite[Def.2.7]{biskoh}.
The zeta function with respect to $\Omega _l$ is defined as 
$$
\zeta _l(s)=\frac{-1}{2\Gamma (s)}\int_0^{\infty }u^{s-1}
\TTR _s(N_{u,l}\exp (-B_{u,l}^2))du.
$$
By Thm.2.9 in \cite{biskoh}, we have 
{\allowdisplaybreaks
\begin{align*}
\frac{d}{dl}\TTR _s(N_{u,l}\exp (-B_{u,l}^2))&=u\frac{d}{du}(M_u
\exp (-B_{u,l}^2))  \\
&\hskip 1pc -\overline{\partial }\left(\left.\frac{d}{db}\TTR _s
([B_{u,l}^{\prime }, N_{u,l}]\exp (-B_{u,l}^2-bM_u))\right|_{b=0}\right)  \\
&\hskip 1pc -\partial \left(\left.\frac{d}{db}\TTR _s
([B_{u,l}^{\prime \prime }, N_{u,l}]\exp (-B_{u,l}^2-bM_u))\right|_{b=0}
\right)  \\
&\hskip 1pc -\overline{\partial }\partial \left(\left.\frac{d}{db}\TTR _s
(N_{u,l}\exp (-B_{u,l}^2-bM_u))\right|_{b=0}\right).
\end{align*}}
We set 
{\allowdisplaybreaks
\begin{align*}
\theta _0(s)&=\frac{-1}{2\Gamma (s)}\int_0^{\infty }u^s\frac{d}{du}
\TTR _s(M_u\exp (-B_{u,l}^2))du,  \\
\theta _1(s)&=\frac{-1}{2\Gamma (s)}\int_0^{\infty }u^{s-1}\left(\left.
\frac{d}{db}\TTR _s([B_{u,l}^{\prime }, N_{u,l}]\exp (-B_{u,l}^2
-bM_u))\right|_{b=0}\right)du,  \\
\theta _2(s)&=\frac{-1}{2\Gamma (s)}\int_0^{\infty }u^{s-1}\left(\left.
\frac{d}{db}\TTR _s([B_{u,l}^{\prime \prime }, N_{u,l}]\exp (-B_{u,l}^2
-bM_u))\right|_{b=0}\right)du,  \\
\theta _3(s)&=\frac{-1}{2\Gamma (s)}\int_0^{\infty }u^{s-1}\left(\left.
\frac{d}{db}\TTR _s(N_{u,l}\exp (-B_{u,l}^2-bM_u))\right|_{b=0}\right)du.
\end{align*}}
These integrals converge when $\RE (s)$ is sufficiently big and are 
meromorphically prolongable to the whole plane.
If we differentiate $\zeta _l(s)$ with the parameter $l$, then we have 
$$
\frac{d}{dl}\zeta _l(s)=\theta _0(s)-\overline{\partial }\theta _1(s)
-\partial \theta _2(s)-\overline{\partial }\partial \theta _3(s).
$$
It follows from the definition of $\theta _0(s)$ that 
$\theta _0^{\prime }(0)$ is equal to the constant term of the power 
series expansion of $\TTR _s(M_u\exp (-B_{u,l}^2))$ at $u=0$.
By Thm.3.17 and Thm.3.22 in \cite{biskoh}, we have 
$$
\theta _0^{\prime }(0)=-d\left(\left.\frac{d}{db}\TTR _s(-uM_u
\exp (-B_{u,l}^2-b\frac{dB_{u,l}}{du}))\right|_{b=0}\right).
$$

The same argument as the proof of Thm.2.2 in \cite{bisgs} yields the 
following:
\begin{gather*}
\left.\frac{d}{db}\TTR _s(-uM_u\exp (-B_{u,l}^2-b\frac{dB_{u,l}}{du}))
\right|_{b=0}\in \underset{p}{\oplus }(\mathscr E^{p-2,p-1}(N)\oplus 
\mathscr E^{p-1,p-2}(N)), \\
\theta _1^{\prime }(0)\in \underset{p}{\oplus }\mathscr E^{p-1,p-2}(N),  \\
\theta _2^{\prime }(0)\in \underset{p}{\oplus }\mathscr E^{p-2,p-1}(N),  \\
\theta _3^{\prime }(0)\in \underset{p}{\oplus }\mathscr E^{p-2,p-2}(N).
\end{gather*}
If $\pi :\underset{p}{\oplus }(\mathscr E^{p-2,p-1}(N)\oplus 
\mathscr E^{p-1,p-2}(N))\to \mathscr D_2(N)$ denotes the canonical 
projection and 
$$
\beta (l)=\pi \left(\left.\frac{d}{db}\TTR _s(-uM_u
\exp (-B_{u,l}^2-b\frac{dB_{u,l}}{du}))\right|_{b=0}
+\theta _1^{\prime }(0)+\theta _2^{\prime }(0)
+\frac{\partial -\overline{\partial }}{2}\theta _3^{\prime }(0)\right)
\in \mathscr D_2(N),
$$
then we have 
$$
d_{\mathscr D}\beta (l)=\frac{d}{dl}\zeta _l^{\prime }(0).
$$
Hence if we set 
$$
\mu (\overline{E}, \Omega , \Omega ^{\prime })=\Theta _2\left(
\int_0^1\beta (l)dl\right)\in \mathscr A_2(N), 
$$
where $\Theta _2:\mathscr D_2(N)\to \mathscr A_2(N)$ is the renormalization 
operator defined in \S 1.6, then we have 
{\allowdisplaybreaks
\begin{align*}
d_{\mathscr A}\mu (\overline{E}, \Omega , \Omega ^{\prime })
&=\Theta _1\left(\int_0^1d_{\mathscr D}\beta (l)dl\right)  \\
&=\Theta _1\left(\int_0^1\frac{d}{dl}\zeta _l^{\prime }(0)dl\right)  \\
&=\Theta _1\left(\zeta _1^{\prime }(0)-\zeta _0^{\prime }(0)\right)  \\
&=T(\overline{E}, \varphi, \Omega )-
T(\overline{E}, \varphi, \Omega ^{\prime }).
\end{align*}}

We next discuss compatibility of $T(\overline{E}, \varphi , \Omega )$ 
and $\mu (\overline{E}, \Omega , \Omega ^{\prime })$ with the 
pull back for closed immersions. 
Consider the following cartesian square:
$$
\begin{CD}
M^{\prime } @>{j}>>  M \\
@VV{\varphi ^{\prime }}V @VV{\varphi }V  \\
N^{\prime } @>{i}>> N,
\end{CD}
$$
where $i$ and $j$ are closed immersions and $\varphi $ is 
a K\"{a}hler fibration with a metric $h_{\varphi }$ on $T\varphi $.
Then it follows that $T\varphi ^{\prime }\simeq j^*T\varphi $, 
therefore a hermitian metric $h_{\varphi ^{\prime }}$ on 
$T\varphi ^{\prime }$ with which $\varphi ^{\prime }$ becomes 
a K\"{a}hler fibration is induced from $h_{\varphi }$.
If $\Omega $ is a K\"{a}hler form with respect to $h_{\varphi }$, 
then the pull back $j^*\Omega $ is a K\"{a}hler form 
with respect to $h_{\varphi ^{\prime }}$. 

Let $\overline{E}$ be a $\varphi $-acyclic hermitian vector bundle 
on $M$.
Then it is obvious that $B_{u,l}$, $N_{u,l}$ and $M_u$ for 
$\overline{E}$ as mentioned above are compatible with the pull back for 
the immersions $i$ and $j$.
Hence the forms $T(\overline{E}, \varphi , \Omega )$ 
and $\mu (\overline{E}, \Omega , \Omega ^{\prime })$ admit 
the naturality for $i$ and $j$.

Summarizing results obtained so far, we have the following:

\vskip 1pc
\begin{prop} \ \ 
Let $\varphi :M\to N$ and $h_{\varphi }$ be as in the last subsection.
Let $\Omega $ and $\Omega ^{\prime }$ be K\"{a}hler forms with respect 
to $h_{\varphi }$.
Let $\overline{E}$ be a $\varphi $-acyclic hermitian vector bundle 
on $M$.
Then there is a form 
$\mu (\overline{E}, \Omega , \Omega ^{\prime })\in \mathscr A_2(N)$ 
such that 
$$
d_{\mathscr A}\mu (\overline{E}, \Omega , \Omega ^{\prime })
=T(\overline{E}, \varphi, \Omega )-
T(\overline{E}, \varphi, \Omega ^{\prime }).
$$
Furthermore, for the cartesian square as above, we have 
\begin{align*}
i^*T(\overline{E}, \varphi , \Omega )&=
T(j^*\overline{E}, \varphi ^{\prime }, j^*\Omega ),  \\*
i^*\mu (\overline{E}, \Omega , \Omega ^{\prime })&=
\mu (j^*\overline{E}, j^*\Omega , j^*\Omega ^{\prime }).
\end{align*}
\end{prop}

\vskip 1pc
\subsection{Higher analytic torsion forms for cubes} \ \ 
In this subsection we introduce the higher analytic torsion form of an 
exact metrized $n$-cube defined by Roessler \cite{roessler}. 
To do this we need another presentation of higher Bott-Chern forms. 

Let $M$ be a compact complex algebraic manifold.
For $u_i\in \mathscr D_1(M)$, 
we set 
$$
C_n(u_1, \cdots , u_n)=\frac{1}{2^n}\sum_{\sigma \in \mathfrak S_n}
(-1)^{\SGN \sigma }u_{\sigma (1)}\bullet \left(u_{\sigma (2)}\bullet 
\left(\cdots u_{\sigma (k)}\cdots \right)\right).
$$
Similarly, for $u_i\in \mathscr A_1(M)$, we set 
$$
C^{GS}_n(u_1, \cdots , u_n)=\sum_{\sigma \in \mathfrak S_n}
(-1)^{\SGN \sigma }u_{\sigma (1)}\bullet \left(u_{\sigma (2)}\bullet 
\left(\cdots u_{\sigma (k)}\cdots \right)\right), 
$$
where $\bullet $ is the multiplication on $\mathscr A_*(M)$ defined 
in \S 4.3.

\vskip 1pc
\begin{prop} \ \ 
It follows that  
$$
C_n(u_1, \cdots , u_n)=\frac{(-1)^n}{2}\sum_{i=1}^n(-1)^i
S_n^i(u_1, \cdots , u_n).
$$
Hence for an exact metrized $n$-cube $\mathcal F$, we have 
$$
\CC _n(\mathcal F)=\frac{1}{(2\pi \sqrt{-1})^nn!}\int_{(\mathbb P^1)^n}
\CC _0(\TTR _n\lambda \mathcal F)C_n(\log |z_1|^2, \cdots , \log |z_n|^2)
$$
and 
\begin{align*}
\CCGS _n(\mathcal F)&=\frac{(-1)^n}{(2\pi \sqrt{-1})^{n-1}n!}
\int_{(\mathbb P^1)^n}\CCGS _0(\TTR _n\lambda \mathcal F)
\sum_{i=1}^n(-1)^iS_n^i  \\
&=\frac{1}{n!}\int_{(\mathbb P^1)^n}
\CCGS _0(\TTR _n\lambda \mathcal F)C^{GS}_n
(\log |z_1|^2, \cdots , \log |z_n|^2).
\end{align*}
\end{prop}

\vskip 1pc
Let $\varphi :M\to N$ be a smooth projective morphism of 
compact complex algebraic manifolds and $h_{\varphi }$ a hermitian 
metric on $T\varphi $ such that $(\varphi , h_{\varphi })$ is 
a K\"{a}hler fibration.
We assume that any exact metrized $n$-cube on $M$ we deal with is 
made of $\varphi $-acyclic hermitian vector bundles. 
Then $\lambda \mathcal F$ for such a cube $\mathcal F$ is also made 
of $\varphi $-acyclic vector bundles and there is a canonical isomorphism 
$$
\varphi _*(\TTR _n\lambda \mathcal F)\simeq 
\lambda \TTR _n\varphi _*\mathcal F.
$$
When we put the $L_2$-metrics on the both sides, however, this isomorphism 
does not preserve the metrics.
In \cite[\S 3.1]{roessler}, Roessler has constructed a hermitian vector 
bundle $\overline{h(\mathcal F)}$ connecting these metrics.
Namely, $\overline{h(\mathcal F)}$ is a hermitian vector bundle on 
$N\times (\mathbb P^1)^{n+1}$ satisfying the following conditions: 
$$
\overline{h(\mathcal F)}\vert _{X\times \{0\}\times (\mathbb P^1)^n}=
\varphi _*(\TTR _n\lambda \mathcal F), \ \ 
\overline{h(\mathcal F)}\vert _{X\times \{\infty \}\times (\mathbb P^1)^n}=
\lambda \TTR _n\varphi _*\mathcal F
$$
and 
$$
\overline{h(\mathcal F)}\vert _{X\times (\mathbb P^1)^i\times \{0\}\times 
(\mathbb P^1)^{n-i}}=\overline{h(\partial _i^0\mathcal F)}, \ \ 
\overline{h(\mathcal F)}\vert _{X\times (\mathbb P^1)^i\times \{\infty \}
\times (\mathbb P^1)^{n-i}}=\overline{h(\partial _i^{-1}\mathcal F)}\oplus 
\overline{h(\partial _i^1\mathcal F)}
$$
for $1\leq i\leq n$.
For the precise definition of $\overline{h(\mathcal F)}$, 
see \cite{roessler}.
Let us define a differential form 
$T_1(\mathcal F, \varphi )\in \mathscr A_{n+1}(N)$ as 
$$
T_1(\mathcal F, \varphi )=\frac{(-1)^n}{(2\pi \sqrt{-1})^n(n+1)!}
\int_{(\mathbb P^1)^{n+1}}\CCGS _0(\overline{h(\mathcal F)})
\sum_{i=1}^{n+1}(-1)^iS_{n+1}^i.
$$
Moreover we take a K\"{a}hler form $\Omega $ with respect to $h_{\varphi }$ 
and define $T_2(\mathcal F, \varphi , \Omega )\in \mathscr A_{n+1}(N)$ 
as 
$$
T_2(\mathcal F, \varphi , \Omega )=
\frac{(-1)^{n+1}}{(2\pi \sqrt{-1})^n(n+1)!}
\int_{(\mathbb P^1)^n}\sum_{i=1}^{n+1}(-1)^iS_{n+1}^i(\mathcal F)
$$
where 
$$
S_{n+1}^i(\mathcal F)=S_{n+1}^i
(T(\TTR _n\lambda \mathcal F, \varphi , \Omega ), 
\log |z_1|^2, \cdots , \log |z_n|^2).
$$

\vskip 1pc
\begin{thm}\cite[Thm.3.6]{roessler} \ \ 
For an exact metrized $n$-cube $\mathcal F$ on $M$, we have 
$$
d_{\mathscr A}T_1(\mathcal F, \varphi )+
T_1(\partial \mathcal F, \varphi )=\CCGS _n(\varphi _*\mathcal F)
-\frac{(-1)^n}{(2\pi \sqrt{-1})^{n-1}n!}
\int_{(\mathbb P^1)^n}\CCGS _0(\varphi _*\TTR _n\lambda \mathcal F)
\sum_{i=1}^n(-1)^iS_n^i
$$
and 
\begin{align*}
d_{\mathscr A}T_2(\mathcal F, \varphi , \Omega )+
T_2(\partial \mathcal F, \varphi , \Omega )&=
\frac{(-1)^n}{(2\pi \sqrt{-1})^{n-1}n!}\int_{(\mathbb P^1)^n}
\CCGS _0(\varphi _*\TTR _n\lambda \mathcal F)\sum_{i=1}^n(-1)^iS_n^i \\
&\hskip 1pc -\int_{M/N}\TD ^{GS}(\overline{T\varphi })\CCGS _n(\mathcal F).
\end{align*}
Hence if we set $T(\mathcal F, \varphi , \Omega )=
T_1(\mathcal F, \varphi )+T_2(\mathcal F, \varphi , \Omega )$, 
then we have 
$$
d_{\mathscr A}T(\mathcal F, \varphi , \Omega )+
T(\partial \mathcal F, \varphi , \Omega )=\CCGS _n(\varphi _*\mathcal F)-
\int_{M/N}\TD ^{GS}(\overline{T\varphi })\CCGS _n(\mathcal F).
$$
\end{thm} 

\vskip 1pc
Let us discuss dependence of $T(\mathcal F, \varphi , \Omega )$ 
on the K\"{a}hler form $\Omega $.
We begin with the following lemma.

\vskip 1pc
\begin{lem} \ \ 
For $u_1\in \mathscr A_2(M)$ and $u_i\in \mathscr A_1(M)$ with 
$2\leq i\leq n$, let 
$$
C^{GS}_n(u_1, \cdots , u_n)=\sum_{j=1}^n(-1)^{j+1}
\underset{\sigma (j)=1}{\underset{\sigma \in \mathfrak S_n}
{\sum }}(-1)^{\SGN \sigma }u_{\sigma (1)}\bullet (u_{\sigma (2)}\bullet 
(\cdots u_{\sigma (j)}\cdots )).
$$
Then we have 
\begin{align*}
d_{\mathscr A}&C^{GS}_n(u_1, u_2, \cdots , u_n) \\*
&=C^{GS}_n(d_{\mathscr A}u_1, u_2, \cdots , u_n)+n\sum_{k=2}^n(-1)^k
(d_{\mathscr A}u_k)\bullet C^{GS}_{n-1}(u_1, u_2, \cdots , \widehat{u_k}, 
\cdots , u_n).
\end{align*}
\end{lem}

{\it Proof}: \ 
Since $d_{\mathscr A}(u\bullet v)=d_{\mathscr A}u\bullet v+(-1)^{\DEG u}
u\bullet d_{\mathscr A}v$, we have 
{\allowdisplaybreaks
\begin{align*}
d&_{\mathscr A}C^{GS}_n(u_1, u_2, \cdots , u_n) \\
&=\sum_{j=1}^n
\underset{\sigma (j)=1}{\underset{\sigma \in \mathfrak S_n}{\sum }}
(-1)^{\SGN \sigma }\underset{i<j}{\sum }(-1)^{i+j}d_{\mathscr A}
u_{\sigma (i)}(u_{\sigma (1)}\bullet 
(\cdots \widehat{u_{\sigma (i)}}\cdots u_{\sigma (j)}\cdots ))  \\*
&\hskip 1pc +\sum_{j=1}^n
\underset{\sigma (j)=1}{\underset{\sigma \in \mathfrak S_n}{\sum }}
(-1)^{\SGN \sigma }u_{\sigma (1)}\bullet (\cdots 
d_{\mathscr A}u_{\sigma (j)}\cdots )  \\*
&\hskip 1pc +\sum_{j=1}^n
\underset{\sigma (j)=1}{\underset{\sigma \in \mathfrak S_n}{\sum }}
(-1)^{\SGN \sigma }\underset{j<i}{\sum }(-1)^{i+j+1}d_{\mathscr A}
u_{\sigma (i)}(u_{\sigma (1)}\bullet (\cdots u_{\sigma (j)}\cdots 
\widehat{u_{\sigma (i)}}\cdots ))  \\
&=\sum_{j=1}^n
\underset{\sigma (j)=1}{\underset{\sigma \in \mathfrak S_n}{\sum }}
(-1)^{\SGN \sigma }u_{\sigma (1)}\bullet (\cdots 
d_{\mathscr A}u_{\sigma (j)}\cdots )  \\*
&\hskip 1pc +\sum_{k=2}^n\sum_{i<j}
\underset{\sigma (i)=k}{\underset{\sigma (j)=1}
{\underset{\sigma \in \mathfrak S_n}{\sum }}}(-1)^{\SGN \sigma }(-1)^{i+j}
d_{\mathscr A}u_k(u_{\sigma (1)}\bullet (\cdots \widehat{u_{\sigma (i)}}
\cdots u_{\sigma (j)}\cdots ))  \\*
&\hskip 1pc +\sum_{k=2}^n\sum_{j<i}
\underset{\sigma (i)=k}{\underset{\sigma (j)=1}
{\underset{\sigma \in \mathfrak S_n}{\sum }}}(-1)^{\SGN \sigma }(-1)^{i+j+1}
d_{\mathscr A}u_k(u_{\sigma (1)}\bullet (\cdots u_{\sigma (j)}
\cdots \widehat{u_{\sigma (i)}}\cdots ))  \\
&=C_{n-1}(d_{\mathscr A}u_1, u_2, \cdots , u_n)+n\sum_{k=2}^n(-1)^k
(d_{\mathscr A}u_k)C_{n-1}(u_1, \cdots , \widehat{u_k}, \cdots , u_n),
\end{align*}}
which completes the proof.
\qed

\vskip 1pc
\begin{prop} \ \ 
For an exact metrized $n$-cube $\mathcal F$ on $M$, let 
$\mu (\mathcal F)=\mu (\TTR _n\lambda \mathcal F, \Omega , 
\Omega ^{\prime })$. 
Then we have 
$$
T(\mathcal F, \varphi , \Omega )-T(\mathcal F, \varphi , \Omega ^{\prime })
\equiv \frac{-1}{n!}\int_{(\mathbb P^1)^{n-1}}C^{GS}_n
(\mu (\partial \mathcal F), \log |z_1|^2, \cdots , \log |z_{n-1}|^2))
$$
modulo $\IIm d_{\mathscr A}$. 
\end{prop}

{\it Proof}: \ 
Prop.6.3 implies 
$$
T_2(\mathcal F, \varphi , \Omega )=\frac{1}{(n+1)!}\int_{(\mathbb P^1)^n}
C_{n+1}^{GS}(T(\TTR _n\lambda \mathcal F, \varphi , \Omega ), 
\log |z_1|^2, \cdots , \log |z_n|^2).
$$
Since $T_1(\mathcal F, \varphi )$ does not depend on $\Omega $, 
by Lem.6.5 and Prop.6.2 we have 
{\allowdisplaybreaks
\begin{align*}
T(\mathcal F, &\varphi , \Omega )-T(\mathcal F, \varphi , 
\Omega ^{\prime })=T_2(\mathcal F, \varphi , \Omega )-
T_2(\mathcal F, \varphi , \Omega ^{\prime })  \\
&=\frac{1}{(n+1)!}\int_{(\mathbb P^1)^n}
C^{GS}_{n+1}(T(\TTR _n\lambda \mathcal F, \varphi , \Omega )-
T(\TTR _n\lambda \mathcal F, \varphi , \Omega ^{\prime }), 
\log |z_1|^2, \cdots , \log |z_n|^2)   \\
&=\frac{1}{(n+1)!}\int_{(\mathbb P^1)^n}
C^{GS}_{n+1}(d_{\mathscr A}\mu (\mathcal F), \log |z_1|^2, 
\cdots , \log |z_n|^2)   \\
&=\frac{1}{(n+1)!}\int_{(\mathbb P^1)^n}d_{\mathscr A}
C^{GS}_{n+1}(\mu (\mathcal F), \log |z_1|^2, 
\cdots , \log |z_n|^2)   \\*
&\hskip 1pc -\frac{1}{n!}\sum_{k=1}^n(-1)^{k-1}\int_{(\mathbb P^1)^n}
d_{\mathscr A}\log |z_k|^2C^{GS}_n(\mu (\mathcal F), 
\cdots , \widehat{\log |z_k|^2}, \cdots )   \\
&=\frac{1}{(n+1)!}d_{\mathscr A}\left(\int_{(\mathbb P^1)^n}
C^{GS}_{n+1}(\mu (\mathcal F), \log |z_1|^2, 
\cdots , \log |z_n|^2)\right)   \\*
&\hskip 1pc -\frac{1}{n!}\int_{(\mathbb P^1)^{n-1}}
C^{GS}_n(\mu (\partial \mathcal F), \log |z_1|^2, \cdots , 
\log |z_{n-1}|^2), 
\end{align*}}
which completes the proof.
\qed

\vskip 1pc
\subsection{Definition of direct images} \ \ 
In this subsection, we apply the results obtained so far in this 
section to the arithmetic situation and define a direct image 
morphism in higher arithmetic $K$-theory.
Let $\varphi :X\to Y$ be a smooth projective morphism of proper 
arithmetic varieties. 
We fix an $F_{\infty }$-invariant metric $h_{\varphi }$ on 
$T\varphi (\mathbb C)$ and take an anti-$F_{\infty }$-invariant 
K\"{a}hler form $\Omega $ on $X(\mathbb C)$ with respect to $h_{\varphi }$.
Let $\widehat{S}(\ACYC )$ denote the S-construction of the category 
of $\varphi $-acyclic hermitian vector bundles on $X$. 
Then the direct image of a $\varphi $-acyclic hermitian vector bundle 
with the $L_2$-metric gives a morphism of simplicial sets 
$$
\varphi _*:\widehat{S}(\ACYC )\to \widehat{S}(Y).
$$
Since the natural inclusion $\widehat{S}(\ACYC )\to \widehat{S}(X)$ 
is homotopy equivalent, we have 
$$
\widehat{\pi }_*(\vert \widehat{S}(\ACYC )\vert , \CCGS )\simeq 
\widehat{\pi }_*(\vert \widehat{S}(X)\vert , \CCGS ).
$$

\vskip 1pc
\begin{prop} \ \ 
If $E$ is a degenerate element of $\widehat{S}_{n+1}(\ACYC )$, then 
we have \linebreak
$T(\CUB (E), \varphi ,\Omega )=0$.
\end{prop}

\vskip 1pc
The proof is similar to that of Thm.3.4, so we omit it.
By Prop.6.7, taking the higher analytic torsion forms yields 
the homomorphism 
$$
T( \ \ , \varphi , \Omega ):C_*(\vert \widehat{S}(\ACYC )\vert )\to 
\mathscr A_*(Y).
$$
In particular, the higher analytic torsion form of a pointed cellular map 
$f:S^{n+1}\to \vert \widehat{S}(\ACYC )\vert $ is defined by 
$T(f, \varphi , \Omega )=T(f_*([S^{n+1}]), \varphi , \Omega )$.
To simplify notation, we write $T(f)$ for $T(f, \varphi , \Omega )$ 
if the morphism $\varphi $ and the K\"{a}hler form $\Omega $ 
are specified.

Let us define a homomorphism of complexes 
$\varphi _!:\mathscr A_*(X)\to \mathscr A_*(Y)$ by 
$$
\varphi _!\omega =\int_{X(\mathbb C)/Y(\mathbb C)}\TD ^{GS}
(\overline{T\varphi })\omega .
$$
Prop.6.4 tells that the diagram 
$$
\begin{CD}
C_*(\vert \widehat{S}(\ACYC )\vert ) @>{\varphi _*}>> 
C_*(\vert \widehat{S}(Y)\vert )  \\
@VV{\CCGS }V  @VV{\CCGS }V  \\
\mathscr A_*(X)[1]  @>{\varphi _!}>> \mathscr A_*(Y)[1]
\end{CD}
$$
is commutative up to the homotopy $T( \ \ , \varphi , \Omega )$.
Hence Prop.2.7 yields the following:

\vskip 1pc
\begin{prop} \ \ 
We can define a homomorphism 
$$
\widehat{\varphi }(\Omega )_*:\widehat{\pi }_{n+1}(\vert \widehat{S}(\ACYC )\vert , 
\CCGS )\to \widehat{\pi }_{n+1}(\vert \widehat{S}(Y)\vert , \CCGS )
$$
by $[(f, \omega )]\mapsto [(\varphi _*f, \varphi _!\omega 
+T(f, \varphi , \Omega ))]$.
\end{prop}

\vskip 1pc
If $\Omega ^{\prime }$ is another anti-$F_{\infty }$-invariant K\"{a}hler 
form with respect to $h_{\varphi }$, then it follows from Prop.6.6 that 
$T(f, \varphi , \Omega )\equiv T(f, \varphi , \Omega ^{\prime })$ 
modulo $\IIm d_{\mathscr A}$ for any 
pointed cellular map $f:S^{n+1}\to \vert \widehat{S}(\ACYC )\vert $. 
Hence the homomorphism $\widehat{\varphi }(\Omega )_*$ depends 
only on the hermitian metric $h_{\varphi }$ and does not concern 
the K\"{a}hler form $\Omega $.

The following is the main theorem of this section.

\vskip 1pc
\begin{thm} \ \ 
Let $\varphi :X\to Y$ be a smooth projective morphism of proper 
arithmetic varieties. 
We fix an $F_{\infty }$-invariant metric $h_{\varphi }$ on $T\varphi $ 
such that the restriction of $h_{\varphi }$ to any fiber $\varphi ^{-1}(y)$ 
on $y\in Y(\mathbb C)$ is a K\"{a}hler metric. 
Then we can define a direct image morphism 
$$
\widehat{\varphi }(h_{\varphi })_*:\widehat{K}_n(X)\to \widehat{K}_n(Y)
$$
by 
$$
\widehat{\pi }_{n+1}(\vert \widehat{S}(X)\vert , \CCGS )\simeq 
\widehat{\pi }_{n+1}(\vert \widehat{S}(\ACYC )\vert , \CCGS )
\overset{\widehat{\varphi }(\Omega )_*}{\longrightarrow }
\widehat{\pi }_{n+1}(\vert \widehat{S}(Y)\vert , \CCGS ),
$$
where $\Omega $ is an anti-$F_{\infty }$-invariant K\"{a}hler from on 
$X(\mathbb C)$ with respect to $h_{\varphi }$.
\end{thm}

\vskip 1pc
When $n=0$, the direct image morphism we have defined 
above agrees with the morphism $\varphi _!$ in \cite{gilsoulz}.
In fact, we have the commutative diagram 
$$
\begin{CD}
\widehat{\mathscr K}_0(X) @>{\widehat{\alpha }}>> \widehat{K}_0(X)  \\
@VV{\varphi _!}V  @VV{\widehat{\varphi }(h_{\varphi })_*}V  \\
\widehat{\mathscr K}_0(Y) @>{\widehat{\alpha }}>> \widehat{K}_0(Y), 
\end{CD}
$$
where $\varphi _!$ is defined as 
$$
\varphi _!(\overline{E}, \omega )=(\varphi _*\overline{E}, 
\varphi _!(\omega )-T(\overline{E}))
$$
for a $\varphi $-acyclic hermitian vector bundle $\overline{E}$ on $X$ 
and $\omega \in \widetilde{\mathscr A}_1(X)$.

Prop.2.8 implies that the Chern form map $\CCGS _n$ is 
compatible with the direct image morphism, that is, the diagram 
$$
\begin{CD}
\widehat{K}_n(X) @>{\CCGS _n}>> \mathscr A_*(X)  \\
@VV{\widehat{\varphi }(h_{\varphi })_*}V  @VV{\varphi _!}V  \\
\widehat{K}_n(Y) @>{\CCGS _n}>> \mathscr A_*(Y)
\end{CD}
$$
is commutative.
In particular, we can define a direct image morphism in $KM$-groups 
$$
\widehat{\varphi }(h_{\varphi })_*:KM_n(X)\to KM_n(Y).
$$

Finally, we give a description of the direct image morphism 
by means of the $G$-construction.
Given a pointed cellular map $f:S^n\to \vert \widehat{G}(X)\vert $ 
for $n\geq 1$, we set 
$$
T(f, \varphi , \Omega )=T(\chi _*f_*([S^n]), \varphi , \Omega ), 
$$
where $\chi $ is the morphism of simplicial sets defined in \S 5.3.
Then a homomorphism 
$$
\widehat{\varphi }(\Omega )_*:
\widehat{\pi }_n(\vert \widehat{G}(\ACYC )\vert , \CCGS )\to 
\widehat{\pi }_n(\vert \widehat{G}(Y)\vert , \CCGS )
$$
is defined by $\widehat{\varphi }(\Omega )_*([(f, \omega )])=
[(\varphi _*f, \varphi _!\omega -T(f, \varphi , \Omega ))]$.
The homomorphism 
$$
\widehat{\pi }_n(\vert \widehat{G}(X)\vert , \CCGS )\simeq 
\widehat{\pi }_n(\vert \widehat{G}(\ACYC )\vert , \CCGS )
\overset{\widehat{\varphi }(\Omega )_*}{\longrightarrow } 
\widehat{\pi }_n(\vert \widehat{G}(Y)\vert , \CCGS )
$$
is identified with the direct image morphism defined in Thm.6.9 
by the isomorphism $\widehat{\chi }_*$ in Prop.5.3.

\vskip 1pc
\subsection{The projection formula} \ \ 
In this subsection we prove the projection formula in higher 
arithmetic $K$-theory.
We first consider the case of $\widehat{\mathscr K}_0$-groups. 
Let $\varphi :X\to Y$, $h_{\varphi }$ and $\Omega $ be as in the 
last subsection. 
Let $\overline{E}$ be a hermitian vector bundle on $Y$ and 
$\overline{F}$ a $\varphi $-acyclic hermitian vector bundle on $X$.
Then the canonical isomorphism 
$\varphi _*(\varphi ^*\overline{E}\otimes \overline{F})\simeq 
\overline{E}\otimes \varphi _*\overline{F}$ preserves the metrics.

For $\omega \in \widetilde{\mathscr A}_1(Y)$ and 
$\tau \in \widetilde{\mathscr A}_1(X)$, we have 
{\allowdisplaybreaks
\begin{align*}
\varphi _!&(\widehat{\varphi }^*(\overline{E}, \omega )\times 
(\overline{F}, \tau ))   \\
&=\varphi _!((\varphi ^*\overline{E}, 
\varphi ^*\omega )\times (\overline{F}, \tau ))  \\
&=\varphi _!((\varphi ^*\overline{E})\otimes \overline{F}, 
\varphi ^*\omega \wedge \CCGS _0(\overline{F})+\varphi ^*
\CCGS _0(\overline{E})\wedge \tau +\varphi ^*d_{\mathscr A}\omega 
\wedge \tau )  \\
&=(\varphi _*(\varphi ^*\overline{E}\otimes \overline{F}), \eta ),
\end{align*}}
where
{\allowdisplaybreaks
\begin{align*}
\eta &=\int_{X(\mathbb C)/Y(\mathbb C)}\TD ^{GS}(\overline{T\varphi })
\wedge \left(\varphi ^*\omega \wedge \CCGS _0(\overline{F})+\varphi ^*
\CCGS _0(\overline{E})\wedge \tau +\varphi ^*d_{\mathscr A}\omega \wedge 
\tau \right)  \\*
&\hskip 1pc -T(\varphi ^*\overline{E}\otimes \overline{F})  \\
&=(\CCGS _0(\overline{E})+d_{\mathscr A}\omega )\wedge 
\int_{X(\mathbb C)/Y(\mathbb C)}\TD ^{GS}(\overline{T\varphi })
\wedge \tau \\*
&\hskip 1pc +\omega \wedge (\CCGS _0(\varphi _*\overline{F})
-d_{\mathscr A}T(\overline{F}))-T(\varphi ^*\overline{E}
\otimes \overline{F}). 
\end{align*}}
On the other hand, we have 
\begin{align*}
(\overline{E}, \omega )\times \varphi _!
(\overline{F}, \tau )&=(\overline{E}, \omega )\otimes 
(\varphi _*\overline{F}, \varphi _!\tau -T(\overline{F})) \\*
&=(\overline{E}\otimes \varphi _*\overline{F}, \eta ^{\prime }),
\end{align*}
where
\begin{align*}
\eta ^{\prime }&=(\CCGS _0(\overline{E})+d_{\mathscr A}\omega )\wedge 
\int_{X(\mathbb C)/Y(\mathbb C)}\TD ^{GS}(\overline{T\varphi })
\wedge \tau +\omega \wedge \CCGS _0(\varphi _*\overline{F})  \\
&\hskip 1pc -(\CCGS _0(\overline{E})+d_{\mathscr A}\omega )\wedge 
T(\overline{F}).
\end{align*}
Comparing these identities, we have 
$$
\eta -\eta ^{\prime }=-T(\varphi ^*\overline{E}\otimes \overline{F})
+\CCGS _0(\overline{E})\wedge T(\overline{F})+
d_{\mathscr A}(\omega \wedge T(\overline{F})).
$$
Hence the projection formula in $\widehat{\mathscr K}_0$-groups can be 
reduced to the following proposition:

\vskip 1pc
\begin{prop} \ \ 
Under the above notations, we have 
$$
T(\varphi ^*\overline{E}\otimes \overline{F})
=\CCGS _0(\overline{E})\wedge T(\overline{F}).
$$
\end{prop}

{\it Proof}: \  
Let $\mathscr E$ be an infinite dimensional vector bundle on $N$ 
consisting of smooth sections of 
$\Lambda ^*T^{* (1,0)}\varphi \otimes \overline{F}$.
Let $B_u$ and $N_u$ denote the Bismut superconnection and 
the number operator on $\mathscr E$ respectively. 
Let $\mathscr E^{\prime }$ be an infinite dimensional vector bundle 
consisting of smooth sections of $\Lambda ^*T^{* (1,0)}\varphi \otimes 
(\varphi ^*\overline{E}\otimes \overline{F})$ and the Bismut 
superconnection and the number operator on $\mathscr E^{\prime }$ are 
denoted by $B^{\prime }_u$ and $N^{\prime }_u$ respectively.
Then we have a canonical isometry 
$\mathscr E^{\prime }\simeq \overline{E}\otimes \mathscr E$ and 
under this identification, we have 
$B^{\prime }_u=1\otimes B_u+\nabla _{\overline{E}}\otimes 1$ and 
$N^{\prime }_u=1\otimes N_u$.
Hence we have $\exp (-{B^{\prime }_u}^2)=\exp(-K_{\overline{E}})
\wedge \exp (-B_u^2)$, therefore 
$$
\TTR _s(N^{\prime }_u\exp (-{B^{\prime }_u}^2))
=\CC _0(\overline{E})\wedge \TTR _s(N_u\exp (-B_u^2)).
$$
Substituting this identity into the definition of 
$T(\varphi ^*\overline{E}\otimes \overline{F})$ yields the desired 
identity.
\qed

\vskip 1pc
Let us move on to the higher case. 
We assume $n, m\geq 1$.
Consider the following diagram:
$$
\begin{CD}
\widehat{G}(Y)\wedge \widehat{G}(\ACYC )  
@>{m^G(\varphi ^*\wedge 1)}>> \widehat{G}^{(2)}(\ACYC ) \\
@VV{1\wedge \varphi _*}V @VV{\varphi _*}V  \\
\widehat{G}(Y)\wedge \widehat{G}(Y) @>{m^G}>> \widehat{G}^{(2)}(Y).
\end{CD}
$$
This diagram is commutative up to homotopy.
Let $f:S^n\to \vert \widehat{G}(Y)\vert $ and 
$g:S^m\to \vert \widehat{G}(\ACYC )\vert $ be pointed cellular maps.
For $\omega \in \widetilde{\mathscr A}_{n+1}(Y)$ and 
$\tau \in \widetilde{\mathscr A}_{m+1}(X)$, we have 
{\allowdisplaybreaks
\begin{align*}
\varphi (\Omega )_*(\varphi ^*(f, \omega )\times (g, \tau ))&=
\varphi (\Omega )_*(\varphi ^*f\times g, (-1)^n
\varphi ^*\CCGS _n(f)\bullet \tau +\varphi ^*\omega \bullet \CCGS _m(g) \\*
&\hskip 5pc +(-1)^nd_{\mathscr A}\varphi ^*
\omega \bullet \tau +(-1)^n\varphi ^*\CCGS _n(f)\vartriangle 
\CCGS _m(g))  \\
&=(\varphi _*(\varphi ^*f\times g), \eta ),
\end{align*}}
where
\begin{align*}
\eta &=(-1)^n(\CCGS _n(f)+d_{\mathscr A}\omega )\bullet 
\int_{X(\mathbb C)/Y(\mathbb C)}\TD ^{GS}(\overline{T\varphi })
\wedge \tau +\omega \bullet (\CCGS _m(\varphi _*g)-d_{\mathscr A}T(g))  \\*
&\hskip 2pc +(-1)^n\CCGS _n(f)\vartriangle (\CCGS _m(\varphi _*g)-
d_{\mathscr A}T(g))-T(\varphi ^*f\times g).
\end{align*}
On the other hand, we have 
\begin{align*}
(f, \omega )\times \varphi (h_{\varphi })_*(g, \tau )&=(f, \omega )
\times (\varphi _*g, \varphi _!\tau -T(g)) \\
&=(f\times \varphi _*g, \eta ^{\prime }),
\end{align*}
where
\begin{align*}
\eta ^{\prime }&=(-1)^n(\CCGS _n(f)+d_{\mathscr A}\omega )\bullet 
\int_{X(\mathbb C)/Y(\mathbb C)}\TD ^{GS}(\overline{T\varphi })
\wedge \tau -(-1)^n(\CCGS _n(f)+d_{\mathscr A}\omega )\bullet T(g)  \\
&\hskip 2pc +\omega \bullet \CCGS _n(\varphi _*g)+(-1)^n\CCGS _n(f)
\vartriangle \CCGS _m(\varphi _*g).
\end{align*}
A homotopy from $\varphi _*(\varphi ^*f\times g)$ to $f\times \varphi _*g$ 
is given by the canonical isometry $\varphi _*(\varphi ^*\overline{E}
\otimes \overline{F})\simeq \overline{E}\times \varphi _*\overline{F}$ 
of hermitian vector bundles.
Hence we can show the identity 
$$
[(\varphi _*(\varphi ^*f\times g), 0)]=[(f\times \varphi _*g, 0)]
$$ 
in $\widehat{\pi }_{n+m}(\vert \widehat{G}^{(2)}(Y)\vert , \CCGS )$ 
in the same way as Lem.5.6.
Furthermore, we have 
$$
\eta -\eta ^{\prime }\equiv (-1)^{n+1}\CCGS _n(f)\vartriangle d_{\mathscr A}
T(g)-T(\varphi ^*f\times g)+(-1)^n\CCGS _n(f)\bullet T(g)
$$
modulo $\IIm d_{\mathscr A}$. 
Hence the projection formula in higher arithmetic $K$-theory can be 
reduced to the following proposition:

\vskip 1pc
\begin{prop} \ \ 
For an exact metrized $n$-cube $\mathcal F$ on $Y$ and 
an exact metrized $m$-cube $\mathcal G$ made of $\varphi $-acyclic 
hermitian vector bundles on $X$, we have  
\begin{align*}
d_{\mathscr A}(\CCGS _n(\mathcal F)\vartriangle T(\mathcal G))= 
&-T(\varphi ^*\mathcal F\otimes \mathcal G)
+(-1)^n\CCGS _n(\mathcal F)\bullet T(\mathcal G)  \\
&+\CCGS _{n+1}(\partial \mathcal F)\vartriangle T(\mathcal G)
+(-1)^{n-1}\CCGS _n(\mathcal F)\vartriangle d_{\mathscr A}T(\mathcal G).
\end{align*}
\end{prop}

{\it Proof}: \ 
We will prove the following identities:
{\allowdisplaybreaks
\begin{align*}
d_{\mathscr A}(\CCGS _n(\mathcal F)\vartriangle T_1(\mathcal G))=
&-T_1(\varphi ^*\mathcal F\times \mathcal G)
+(-1)^n\CCGS _n(\mathcal F)\bullet T_1(\mathcal G)  \\*
&+\CCGS _{n-1}(\partial \mathcal F)\vartriangle T_1(\mathcal G)
+(-1)^{n+1}\CCGS _n(\mathcal F)\vartriangle d_{\mathscr A}T_1(\mathcal G), \\
d_{\mathscr A}(\CCGS _n(\mathcal F)\vartriangle T_2(\mathcal G))=
&-T_2(\varphi ^*\mathcal F\times \mathcal G)
+(-1)^n\CCGS _n(\mathcal F)\bullet T_2(\mathcal G)  \\*
&+\CCGS _{n-1}(\partial \mathcal F)\vartriangle T_2(\mathcal G)
+(-1)^{n+1}\CCGS _n(\mathcal F)\vartriangle d_{\mathscr A}T_2(\mathcal G).
\end{align*}}
These two identities can be proved in the same way, so we will prove 
only the latter one.

The following lemma is the renormalized version of Prop.4.5.

\vskip 1pc
\begin{lem} \ \ 
For $u_1, \cdots , u_n, v_1, \cdots , v_m\in \mathscr A_1(M)$, we have 
\begin{small}
\begin{align*}
d_{\mathscr A}&\left(\sum_{i=1}^n(-1)^iS_n^i(u_1, \cdots , u_n)\vartriangle 
\sum_{j=1}^m(-1)^jS_m^j(v_1, \cdots , v_m)\right)  \\
&=\sfrac{(-1)^{n+1}}{2\pi \sqrt{-1}}\sbinom{n+m}{n}^{-1}\sum_{k=1}^{n+m}
(-1)^kS_{n+m}^k(u_1, \cdots , u_n, v_1, \cdots , v_m)  \\
&\hskip 1pc -\sfrac{n}{4\pi \sqrt{-1}}
\underset{1\leq j\leq m}{\sum_{1\leq i\leq n-1}}
(-1)^{i+j}a^{n-1,m}_{i,j}\sum_{\alpha =1}^n(-1)^{\alpha }
\partial \overline{\partial }u_{\alpha }
S_{n-1}^i(u_1, \cdots , \widehat{u_{\alpha }}, \cdots , u_n)\wedge 
S_m^j(v_1, \cdots , v_m)  \\*
&\hskip 1pc +\sfrac{(-1)^nm}{4\pi \sqrt{-1}}
\underset{1\leq j\leq m-1}{\sum_{1\leq i\leq n}}
(-1)^{i+j}a^{n,m-1}_{i,j}S_n^i(u_1, \cdots , u_n)\wedge 
\sum_{\beta =1}^m(-1)^{\beta }\partial \overline{\partial }v_{\beta }
S_{m-1}^j(v_1, \cdots , \widehat{v_{\beta }}, \cdots , v_m)  \\
&\hskip 1pc +(-1)^n\underset{1\leq j\leq m}{\sum_{1\leq i\leq n}}
(-1)^{i+j}S_n^i(u_1, \cdots , u_n)\bullet S_m^j(v_1, \cdots , v_m).
\end{align*}
\end{small}
\end{lem}

\vskip 1pc
Let us start proving the latter identity.
For $t<s$, let $\pi _1:(\mathbb P^1)^s\to (\mathbb P^1)^t$ denote 
the projection defined by $\pi _1(x_1, \cdots , x_s)=(x_1, \cdots , x_t)$ 
and $\pi _2:(\mathbb P^1)^s\to (\mathbb P^1)^t$ denote the projection 
defined by $\pi _2(x_1, \cdots , x_s)=(x_{s-t+1}, \cdots , x_s)$.       
Then Lem.6.12 implies 
\begin{small} 
{\allowdisplaybreaks
\begin{align*}
&d_{\mathscr A}(\CCGS _n(\mathcal F)\vartriangle T_2(\mathcal G)) \\*
&=\sfrac{(-1)^{n+m+1}}{(2\pi \sqrt{-1})^{n+m-1}n!(m+1)!}
\int_{(\mathbb P^1)^{n+m}}\pi _1^*\CCGS _0(\TTR _n\lambda \mathcal F)
d_{\mathscr A}\left(\underset{1\leq j\leq m+1}
{\underset{1\leq i\leq n}{\sum }}(-1)^{i+j}\pi _1^*S_n^i\vartriangle 
\pi _2^*S_{m+1}^j(\mathcal G)\right)  \\
&=\sfrac{(-1)^{n+m}}{(2\pi \sqrt{-1})^{n+m}(n+m+1)!}
\int_{(\mathbb P^1)^{n+m}}\pi _1^*\CCGS _0(\TTR _n\lambda \mathcal F)
\sum _{k=1}^{n+m+1}(-1)^kS_{n+m+1}^k(T(\TTR _m\lambda \mathcal G), 
\log |t_1|^2, \cdots , \log |t_{n+m}|^2)  \\*
&\hskip 1pc +\sfrac{(-1)^{n+m}}{2(2\pi \sqrt{-1})^{n+m-1}(n-1)!(m+1)!}
\int_{(\mathbb P^1)^{n+m-1}}\pi _1^*\CCGS _0(\TTR _{n-1}\lambda \partial 
\mathcal F)\underset{1\leq j\leq m+1}{\underset{1\leq i\leq n-1}{\sum }}
(-1)^{i+j}a^{n-1,m+1}_{i,j}\pi _1^*S_{n-1}^i\wedge 
\pi _2^*S_{m+1}^j(\mathcal G) \\*
&\hskip 1pc +\sfrac{(-1)^{m+1}}{2(2\pi \sqrt{-1})^{n+m-1}n!m!}
\int_{(\mathbb P^1)^{n+m}}\pi _1^*\CCGS _0(\TTR _n\lambda \mathcal F)
\underset{1\leq j\leq m}{\underset{1\leq i\leq n}{\sum }}
(-1)^{i+j}a^{n,m}_{i,j}\pi _1^*S_n^i\wedge 
\pi _2^*\left(d_{\mathscr A}T(\TTR _m\lambda \mathcal G)S_m^j\right) \\*
&\hskip 1pc +\sfrac{(-1)^m}{2(2\pi \sqrt{-1})^{n+m-1}n!m!}
\int_{(\mathbb P^1)^{n+m-1}}\pi _1^*\CCGS _0(\TTR _n\lambda \mathcal F)
\underset{1\leq j\leq m}{\underset{1\leq i\leq n}{\sum }}
(-1)^{i+j}a^{n,m}_{i,j}\pi _1^*S_n^i\wedge 
\pi _2^*S_m^j(\partial \mathcal G) \\*
&\hskip 1pc +\sfrac{(-1)^{m+1}}{(2\pi \sqrt{-1})^{n+m-1}n!(m+1)!}
\int_{(\mathbb P^1)^{n+m}}\pi _1^*\CCGS _0(\TTR _n\lambda \mathcal F)
\underset{1\leq j\leq m+1}{\underset{1\leq i\leq n}{\sum }}
(-1)^{i+j}\pi _1^*S_n^i\bullet \pi _2^*S_{m+1}^j(\mathcal G).
\end{align*}}
\end{small}
By Prop.6.10, we have 
{\allowdisplaybreaks
\begin{align*}
\pi _1^*&\CCGS _0(\TTR _n\lambda \mathcal F)S_{n+m+1}^k
(T(\TTR _m\lambda \mathcal G), \log |t_1|^2, \cdots , \log |t_{n+m}|^2)  \\
&=S_{n+m+1}^k(\pi _1^*\CCGS _0(\TTR _n\lambda \mathcal F)\wedge \pi _2^*
T(\TTR _m\lambda \mathcal G), \log |t_1|^2, \cdots , \log |t_{n+m}|^2) \\
&=S_{n+m+1}^k(T(\TTR _{n+m}\lambda (\varphi ^*\mathcal F\otimes \mathcal G), 
\log |t_1|^2, \cdots , \log |t_{n+m}|^2) \\
&=S_{n+m+1}^k(\varphi ^*\mathcal F\otimes \mathcal G).
\end{align*}}
Moreover, we have 
{\allowdisplaybreaks
\begin{align*}
d_{\mathscr A}T_2(\mathcal G)&=\sfrac{(-1)^m}{(2\pi \sqrt{-1})^m(m+1)!}
\int_{(\mathbb P^1)^m}\sum_{j=1}^{m+1}(-1)^jdS_{m+1}^j(\mathcal G)  \\
&=\sfrac{(-1)^m}{(2\pi \sqrt{-1})^{m-1}m!}
\int_{(\mathbb P^1)^m}\sum_{j=1}^m(-1)^j\left(d_{\mathscr A}
T(\TTR _m\lambda \mathcal G)S_m^j-S_m^j(\partial \mathcal G)\right).
\end{align*}}
Hence we have 
\begin{small}
{\allowdisplaybreaks
\begin{align*}
&d_{\mathscr A}(\CCGS _n(\mathcal F)\vartriangle T_2(\mathcal G))  \\
&=\sfrac{(-1)^{n+m}}{(2\pi \sqrt{-1})^{n+m}(n+m+1)!}
\int_{(\mathbb P^1)^{n+m}}\sum _{k=1}^{n+m+1}(-1)^k
S_{n+m+1}^k(\varphi ^*\mathcal F\otimes \mathcal G)  \\*
&\hskip 1pc +\sfrac{(-1)^{n+m}}{2(2\pi \sqrt{-1})^{n+m-1}(n-1)!(m+1)!}
\int_{(\mathbb P^1)^{n+m-1}}\underset{1\leq j\leq m+1}
{\underset{1\leq i\leq n-1}{\sum }}(-1)^{i+j}a^{n-1,m+1}_{i,j}
\pi _1^*(\CCGS _0(\TTR _{n-1}\lambda \partial \mathcal F)S_{n-1}^i)
\wedge \pi _2^*S_{m+1}^j(\mathcal G) \\*
&\hskip 1pc +\sfrac{(-1)^{m+1}}{2(2\pi \sqrt{-1})^{n+m-1}n!m!}
\int_{(\mathbb P^1)^{n+m}}\underset{1\leq j\leq m}
{\underset{1\leq i\leq n}{\sum }}(-1)^{i+j}a^{n,m}_{i,j}
\pi _1^*\CCGS _0(\TTR _n\lambda \mathcal F)\pi _1^*S_n^i\wedge 
\pi _2^*\left(d_{\mathscr A}T(\TTR _m\lambda \mathcal G)S_m^j-
S_m^j(\partial \mathcal G)\right) \\*
&\hskip 1pc +\sfrac{(-1)^{m+1}}{(2\pi \sqrt{-1})^{n+m-1}n!(m+1)!}
\int_{(\mathbb P^1)^{n+m}}\underset{1\leq j\leq m+1}
{\underset{1\leq i\leq n}{\sum }}(-1)^{i+j}\pi _1^*
(\CCGS _0(\TTR _n\lambda \mathcal F)S_n^i)\bullet 
\pi _2^*S_{m+1}^j(\mathcal G)   \\
&=-T_2(\varphi ^*\mathcal F\otimes \mathcal G)
+\CCGS _{n-1}(\partial \mathcal F)\vartriangle T_2(\mathcal G)  \\*
&\hskip 1pc +(-1)^{n+1}\CCGS _n(\mathcal F)\vartriangle 
d_{\mathscr A}T_2(\mathcal G)+
(-1)^n\CCGS _n(\mathcal F)\bullet T_2(\mathcal G),
\end{align*}}
\end{small}
which completes the proof.
\qed
 
\vskip 1pc
Let us consider the case of $n=0$ and $m>0$.
Let $(\overline{E}, \omega )$ be a pair of a hermitian vector 
bundle on $Y$ and $\omega \in \widetilde{\mathscr A}_1(Y)$ and let 
$(g, \tau)$ be a pair of a pointed cellular map 
$g:S^m\to \vert \widehat{G}(X)\vert $ and 
$\tau \in \widetilde{\mathscr A}_{m+1}(X)$.
Then we have 
{\allowdisplaybreaks
\begin{align*}
\widehat{\varphi }(\Omega )_*&(\widehat{\varphi }^*(\overline{E}, \omega )
\times (g, \tau)) \\
&=(\varphi _*(\varphi ^*\overline{E}\otimes g), \varphi _!(\varphi ^*\omega 
\bullet (\CCGS _m(g)+d_{\mathscr A}\tau ))+
\varphi _!(\varphi ^*\CCGS _0(\overline{E})\bullet \tau )
-T(\varphi ^*\overline{E}\otimes g)   \\
&=(\varphi _*(\varphi ^*\overline{E}\otimes g), \omega \bullet 
\varphi _!(\CCGS _m(g)+d_{\mathscr A}\tau ))+\CCGS _0(\overline{E})
\bullet (\varphi _!\tau -T(g)).
\end{align*}
On the other hand, we have 
{\allowdisplaybreaks
\begin{align*}
(\overline{E}, \omega )&\times \widehat{\varphi }(\Omega )_*(g, \tau)    \\
&=(\overline{E}\otimes \varphi _*g, \CCGS _0(\overline{E})\bullet 
(\varphi _!\tau-T(g))+\omega \bullet \CCGS _m(\varphi _*g)+\omega \bullet 
(\varphi _!\tau-T(g))   \\
&=(\overline{E}\otimes \varphi _*g, \CCGS _0(\overline{E})\bullet 
(\varphi _!\tau-T(g))+\omega \bullet d_{\mathscr A}\varphi _!\tau
+\omega \bullet 
\varphi _!\CCGS _m(g)).   
\end{align*}
Hence we have   
$$
\widehat{\varphi }(\Omega )_*(\widehat{\varphi }^*(\overline{E}, \omega )
\times (g, \tau))=(\overline{E}, \omega )\times 
\widehat{\varphi }(\Omega )_*(g, \tau).
$$

In the case of $n>0$ and $m=0$, we can also prove the projection 
formula for the pairing $\widehat{K}_n\times \widehat{\mathscr K}_0
\to \widehat{K}_n$ in the same way.
Hence we have the following theorem:

\vskip 1pc
\begin{thm} \ \ 
Let $\varphi :X\to Y$ be a projective smooth morphism of proper 
arithmetic varieties.
Let $h_{\varphi }$ be an $F_{\infty }$-invariant hermitian metric 
on the relative tangent bundle $T\varphi (\mathbb C)$ that induces 
a K\"{a}hler metric on any fiber of $\varphi (\mathbb C)$.
Then for $y\in \widehat{K}_n(Y)$ and $x\in \widehat{K}_m(X)$, 
we have 
$$
\widehat{\varphi }(h_{\varphi })_*(\widehat{\varphi }^*y\times x)=y\times 
\widehat{\varphi }(h_{\varphi })_*(x).
$$
\end{thm}

\appendix
\vskip 2pc
\section{Some identities satisfied by binomial coefficients}
\vskip 1pc

\vskip 1pc
\begin{lem} \ \ 
{\rm (1)} \ For $0\leq k\leq i$, we have 
\begin{align*} 
(n-i)&\sum_{\alpha =0}^{k-1}\sbinom{n+m-i-j+1}{n-\alpha }
\sbinom{i+j-1}{\alpha }+i\sum_{\alpha =0}^k\sbinom{n+m-i-j}{n-\alpha }
\sbinom{i+j}{\alpha }  \\*
&=(n+m)\sum_{\alpha =0}^{k-1}\sbinom{n+m-i-j}{n-1-\alpha }
\sbinom{i+j-1}{\alpha }+(i-k)\sbinom{n+m-i-j}{n-k}\sbinom{i+j-1}{k}.
\end{align*}
In particular, we have 
$$
(n-i)\sum_{\alpha =0}^{i-1}\sbinom{n+m-i-j+1}{n-\alpha }
\sbinom{i+j-1}{\alpha }+i\sum_{\alpha =0}^i\sbinom{n+m-i-j}{n-\alpha }
\sbinom{i+j}{\alpha }=(n+m)\sum_{\alpha =0}^{i-1}
\sbinom{n+m-i-j}{n-1-\alpha }\sbinom{i+j-1}{\alpha }.
$$

\noindent
{\rm (2)} \ For $0\leq k\leq i$, we have 
\begin{align*}
(m-j)&\sum_{\alpha =k}^{i-1}\sbinom{n+m-i-j+1}{n-\alpha }
\sbinom{i+j-1}{\alpha }+j\sum_{\alpha =k}^{i-1}\sbinom{n+m-i-j}{n-\alpha }
\sbinom{i+j}{\alpha }  \\*
&=(n+m)\sum_{\alpha =k}^{i-1}\sbinom{n+m-i-j}{n-\alpha }
\sbinom{i+j-1}{\alpha }-(i-k)\sbinom{n+m-i-j}{n-k}\sbinom{i+j-1}{k-1}.
\end{align*}
In particular, we have 
$$
(m-j)\sum_{\alpha =0}^{i-1}\sbinom{n+m-i-j+1}{n-\alpha }
\sbinom{i+j-1}{\alpha }+j\sum_{\alpha =0}^{i-1}\sbinom{n+m-i-j}{n-\alpha }
\sbinom{i+j}{\alpha }=(n+m)\sum_{\alpha =0}^{i-1}
\sbinom{n+m-i-j}{n-\alpha }\sbinom{i+j-1}{\alpha }.
$$
\end{lem}

{\it Proof}: \ 
We will prove them by induction on $k$.
When $k=0$, the claim (1) is trivial. 
If the claim (1) holds for $k-1$, then we have 
{\allowdisplaybreaks
\begin{align*}
&(n-i)\sum_{\alpha =0}^{k-1}\sbinom{n+m-i-j+1}{n-\alpha }
\sbinom{i+j-1}{\alpha }+i\sum_{\alpha =0}^k
\sbinom{n+m-i-j}{n-\alpha }\sbinom{i+j}{\alpha }  \\
&=(n+m)\sum_{\alpha =0}^{k-2}\sbinom{n+m-i-j}{n-1-\alpha }
\sbinom{i+j-1}{\alpha }+(i-k+1)\sbinom{n+m-i-j}{n-k+1}\sbinom{i+j-1}{k-1} \\*
&\hskip 1pc +(n-i)\sbinom{n+m-i-j+1}{n-k+1}\sbinom{i+j-1}{k-1}
+i\sbinom{n+m-i-j}{n-k}\sbinom{i+j}{k}  \\
&=(n+m)\sum_{\alpha =0}^{k-2}\sbinom{n+m-i-j}{n-1-\alpha }
\sbinom{i+j-1}{\alpha }+(n-k+1)\sbinom{n+m-i-j}{n-k+1}\sbinom{i+j-1}{k-1} \\*
&\hskip 1pc +n\sbinom{n+m-i-j}{n-k}\sbinom{i+j-1}{k-1}
+i\sbinom{n+m-i-j}{n-k}\sbinom{i+j-1}{k}  \\
&=(n+m)\sum_{\alpha =0}^{k-1}\sbinom{n+m-i-j}{n-1-\alpha }
\sbinom{i+j-1}{\alpha }-(i+j-k)\sbinom{n+m-i-j}{n-k}\sbinom{i+j-1}{k-1} \\*
&\hskip 1pc +i\sbinom{n+m-i-j}{n-k}\sbinom{i+j-1}{k}  \\
&=(n+m)\sum_{\alpha =0}^{k-1}\sbinom{n+m-i-j}{n-1-\alpha }
\sbinom{i+j-1}{\alpha }+(i-k)\sbinom{n+m-i-j}{n-k}\sbinom{i+j-1}{k}. 
\end{align*}}
Hence the claim (1) holds for $k$.

The claim (2) for $k=i$ is trivial. 
If (2) holds for $k+1$, then we have 
{\allowdisplaybreaks
\begin{align*}
&(m-j)\sum_{\alpha =k}^{i-1}\sbinom{n+m-i-j+1}{n-\alpha }
\sbinom{i+j-1}{\alpha }+j\sum_{\alpha =k}^{i-1}\sbinom{n+m-i-j}{n-\alpha }
\sbinom{i+j}{\alpha }  \\
&=(n+m)\sum_{\alpha =k+1}^{i-1}\sbinom{n+m-i-j}{n-\alpha }
\sbinom{i+j-1}{\alpha }-(i-k-1)\sbinom{n+m-i-j}{n-k-1}\sbinom{i+j-1}{k}  \\*
&\hskip 1pc +(m-j)\sbinom{n+m-i-j+1}{n-k}\sbinom{i+j-1}{k}
+j\sbinom{n+m-i-j}{n-k}\sbinom{i+j}{k}  \\
&=(n+m)\sum_{\alpha =k+1}^{i-1}\sbinom{n+m-i-j}{n-\alpha }
\sbinom{i+j-1}{\alpha }+(m-j-i+k+1)\sbinom{n+m-i-j}{n-k-1}\sbinom{i+j-1}{k}\\*
&\hskip 1pc +m\sbinom{n+m-i-j}{n-k}\sbinom{i+j-1}{k}
+j\sbinom{n+m-i-j}{n-k}\sbinom{i+j-1}{k-1}  \\
&=(n+m)\sum_{\alpha =k+1}^{i-1}\sbinom{n+m-i-j}{n-\alpha }
\sbinom{i+j-1}{\alpha }+(n+m-k)\sbinom{n+m-i-j}{n-k}\sbinom{i+j-1}{k}  \\*
&\hskip 1pc +j\sbinom{n+m-i-j}{n-k}\sbinom{i+j-1}{k-1}  \\
&=(n+m)\sum_{\alpha =k}^{i-1}\sbinom{n+m-i-j}{n-\alpha }
\sbinom{i+j-1}{\alpha }-(i-k)\sbinom{n+m-i-j}{n-k}\sbinom{i+j-1}{k-1},
\end{align*}}
hence the claim (2) holds for $k$.
\qed

\vskip 1pc
\begin{lem} \ \ 
We have 
$$
\sbinom{n+m+l}{n}^{-1}\sum_{i=0}^n\sbinom{n+m-i-j}{n-i}
\sbinom{i+j-1}{i}\sum_{\alpha =0}^{i+j-1}
\sbinom{n+m+l-i-j-k+1}{n+m-\alpha }\sbinom{i+j+k-1}{\alpha }
=\sum_{\alpha =0}^{j-1}\sbinom{m+l-j-k+1}{m-\alpha }\sbinom{j+k-1}{\alpha }.
$$
\end{lem}

{\it Proof}: \ 
Let $F_n$ denote the left hand side of the above.
Then we have 
{\allowdisplaybreaks
\begin{align*}
F_n&=\sbinom{n+m+l}{n}^{-1}\sum_{i=1}^n\left((1-\frac{i-1}{n})
\sbinom{n+m-i-j+1}{n-i+1}\sbinom{i+j-2}{i-1}\sum_{\alpha =0}^{i+j-2}
\sbinom{n+m+l-i-j-k+2}{n+m-\alpha }\sbinom{i+j+k-2}{\alpha }\right. \\*
&\hskip 1pc +\left.\frac{i}{n}\sbinom{n+m-i-j}{n-i}\sbinom{i+j-1}{i}
\sum_{\alpha =0}^{i+j-1}\sbinom{n+m+l-i-j-k+1}{n+m-\alpha }
\sbinom{i+j+k-1}{\alpha }\right)  \\
&=\frac{1}{n}\sbinom{n+m+l}{n}^{-1}\sum_{i=1}^n\sbinom{n+m-i-j}{n-i}
\sbinom{i+j-2}{i-1}\left((n+m-i-j+1)\sum_{\alpha =0}^{i+j-2}
\sbinom{n+m+l-i-j-k+2}{n+m-\alpha }\sbinom{i+j+k-2}{\alpha }\right. \\*
&\hskip 1pc \left.+(i+j-1)\sum_{\alpha =0}^{i+j-1}
\sbinom{n+m+l-i-j-k+1}{n+m-\alpha }\sbinom{i+j+k-1}{\alpha }\right).
\end{align*}
By Lem.A.1, we have 
\begin{align*} 
F_n&=\frac{n+m+l}{n}\sbinom{n+m+l}{n}^{-1}\sum_{i=1}^n\sbinom{n+m-i-j}{n-i}
\sbinom{i+j-2}{i-1}\sum_{\alpha =0}^{i+j-2}
\sbinom{n+m+l-i-j-k+1}{n+m-1-\alpha }\sbinom{i+j+k-1}{\alpha }  \\
&=\sbinom{n+m+l-1}{n-1}^{-1}\sum_{i=0}^{n-1}\sbinom{n+m-i-j-1}{n-i-1}
\sbinom{i+j-1}{i}\sum_{\alpha =0}^{i+j-1}\sbinom{n+m+l-i-j-k}{n+m-1-\alpha }
\sbinom{i+j+k}{\alpha }  \\
&=F_{n-1}
\end{align*}}
Hence $F_n=F_0$, that is, 
$$
F_n=\sum_{\alpha =0}^{j-1}\sbinom{m+l-j-k+1}{m-\alpha }
\sbinom{j+k-1}{\alpha },
$$
which completes the proof.
\qed

\vskip 1pc
\begin{lem} \ \ 
If $0\leq i\leq n$ and $1\leq j\leq m$, 
$$
\sum_{\alpha =0}^i\sbinom{n+m-\alpha -j}{n-\alpha }
\sbinom{\alpha +j-1}{\alpha }=\sum_{\alpha =0}^i\sbinom{n+m-i-j}{n-\alpha }
\sbinom{i+j}{\alpha }.
$$
In particular, we have 
$$
\sum_{\alpha =0}^n\sbinom{n+m-\alpha -j}{n-\alpha }
\sbinom{\alpha +j-1}{\alpha }=\sum_{\alpha =0}^n\sbinom{m-j}{n-\alpha }
\sbinom{n+j}{\alpha }=\sbinom{n+m}{n}.
$$
\end{lem}

{\it Proof}: \ 
We prove the lemma by induction on $i$. 
When $i=0$, the statement of the lemma is clear. 
If the identity holds for $i-1$, then we have 
{\allowdisplaybreaks
\begin{align*}
\sum_{\alpha =0}^i\sbinom{n+m-\alpha -j}{n-\alpha }
\sbinom{\alpha +j-1}{\alpha }&=\sum_{\alpha =0}^{i-1}
\sbinom{n+m-i-j+1}{n-\alpha }\sbinom{i+j-1}{\alpha }
+\sbinom{n+m-i-j}{n-i}\sbinom{i+j-1}{i}  \\
&=\sum_{\alpha =0}^{i-1}\sbinom{n+m-i-j}{n-\alpha }\sbinom{i+j-1}{\alpha }
+\sum_{\alpha =0}^{i-1}\sbinom{n+m-i-j}{n-1-\alpha }\sbinom{i+j-1}{\alpha }
+\sbinom{n+m-i-j}{n-i}\sbinom{i+j-1}{i}  \\
&=\sum_{\alpha =0}^i\sbinom{n+m-i-j}{n-\alpha }\sbinom{i+j-1}{\alpha }
+\sum_{\alpha =1}^i\sbinom{n+m-i-j}{n-\alpha }\sbinom{i+j-1}{\alpha -1} \\
&=\sum_{\alpha =0}^i\sbinom{n+m-i-j}{n-\alpha }\sbinom{i+j}{\alpha }, 
\end{align*}}
which completes the proof.
\qed

\end{document}